\documentclass{article}
\usepackage{amsmath,amsthm,mathtools}
\usepackage{dsfont}
\usepackage{fullpage}
\usepackage[utf8]{inputenc}     
\usepackage[T1]{fontenc}
\usepackage[charter]{mathdesign}
\usepackage{graphicx}
\usepackage[english]{babel}
\usepackage{float}
\usepackage{pifont}
\usepackage{subcaption}
\usepackage{tikz}
\usetikzlibrary{patterns}
\usepackage{cleveref}

\newtheorem{thm}{Theorem}[section]
\newtheorem{prop}[thm]{Proposition}
\newtheorem{lem}[thm]{Lemma}

\newtheorem{claim}[thm]{Claim}
\newtheorem{hyp}[thm]{Hypothesis}
\newtheorem{remark}[thm]{Remark}

\crefname{prop}{Proposition}{Propositions}
\crefname{lem}{Lemma}{Lemmas}
\crefname{cor}{Corollary}{Corollaries}
\crefname{hyp}{Hypothesis}{Hypotheses}

\def\eps{{\varepsilon}}

\def\D{{\mathcal{D}}}
\newcommand{\opd}[1]{\D[\,{#1}\,]}
\newcommand{\opds}[1]{\D_s[\,{#1}\,]}
\newcommand{\opdns}[1]{\D_{ns}[\,{#1}\,]}
\newcommand{\opdr}[1]{\D_{J_{R}}[\,{#1}\,]}

\renewcommand{\epsilon}{\varepsilon}
\newcommand{\undu}{\underline{u}}

\author{Emeric Bouin \footnote{CEREMADE, UMR CNRS 7534, Université Paris-Dauphine, Université PSL, Place du Maréchal de Lattre de Tassigny, 75775 Paris cedex 16 (France), \texttt{bouin@ceremade.dauphine.fr}}
\and Jérôme Coville \footnote{UR 546 Biostatistique et Processus Spatiaux, INRAE - Centre de Recherche PACA, 228 route de l'Aérodrome, CS 40509, Domaine Saint Paul - Site Agroparc, 84914 Avignon cedex 9 (France), \texttt{jerome.coville@inrae.fr}, current address: Institut Camille Jordan -- UMR 5208, Université Claude Bernard Lyon 1, CNRS, Centrale Lyon, INSA Lyon, Université Jean Monnet, 69622 Villeurbanne (France)}
\and Guillaume Legendre \footnote{CEREMADE, UMR CNRS 7534, Université Paris-Dauphine, Université PSL, Place du Maréchal de Lattre de Tassigny, 75775 Paris cedex 16 (France), \texttt{guillaume.legendre@ceremade.dauphine.fr}}}

\begin{document}
\title{Sharp exponent of acceleration in general non-local equations with a weak Allee effect}
\maketitle




\begin{abstract}
We study an acceleration phenomenon arising in monostable inte\-gro-differential equations with a weak Allee effect. Previous works have shown its occurrence and have given correct upper bounds on the rate of expansion in some particular cases, but precise lower bounds were still missing. In this paper, we provide a sharp lower bound for this acceleration rate, which is valid for a large class of dispersion operators. Our results manage to cover fractional Laplace operators and standard convolutions in a unified way, which is new in the literature. An important result of the paper is a general flattening estimate of independent interest: this phenomenon appears regularly in acceleration situations, but getting quantitative estimates is in general an open question. With this estimate at hand, we construct a subsolution that captures the expected behaviour of the accelerating solution (rates of expansion and flattening) and identifies several regimes that appear in the dynamics depending on the parameters of the problem.
\end{abstract}

\maketitle

\section{Introduction}\label{sec:Intro}
In this paper, we are interested in describing quantitatively the propagation phenomenon in the following (non-local) integro-differential equation, complemented with an initial condition:
\begin{align}
&\partial_tu(t,x) =\opd{u}(t,x) +f(u(t,x)) \quad \text{ for all } t>0\text{ and }x\text{ in }\mathbb{R},\label{eq:main-gen}\\
&u(0,x)= u_0(x)\quad \text{ for all } x\text{ in }\mathbb{R},\label{bc:main-gen}
\end{align}
where the function $u$ represents a density of population and thus takes its values in $[0,1]$, the function $f$ is a nonlinearity to be specified, accounting for a reaction process, the nonnegative function $u_0$ is the initial density, and the dispersal operator $\opd{\cdot}$ is defined by
\[
\opd{u}(t,x):=\mathrm{P.V.}\left(\int_{\mathbb{R}}[u(t,y)-u(t,x)] J(x-y)\,dy\right),
\]
where $\mathrm{P.V.}$ stands for principal value, the dispersal kernel $J$ being a nonnegative function satisfying the following properties.

\begin{hyp}\label{hyp:J}
Let $s$ be a positive real number. The kernel $J$ is nonnegative, symmetric and such that there exist positive constants $\mathcal{J}_0$, $\mathcal{J}_1$ and $R_0\ge 1$ verifying
\[
\int_{-1}^1J(z)z^2\,dz\le 2\mathcal{J}_1 \text{ and } \frac{{\mathcal{J}_0}^{-1}}{|z|^{1+2s}} \mathds{1}_{\{|z|\ge R_0\}}(z) \leq J(z) \mathds{1}_{\{|z|\ge 1\}}(z) \leq \frac{\mathcal{J}_0}{|z|^{1+2s}}.
\]
\end{hyp}

The operator $\opd{\cdot}$ describes the dispersion process of individuals. Rough\-ly speaking, the value $J(x-y) dy$ is the probability of a jump from position $x$ to position $y$, which makes the tails of the dispersal kernel $J$ of crucial importance when quantifying the dynamics of the population. As a matter of fact, the parameter $s$ will appear in the rates we obtain. One may readily notice that the hypothesis on $J$ allows us to cover both main types of integro-differential operators usually considered in the literature: the fractional Laplace operator $(-\Delta)^s$ on the one hand, and a standard convolution operator with an integrable kernel, often written $J \star u - u$, on the other. This universality is one of the main contributions of the present paper. 

\smallskip

In what follows, we assume that the nonlinearity satisfies the following hypothesis.
\begin{hyp}\label{hyp:f}
The function $f$ belongs to $\mathscr{C}^{1}([0,1],\mathbb{R})$ and is of the monostable type, in the sense that
\[
\begin{cases}
f(0)=f(1)=0,\\
f(z)>0\text{ for all }z\text{ in }(0,1)\\
f'(1)<0,\\
r \leq \displaystyle{\lim_{z\to 0^+} \frac{f(z)}{z^{\beta}}\leq R},
\end{cases}
\]
for some real numbers $0 < r \leq R $ and $\beta > 1$.
\end{hyp}

The parameter $\beta$ in the nonlinearity above quantifies the possibility of a weak \textit{Allee effect} that the population has to overcome. A biological description and discussion about the origin and relevance of such a phenomenon can be found in the book by Courchamp \textit{et al.} \cite{Courchamp2008}, and also in \cite{Allee1938,Berec2007,Dennis1989}. In rough terms, the Allee effect states that the growth rate per capita of the population is positively correlated to the population density, meaning that a population with too few individuals will not be fit enough to persist and grow. It is said to be \emph{weak} whenever the growth rate of a very small population is eventually extremely small, but still positive, as opposed to a \emph{strong} effect, leading to negative growth rates for small populations. In the sequel, we take $\beta >1$, which yields small growth rates for small densities, and assume that the initial datum satisfies the following hypothesis.
\begin{hyp}\label{hyp:u0}
The initial datum $u_0$ belongs to $\mathscr{C}_{bu}(\mathbb{R})$, the space of bounded and uniformly continuous functions on $\mathbb{R}$, and is such that\\ $a \mathds{1}_{(-\infty,b]}\le u_0\le \mathds{1}_{(-\infty,c]}$ for some real numbers $0<a<1$ and~$b<c$.
\end{hyp}

Under Hypotheses \ref{hyp:J}, \ref{hyp:f} and \ref{hyp:u0}, the existence of a global \emph{mild} solution to the Cauchy problem \eqref{eq:main-gen}-\eqref{bc:main-gen}, that is, a weak solution $u$ in $\mathscr{C}([0,+\infty),\mathscr{C}_{bu}(\mathbb{R}))$, is actually a non-trivial question in general. A particular situation, however, is the case of an integrable dispersal kernel, for which the existence of a strong solution in $\mathscr{C}^{1}([0,+\infty),\mathscr{C}_{bu}(\mathbb{R}))$ is a consequence of the Cauchy--Lipschitz theorem. When the kernel is not integrable, one can construct a mild solution splitting the dispersal operator as follows: $-\opd{\cdot}=-\opds{\cdot}-\opdns{\cdot}$, where $\opds{\cdot}$ and $\opdns{\cdot}$ are associated to the kernels $J\mathds{1}_{\{|z|<1\}}$ and $J\mathds{1}_{\{|z|\ge 1\}}$, respectively.
Assuming that $-\opds{\cdot}$ is the infinitesimal generator of a strongly continuous semigroup on $\mathscr{C}_{bu}(\mathbb{R})$, the Cauchy problem \eqref{eq:main-gen}-\eqref{bc:main-gen} possesses a global \emph{mild} solution explicitly given by the Duhamel principle,
\[
u(t,x)= (T^s_tu_0)(x) +\int_{0}^tT^s_{t-\tau}\Big(\opdns{u}(\tau,x)+f(u(\tau,x))\Big)\,d\tau.
\]
 where $T_t^s$ denotes the semigroup generated by the operator $-\opds{\cdot}$. Conditions for which $-\opds{\cdot}$ is the infinitesimal generator of a strongly continuous semigroup on $\mathscr{C}_{bu}(\mathbb{R})$ can be derived from the properties of the associated heat kernel, which is a classical problem in the literature of L\'evy processes, explored, for example, in  \cite{Kolokoltsov2000,Chen2008a,Barlow2009,Chen2011,Bogdan2020,Knopova2021}.

To avoid discussing supplementary assumptions on $J$ and/or a dichotomy between integrable and non-integra\-ble kernels, we will also assume the following.  
\begin{hyp}\label{hyp:D}
The kernel $J$ is such that $-\opds{\cdot}$ is the infinitesimal generator of a strongly continuous semigroup on $\mathscr{C}_{bu}(\mathbb{R})$.
\end{hyp}
Note that such an assumption covers a large class of dispersal operators, ranging from convolutions to the fractional Laplacian. It is worth noticing that, for non-integrable kernels, a sufficient condition is given in \cite{Chen2011} and essentially corresponds to having the following behaviour near the origin:   
$$\frac{C^{-1}}{|z|^{1+2s_1}}\mathds{1}_{\{\{|z|<1\}}(z)\le  J(z)\mathds{1}_{\{\{|z|<1\}}(z)\le  \frac{C}{|z|^{1+2s_2}}\mathds{1}_{\{\{|z|<1\}}(z), $$ with $C>0$ and 
$0<s_1\le s_2<1$.

In the rest of the paper, if not explicitly mentioned, a solution to the Cauchy problem \eqref{eq:main-gen}-\eqref{bc:main-gen} will always refer to a global mild solution in $\mathscr{C}([0,+\infty),\mathscr{C}_{bu}(\mathbb{R}))$.

To follow the propagation of the population modelled by system \eqref{eq:main-gen}-\eqref{bc:main-gen}, we define some particular position of the level set of height $\lambda$, with $\lambda$ a real number in $(0,1)$, of a solution $u$ to the problem, that is
\[
x_\lambda(t) := \inf\left\{ x \in \mathbb{R}, \, u(t,x) \leq \lambda \right\}\text{ for all }t>0.
\]
When the initial data $u_0$ is nonincreasing, this point can be alternatively defined using the following quantity 
\[
x_\lambda(t) := \sup\left\{ x \in \mathbb{R}, \, u(t,x) \geq \lambda \right\}\text{ for all }t>0.
\]

\subsection*{Existing works and previous results}
Let us review the existing literature in order to position our work. Propagation phenomena in reaction-diffusion and integro-differential equations have been the object of intense studies in the last decades. Starting from the work of Fisher on the propagation of an advantageous gene in a population \cite{Fisher1937} and its analysis by Kolmogorov, Petrovski and Piskunov \cite{Kolmogorov1937} and related works by \textit{e.g.} Aronson and Weinberger \cite{Aronson1978}, the quantitative description of spreading gave birth to various mathematical tools and techniques such as travelling waves, accelerating profiles, transition fronts, among many others. 

When $\beta = 1$ and the nonlinearity $f$ is such that $f(z)\le f'(0)z$, meaning it is a Fisher-KPP nonlinearity, it is known that solutions to problem \eqref{eq:main-gen}-\eqref{bc:main-gen}  {exhibit some propagation phenomenon}: starting with a nonnegative non-trivial compactly supported initial datum, the corresponding solution $u$ converges to $1$ locally uniformly in space as time tends to infinity. This is referred to as the {\it hair-trigger effect} \cite{Aronson1978}. Moreover, in many cases, this convergence can be precisely characterised. For instance, when the dispersal kernel is continuous and exponentially bounded, that is when there exists a positive real number $\mu$ such that
\begin{equation*}
\int_{\mathbb{R}} J(z) e^{\mu \vert z \vert} \, dz < +\infty,
\end{equation*}
travelling waves are known to exist and solutions to the Cauchy problem typically propagate at constant speed, see \cite{Schumacher1980,Weinberger1982,Carr2004,Coville2007a,Coville2008a,Lutscher2005,Yagisita2009}. On the other hand, when the kernel $J$ possesses heavy tails, in the sense that, for any positive real number $\eta$, one has
\begin{equation*}
\lim_{\vert z \vert \to+\infty} J(z) e^{\eta \vert z \vert} \, = +\infty,
\end{equation*}
travelling waves do not exist and the solutions exhibit an acceleration phenomenon, see \cite{Medlock2003,Yagisita2009,Garnier2011}. Garnier \cite{Garnier2011} gave the first acceleration estimates for convolution operators, showing that $x_\lambda$ spreads super-linearly in time and providing explicit upper and lower bounds of order $J^{-1}(e^{-\rho t})$ for some positive real number $\rho$. The first author with Garnier, Henderson and Patout \cite{Bouin2018} next refined these bounds by providing sharp estimates for the level sets of the solution. In particular, when $J$ is integrable and has an algebraic decay at infinity, the level set of height $\lambda$ expands exponentially fast: $x_\lambda(t) \sim e^{\rho t}$ for an explicit value $\rho$. In the same period of time, a group around Cabré and Roquejoffre \cite{Cabre2009,Cabre2013,Meleard2015} studied the fractional Fisher-KPP equation, concluding to a similar exponential propagation behaviour, using various techniques. The same group later explored the acceleration phenomena occurring in some heterogeneous media, obtaining exponential propagation in periodic media and revealing an interplay between some properties of the semigroup of the linearised operator and the exponential growth generated by the nonlinearity \cite{Cabre2012,Coulon2012} (see also \cite{Meleard2015}). Note that a related, but different, acceleration phenomenon for positive solutions of a local Cauchy problem appears in reaction-diffusion equations of KPP-type when playing with the tails of the initial datum \cite{Hamel2010}. We emphasise that, in the present work, acceleration is solely due to the structure of the dispersal operator.

When an Allee effect is introduced, the study of propagation becomes more subtle.  Coville \textit{et al.} \cite{Coville2007a,Coville2008a,Coville2007d} proved the existence of travelling fronts when the dispersal kernel is exponentially bounded and the Cauchy problem typically does not lead to acceleration \cite{Zhang2012}. When not, the competition between the heavy tails of the kernel and the Allee effect leads to intense discussions. Gui and Huan \cite{Gui2015} first investigated the existence or not of travelling front solutions for equations with a fractional Laplace operator and a weak Allee effect. They obtained existence (and thus propagation at a fixed speed) when $s$ belongs to $\left(\frac{1}{2},1\right)$ and $\beta\ge 2$ satisfy the relation $\tfrac{\beta}{2s(\beta-1)}\le 1$, this being equivalent to $$(2s-1)(\beta-1)\ge 1. $$ They also showed that this condition is sharp, in the sense that no front solution can exist when $\tfrac{\beta}{2s(\beta-1)}>1$. However, in this context, neither a description of acceleration nor a precise rate of acceleration were given.

In the same spirit, for algebraically decaying integrable kernels with a finite first moment, Alfaro and Coville \cite{Alfaro2017} showed that the curve $ (2s-1)(\beta-1)=1 $ is the separatrix between existence and nonexistence of travelling waves for convolution type equations, that is front solutions exist if and only if $\tfrac{\beta}{2s(\beta-1)}\le 1$. In addition, they proved that this curve provides the exact separation between nonaccelerated and accelerated solutions to the Cauchy problem and gave the first upper and lower bounds on the expansion of the level sets, revealing a delicate interplay between the tails of the dispersal kernel and the nonlinearity parameter $\beta$. More precisely, in such a situation, solutions spread at least as $t^{\frac{1}{2s(\beta-1)}}$ and at most as $t^{\frac{\beta}{2s(\beta-1)}}$ when $\tfrac{\beta}{2s(\beta-1)} >1$. Let us mention that Alfaro \cite{Alfaro2017a} concurrently described the interaction between a heavy-tailed initial datum and the Allee effect in local reaction-diffusion equations. Alfaro and Giletti \cite{Alfaro2017b,Alfaro2019} later described a similar interplay in local porous media reaction-diffusion equations, in which some particular self-similar solutions of the equation play a role.  For a recent contribution regarding the Allee effect in diffusive equations with mixed type dispersal, see \cite{Dipierro2024}.

In a recent collaboration between Gui, Zhao and the second author \cite{Coville2021}, the constraint $s>\frac{1}{2}$ was lifted in the study of acceleration with a fractional Laplace diffusion. Nonmatching upper and lower bounds on the expansion of the level sets of solutions, similar to those given in \cite{Alfaro2017}, were obtained. More explicitly, it was shown that, when $s$ belongs to $(0,1)$ and $\tfrac{\beta}{2s(\beta-1)}>1$, solutions spread at most as $t^{\frac{\beta}{2s(\beta-1)}}$ and at least as $t^{\min\{\frac{1}{2s}, \frac{1}{2s(\beta-1)}\}}$, but the determination of the exact speed of the level sets was left as an open question. It was recently answered by the authors in an early version of the present work \cite{Bouin2021-hal} for a generic integral dispersal operator and, independently, by Zhang and Zlato\v{s} \cite{Zhang2023} for the fractional Laplace operator.

In order to provide a clear picture of the existing results in both the convolution and the fractional Laplacian cases, we have summarised in Figure \ref{fig:diagrams} the already known behaviours recalled above. It is worth mentioning that the derivation of the refined lower bound obtained for the fractional Laplace operator crucially involves a precise knowledge of the form of the associated heat kernel. However, such understanding is lacking for most of the other operators considered here, even when the kernel $J$ is integrable.

\begin{center}
\begin{figure}[H]\centering
\begin{subfigure}[b]{0.4\textwidth}
\centering\resizebox{!}{140pt}{\begin{tikzpicture}[scale=0.9]


\fill[pattern=north west lines,pattern color=black!50!green, opacity=0.5] 
plot [domain=3.429:10] (\x, {(2*(\x/6))/(2*(\x/6)-1)})
--(10,8)--(3.429,8)
-- cycle;

\fill[pattern=north west lines, pattern color=blue, opacity=0.9] 
(3,2)--(3,8)--(3.429,8)
-- plot [domain=3.429:10] (\x, {(2*(\x/6))/(2*(\x/6)-1)})
-- (10,1.3)
-- plot [domain=3:10] (\x, {1+1/(2*(\x/6))})
-- cycle;

\fill[pattern=north west lines, pattern color=violet, opacity=0.5] 
(3,1)--(3,2)
-- plot [domain=3:10] (\x, {1+1/(2*(\x/6))})
-- (10,1)
-- (6.6,1) -- (6.6,1.3) -- (4.4,1.3) -- (4.4,1) 
-- cycle;

\draw[->,line width=1.2, color=black](9.5,1)--(10.5,1) node[below]{$s$};
\draw[->,line width=1.2, color=black](0.5,1)--(0.5,8.5) node[right]{$\beta$};
\draw[line width=1, dashed](3,1)--(3,8.5) node[right]{$s=\frac12$};
\draw[line width=1, dashed](6,1.25)--(6,8.5) node[right]{$s=1$};

%
%
\draw [domain=3.429:10,line width=3, samples=150, black!50!green] plot (\x, {(2*(\x/6))/(2*(\x/6)-1)}) node[anchor=south east]{$\qquad \beta = \frac{2s}{2s-1}$};

\draw [domain=3:10,line width=1, samples=150, violet] plot (\x, {1+1/(2*(\x/6))}) ;

\draw [line width=3,orange] (6.6,1)--(10,1) node[anchor=north east]{$\beta = 1$};
\draw [line width=3,orange] (0.5,1)--(4.4,1) ;


\node at (5.5,1)[line width=6,orange]{$x_\lambda(t) \asymp_\lambda e^{\rho t}$};
\node at (1.8,5)[line width=6,black]{\ding{204}};
\node at (3.7,3)[line width=6,blue]{\ding{202}};
\node at (3.7,1.4)[line width=6,violet]{\ding{205}};
\node at (7.5,5.2)[line width=15,black!50!green]{$x_\lambda(t) \asymp t$};
\node at (5,5)[line width=6,black!50!green]{\ding{203}};
\end{tikzpicture}}
\caption{The convolution case}
\label{fig:diagramconvolee}
\end{subfigure}
\qquad 
\begin{subfigure}[b]{0.45\textwidth}
\centering\resizebox{!}{140pt}{\begin{tikzpicture}[scale=0.9]


\fill[pattern=north west lines,pattern color=black!50!green, opacity=0.5] 
plot [domain=3.429:6] (\x, {(2*(\x/6))/(2*(\x/6)-1)})
--(6,8)--(3.429,8)
-- cycle;

\fill[pattern=north west lines, pattern color=blue, opacity=0.5] 
(0.5,2) -- (0.5,8) --(3.429,8)
-- plot [domain=3.429:6] (\x, {(2*(\x/6))/(2*(\x/6)-1)})
-- (6,1.5)
-- plot [domain=6:3] (\x, {1+1/(2*(\x/6))})
-- cycle;

\fill[pattern=north west lines, pattern color=blue, opacity=0.5] 
(0.5,2)--(3,2)
-- plot [domain=3:6] (\x, {1+1/(2*(\x/6))})
-- (6,1)
-- (4.4,1) -- (4.4,1.3) -- (2,1.3) -- (2,1) 
-- (0.5,1)
-- cycle;
%
%
\draw[->,line width=1.2, color=black](5.5,1)--(6.25,1) node[below]{$s$};
\draw[->,line width=1.2, color=black](0.5,1)--(0.5,8.5) node[right]{$\beta$};
\draw[line width=1, dashed](3,1.5)--(3,8.5) node[right]{$s=\frac12$};
\draw[line width=1.5](6,1)--(6,8.5) node[right]{$s=1$};
%
%
\draw [domain=6:3.429,line width=3, samples=150, black!50!green] plot (\x, {(2*(\x/6))/(2*(\x/6)-1)}) node[anchor=north west]{$\beta = \frac{2s}{2s-1}$};

\draw [line width=3,orange] (4.4,1)--(6,1) node[anchor=north east]{$\beta = 1$};
\draw [line width=3,orange] (0.5,1)--(2,1) ;


\node at (3.2 ,1)[line width=6,orange]{$x_\lambda(t) \asymp_\lambda e^{\rho t}$};
\node at (2,5)[line width=6,blue]{\ding{202}};
\node at (4.8,6.2)[line width=15,black!50!green]{$x_\lambda(t)\asymp  t$};
\node at (5,5)[line width=6,black!50!green]{\ding{203}};
\end{tikzpicture}}
\caption{ The fractional case}
\label{fig:diagramfrac}
\end{subfigure}
\caption{
{\bf Subfigure (a): summary of the existing results  for the convolution case.} In the blue zone \textcolor{blue}{\ding{202}}, only an upper bound has been derived \cite{Alfaro2017}: $x_\lambda(t) \lesssim_\lambda t^{\frac{\beta}{2s(\beta-1)}}$. In the green zone \textcolor{black!50!green}{\ding{203}}, the model enjoys linear propagation with existence of travelling fronts \cite{Coville2006}: $x_\lambda(t) \asymp t$. In the white zone \ding{204}, no estimates are known. In the violet zone \textcolor{violet}{\ding{205}}, nonmatching lower and upper bounds have been obtained \cite{Alfaro2017}: $t^{\frac{1}{2s(\beta-1)} } \lesssim_\lambda x_{\lambda}(t) \lesssim_\lambda t^{\frac{\beta}{2s(\beta-1)} }$. In the orange zone, exponential propagation occurs \cite{Garnier2011,Bouin2018}: $x_\lambda(t)\asymp_\lambda\exp({\rho} t)$.\\
{\bf Subfigure (b): summary of the existing results for the fractional case.} In the blue zone~\textcolor{blue}{\ding{202}}, an upper bound was derived in \cite{Coville2021}, with a matching lower bound found by the authors in an early version of the present work \cite{Bouin2021-hal} and, independently, by Zhang and Zlato\v{s} in \cite{Zhang2023}: $x_\lambda(t) \asymp_\lambda t^{\frac{\beta}{2s(\beta-1)}}$. In the green zone \textcolor{black!50!green}{\ding{203}}, the model enjoys linear propagation with existence of travelling fronts \cite{Gui2015,Coville2021}: $x_\lambda(t) \asymp t$. In the orange zone, exponential propagation occurs \cite{Cabre2013}: $x_\lambda(t)\asymp_\lambda\exp( \rho t)$.}\label{fig:diagrams}
\end{figure}
\end{center}

One of the main purposes of the present article is to fill the blanks in the above diagrams by providing the exact speed of propagation of the level sets. To do so, an approach is proposed that extends the existing estimates to a generic situation in which the dispersal kernel satisfies the properties listed above. Our aim is to highlight the essential features of the system (related to both the kernel and nonlinearity) required to ensure an accurate description of the acceleration phenomenon in this context.

\subsection*{Statement of the main result}
Let us precisely state the main result of the present paper. 
\begin{thm}\label{thm:main}
Assume that the dispersal kernel $J$, the initial datum $u_0$, and the nonlinearity $f$ respectively satisfy Hypotheses \ref{hyp:J}, \ref{hyp:D}, \ref{hyp:u0}, and \ref{hyp:f}, and that the parameters $s > 0$ and $\beta > 1$ are such that
\begin{equation}\label{s-beta rel} 
(2s-1)(\beta-1)<1.
\end{equation}
Then, for any $\lambda$ in $(0,1)$, the level set $x_\lambda$ of a solution $u$ in $\mathscr{C}([0,+\infty),\mathscr{C}(\mathbb{R}))$ to problem \eqref{eq:main-gen}-\eqref{bc:main-gen} accelerates with a rate equal to $\frac{\beta}{2s(\beta-1)}$, that is\footnote{The notation $u \asymp_\lambda v$ means that there exists a constant $C_\lambda$ in $(0,1]$, depending  only on $\lambda$, such that $C_\lambda v \leq u\leq {C_\lambda}^{-1} v$. Similarly, the notation $u\lesssim_\lambda v$ means that there exists a positive constant $C_\lambda$, depending only on $\lambda$, such that $u\leq C_\lambda v$.}
\[
x_\lambda(t) \asymp_\lambda t^{\frac{\beta}{2s(\beta-1)}}.
\]
\end{thm} 
Note that condition \eqref{s-beta rel} above is trivially satisfied when $s \leq \frac12$ and gives birth to a separatrix when not. To give the reader a clear view of the scope of \Cref{thm:main}, we have summarised previous contributions with ours in \Cref{fig:diagram}. The algebraic rates derived in \Cref{thm:main} are in sharp contrast with the situation $\beta = 1$ for which exponential propagation is expected, see for instance \cite{Garnier2011,Cabre2013,Bouin2018}.

\begin{figure}[H]\centering
\resizebox{!}{160pt}{
\begin{tikzpicture}[scale=0.9]


\fill[pattern=north west lines,pattern color=black!50!green, opacity=0.5] 
plot [domain=3.429:10] (\x, {(2*(\x/6))/(2*(\x/6)-1)})
--(10,8)--(3.429,8)
-- cycle;

\fill[pattern=north west lines, pattern color=blue, opacity=0.5] 
(0.5,1)--(0.5,8)--(3.429,8)
-- plot [domain=3.429:10] (\x, {(2*(\x/6))/(2*(\x/6)-1)})
-- (10,1)
-- (6.6,1) -- (6.6,1.3) -- (4.4,1.3) -- (4.4,1) 
-- cycle;

%
%
\draw[->,line width=1.2, color=black](9.5,1)--(10.5,1) node[below]{$s$};
\draw[->,line width=1.2, color=black](0.5,1)--(0.5,8.5) node[right]{$\beta$};
\draw[line width=1, dashed](3,4.2)--(3,8.5) node[right]{$s=\frac12$};
\draw[line width=1, dashed](3,1)--(3,3) ;
%
%
\draw [domain=3.429:10,line width=3, samples=150, black!50!green] plot (\x, {(2*(\x/6))/(2*(\x/6)-1)}) node[anchor=south east]{$\qquad \beta = \frac{2s}{2s-1}$};
\draw [line width=3,orange] (6.6,1)--(10,1) node[anchor=north east]{$\beta = 1$};
\draw [line width=3,orange] (0.5,1)--(4.4,1) ;


\node at (5.5,1)[line width=6,orange]{$x_\lambda(t) \asymp_\lambda e^{\rho t}$};
\node at (2.25,3.5)[line width=6,blue]{$x_\lambda(t) \asymp_\lambda t^{\frac{\beta}{2s(\beta-1)}}$};
\node at (6.5,5.2)[line width=15,black!50!green]{$x_\lambda(t) \asymp t$};
\node at (2,5)[line width=6,blue]{\ding{202}};
\node at (5,5)[line width=6,black!50!green]{\ding{203}};
\end{tikzpicture}}
\caption{\small{In the blue zone \textcolor{blue}{\ding{202}}, sharp lower and upper bounds are provided by \Cref{thm:main}: $x_\lambda(t)\asymp_\lambda t^{\frac{\beta}{2s(\beta-1)}}$. In the green zone \textcolor{black!50!green}{\ding{203}}, the model enjoys a linear propagation with existence of travelling fronts \cite{Alfaro2017,Coville2006,Coville2021,Gui2015}: $x_\lambda(t)\asymp t$. In the orange zone, exponential propagation occurs, by straightforward extension of the work of Bouin \textit{et al.} \cite{Bouin2018}: $x_\lambda(t) \asymp_\lambda \exp(\rho t)$.}}\label{fig:diagram}
\end{figure}


To the best of our knowledge, \Cref{thm:main} provides the first sharp and unified estimate for the propagation of level sets in such a generic context. While correct upper bounds in some particular settings had previously been derived, no precise lower bound existed, except in the specific case of the fractional Laplacian \cite{Zhang2023}. One can observe that condition \eqref{s-beta rel} generalises the ones given in related papers \cite{Alfaro2017,Bouin2021,Coville2021,Gui2015,Zhang2023}. It should be also noted that we obtained the rate of invasion for a convolution operator when $s$ belongs to $(0,\frac12]$, a question which remained open in \cite{Alfaro2017}.

The constructions used in \cite{Alfaro2017,Coville2021} to obtain upper bounds are robust and can be adapted to the general setting considered in the present paper for kernels satisfying Hypothesis \ref{hyp:J}. For the sake of brevity, they will not be duplicated here. Our contribution is thus a generic way of obtaining lower bounds that match the already known upper bounds. 

As already mentioned, Zhang and Zlato\v{s} \cite{Zhang2023} obtained, in the particular case of the fractional Laplace operator in the presence of an ignition or monostable  nonlinearity, lower bounds similar to those appearing in an early version of the present work, submitted to arXiv in May 2021 \cite{Bouin2021} and completed in the current form in November 2021 \cite{Bouin2021-hal}. Although the construction of a subsolution in their work shares some common points with ours, they strongly rely on the homogeneity between small and large jumps of a L\'evy process associated to the fractional Laplacian, a specific feature which is not guaranteed to hold when a general dispersal kernel is considered. In \cite{Zhang2023}, this property is crucially used to construct Lipschitz continuous compactly supported subsolutions, making this approach difficult to extend to kernels whose nature changes drastically between small and large jumps. They also do not try to capture the proper behaviour of the solution at low densities, which turns out to be fundamental for the complete description of the propagation phenomenon. It is worth adding that the use of compactly supported subsolution allows them to handle localised initial data and a multidimensional setting in space, two situations we have chosen not to focus on in this work. Indeed, our principal aim was to derive a sharp estimate with the fewest possible assumptions on the dispersal kernel, in particular with the fewest restrictions on the parameter $s$ in Hypothesis \ref{hyp:J}.

\smallskip

Finally, to get some insight on the behaviour of solutions to problem \eqref{eq:main-gen}-\eqref{bc:main-gen} and complement the theoretical results, we conducted a numerical exploration in the fractional diffusion case. The numerical approximation of non-local diffusion operators is notoriously more challenging than that of their local counterparts, this fact being further complicated by the possibly rapid propagation of the solution. An \textit{ad hoc} approach has to be devised to accurately simulate the propagation in an unbounded domain and over a long time frame in order to track a level set and capture the correct behaviour of the solution at low density. The preliminary numerical results presented here constitute our first attempt at using a numerical method to better understand the acceleration phenomenon.

\subsection*{Comments on the strategy}
The first step in proving \Cref{thm:main} is to investigate how the solution evolves from the initial datum for short and large times, and, in particular, what is the decay at infinity induced by a dispersal kernel with fat tails. Proposition \ref{prop:initialtail} shows that a solution to \eqref{eq:main-gen}-\eqref{bc:main-gen} with an initial datum satisfying Hypothesis \ref{hyp:u0} behaves as $x^{-2s}$ at infinity at time $1$. Another important aspect is to know that a positive solution to problem \eqref{eq:main-gen}-\eqref{bc:main-gen}, with an initial datum satisfying Hypothesis \ref{hyp:u0}, flattens through time with the expected decay, that is for any positive number $C$, there exists a positive time $t_C$ such that the solution $u$ to problem \eqref{eq:main-gen}-\eqref{bc:main-gen} satisfies
\[
\lim_{x\to +\infty }x^{2s}u(t,x)\ge C \quad\text{ for all } t\ge t_C.
\]

\begin{figure}[H]
	\centering
	\begin{subfigure}[b]{0.4\textwidth}
 	\centering
 	\includegraphics[width=\linewidth]{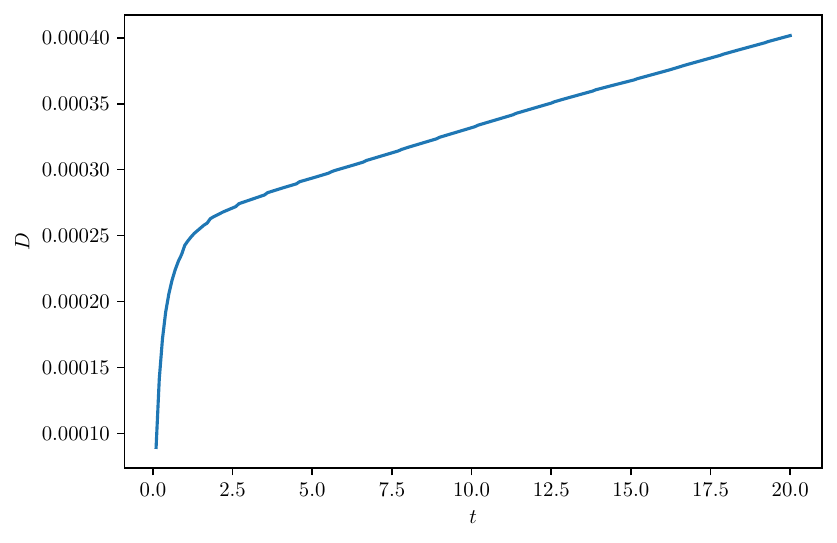}
 	\caption{$s=0.4$.}
 \end{subfigure} %
 \quad
	\begin{subfigure}[b]{0.4\textwidth}
		\centering
		\includegraphics[width=\linewidth]{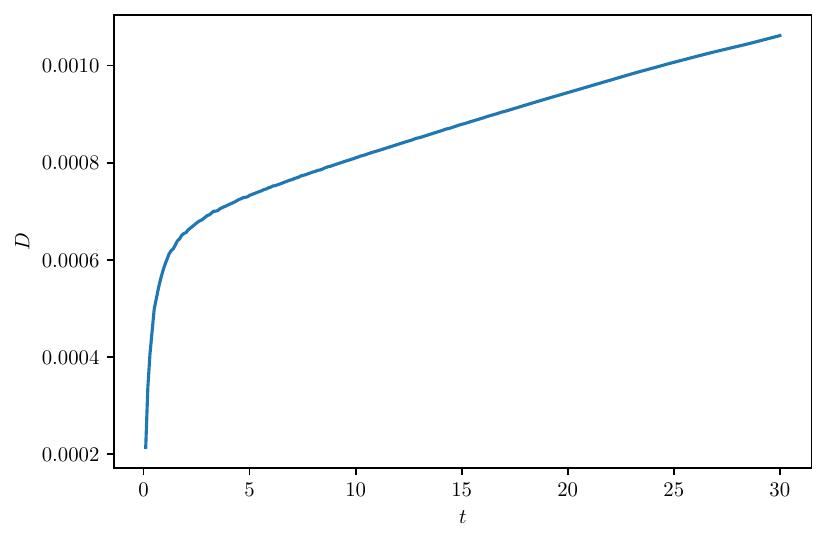}
		\caption{$s=0.5$.}
	\end{subfigure}
\caption{Evolution over time of the constant $D$ obtained by fitting the function $\frac{D}{x^{2s}}$ with the part of the tail between values $10^{-2}$ and $10^{-5}$ of the numerical approximation of the solution to the fractional diffusion-reaction problem, with $\beta=1.5$ and $s=0.4$ (left) or $0.5$ (right).}\label{fig:flat2}
\end{figure}

To illustrate this particular behaviour, we used the method of least-squares to fit the constant in the function $\frac{D}{x^{2s}}$ to the part of the tail between values $10^{-2}$ and $10^{-5}$ of the computed numerical approximation of the solution at different times. The evolution of the fitted constant $D$ is plotted in \Cref{fig:flat2} for a fixed value of $\beta$ and two distinct values of $s$. It is observed that, after a rapid transition, the constant grows linearly with time, which corroborates the flattening with $x^{-2s}$ decay of the the solution at the front. We also refer to \Cref{nose-fit} for various plots showing the correspondence between $u(t,\cdot)$ and $x^{-2s}$ at the edge of the invasion profile.

For particular diffusion operators, such as the fractional Laplace operator, a tail flattening estimate can be obtained through the time and space scaling properties of the associated heat kernel. Although the characterisation of the heat kernel associated to the generator of a L\'evy process is a well-known problem in probability theory and analysis that dates back to the original works of P\'olya \cite{Polya1923} and Blumenthal and Getoor \cite{Blumenthal1960} on $\alpha$-stable processes, characterisations of the kernels inducing such a tail flattening effect have, as far as we know, only been established for some specific classes of processes (see \cite{Bogdan2014,Cygan2017,Grzywny2019,Kaleta2019}) and do not exist for a general one.

In the present paper, the proof of this estimate is based on two ingredients. We first show an invasion property in the context we consider, which is then combined with the construction of a subsolution of the associated linear problem, that mimics the natural scaling behaviour of the original problem. Importantly, this tail flattening property is true for any value of the parameter $s$, and in particular for $s\ge 1$, as shown in \Cref{sec:asymp-esti}. 

The regime $s\ge 1$ is one for which the heat kernel is supposed to behave at large times like a Gaussian diffusion kernel, that is, for any positive integrable initial data, the corresponding solution to the linearised problem, rescaled using the standard self-similar variables, converges to a Gaussian distribution (see \cite{Chasseigne2006} for integrable kernels). This implies that the acceleration and the flattening of a solution to problem \eqref{eq:main-gen}-\eqref{bc:main-gen} cannot be uniquely explained through the scaling properties of the associated heat kernel. In this respect, there is a clear dichotomy between the two regimes $s< 1$ and $s\ge 1$.

Once the tail flattening estimate has been established, the strategy to achieve a lower bound for large times consists in the construction of a new type of subsolution, capturing all the expected dynamics of the solution $u$. In particular, it turns out to be mandatory to identify several regions of space over which the behaviour of $u$ is governed by one specific part of the equation. This appears to be something new compared to previous approaches. Roughly speaking, the dynamics close to $t^{\frac{\beta}{2s(\beta-1)}}$ are due to the nonlinearity only, via the related ordinary differential equation $n'=n^\beta$, the far-field zone is ruled by purely dissipative effects and the solution has the behaviour of that of the linearised equation, \textit{i.e.} $u(t,x) \sim \frac{\kappa t}{x^{2s}}$, and the transition zone between the two sees a subtle interplay occur between the dissipation and the nonlinearity. This trichotomy will be detailed and illustrated in \Cref{sec:strategy}. In relation to what has just been explained, it is interesting to notice that the acceleration exponent is a function of the nonlinearity parameter $\beta$, but not of the way that the tail of the solution flattens with time: it is purely related to the rate of dispersion. This is observed numerically. \Cref{fig:view} presents a schematic view of the expected behaviour of the solution profile. 


\begin{figure}[h]\centering
\begin{tikzpicture}[scale=0.65]
%
\draw[->,line width=0.2, dashed](-8,-0.1)--(8,-0.1) node[below]{$x$};
\draw [->,line width=2, blue] (0.1,2.5)--(2.5,2.5) node[anchor=south]{$\asymp_\lambda t^{\frac{\beta}{2s(\beta-1)}}$};

\draw[line width=1](-6.2,2.5)--(-5.8,2.5) node[above left]{$\lambda\quad$};

\draw[line width=1, dashed](-6,-1.5)--(-6,6.5) node[right]{$ $};
\draw[line width=1, dashed](6,-1.5)--(6,6.5) node[right]{$ $};

%
%
\draw [domain=-7.5:7.5,line width=3, samples=150, red] plot (\x, {(5*exp{-\x})/(1+exp{-\x})}) ;
\fill[blue] (0,2.5) circle (5pt);
\draw[<->,line width=1](-5.2,1)--(5.2,1) node[below]{$\sim t^{\frac{1}{2s}}$};

\end{tikzpicture}
\caption{Schematic view of the expected behaviour of solution at a given time $t$. The solution travels at speed $t^{\frac{\beta}{2s(\beta-1)}}$ and the bulk is $t^{\frac{1}{2s}}$ wide.}
\label{fig:view}
\end{figure}


\subsection*{Further comments and structure of the paper}
From the theory of L\'evy processes \cite{Applebaum2009,Bertoin1996,Sato2013}, it is well known that the semigroup generated by $\opd{\cdot}$ is positive. Consequently, the Cauchy problem associated to the linear equation
\[
\partial_tu(t,x)=\opd{u}(t,x)\text{ for all }t>0\text{ and }x\text{ in }\mathbb{R},
\]
enjoys a parabolic-type comparison principle for mild solutions, that is, if $u_1$ and $u_2$ are two mild solutions to the above evolution equation with $u_1(0,\cdot)\ge u_2(0,\cdot)$, then $u_1(t,\cdot)\ge u_2(t,\cdot)$ for all $t>0$. Such a comparison principle can be extended to the semilinear equation \eqref{eq:main-gen}; we will not provide a proof here, but refer to \cite{Cabre2013, Zhang2023} in the particular case of the fractional Laplace operator. As a consequence of this comparison principle and in view of Hypothesis \ref{hyp:u0}, it is then sufficient to prove \Cref{thm:main} for smooth initial data, say of class $\mathscr{C}^{2}$. Indeed, for any $u_0$ satisfying Hypothesis \ref{hyp:u0}, one can find a monotone function $v_0$ in $\mathscr{C}^{2}(\mathbb{R})$, such that $v_0\le u_0$. It then follows from the comparison principle that, for all $t > 0$ and $x$ in $\mathbb{R}$, $u(t,x)\ge v(t,x)$, and thus $x_\lambda(t)\ge \sup\{x\in\mathbb{R}, v(t,x)\ge \lambda\}$.

In the remainder of the paper, we will assume that the initial datum $u_0$ belongs to $\mathscr{C}^{2}(\mathbb{R})$. For such a smooth initial datum, and since a nonlinearity satisfying Hypothesis \ref{hyp:f} is Lipschitz continuous, a mild solution of the Cauchy problem is then a classical solution, see \cite[chapter 6, Theorem 1.5]{Pazy1983}. Therefore, the considered solution $u$ satisfies problem \eqref{eq:main-gen}-\eqref{bc:main-gen} in the classical sense.

Furthermore, we mention that having a \textit{compactly supported} initial datum would lead to different considerations. In particular, the possibility of invasion is related to the size of the initial datum due to the existence or not of the so-called hair-trigger effect. Depending on the choice of parameters $s$ and $\beta$, it can happen, for a compactly supported initial datum, that the solution gets extinct after a long time, which is referred to as the \textit{quenching phenomenon} \cite{Alfaro2016,Zlatos2005}, and no propagation occurs. We have not chosen to focus on this particular issue in order to concentrate on an accurate description of the acceleration process.

It is important to keep in mind that, from the point of view of applications, having results with assumptions at such a level of generality is of great interest, in particular in ecology where dispersal is a fundamental process which strongly impacts the evolution of species and for which understanding is still partial (see \cite{Nathan2012,Nield2020,Petrovskii2011}). In a sense, by giving access to the correct speed of acceleration for a large class of measures, our results provide a unified view of the consequences of potentially large jumps in the dispersal process. 

It is worth noticing that the numerical graphs in \Cref{fig:flat2} suggest some particular asymptotic behaviour of solutions to problem \eqref{eq:main-gen}-\eqref{bc:main-gen}. For the fractional Laplace operator, we observed numerically the following behaviour: $u(t,x)\sim \tfrac{C_0t}{x^{2s}}$ for large $x$. Such a scaling is indeed satisfied by the subsolution constructed to estimate the speed of level sets from below. However, the supersolution used to control this speed does not. Obtaining rigorously such asymptotic behaviour remains an open question, which requires a more precise description of the supersolution in the spirit of our construction. Some investigations in this direction are currently underway.

Lastly, our approach is rather robust and can be extended to more singular nonlinearities, notably ignition-type ones (see the companion paper \cite{Bouin2021ign}). Let us also mention that an acceleration phenomenon has been identified in some nonlinear reaction-diffusion equations of porous medium type, with or without an Allee effect, see \cite{King2003,Stan2014,Alfaro2017b,Alfaro2019}. We believe that our approach may be of some help in understanding qualitatively these situations.

\smallskip

The paper is organised as follows. In \Cref{sec:asymp-esti}, we derive some estimates on the asymptotic behaviour of solutions to problem \eqref{eq:main-gen}-\eqref{bc:main-gen} and prove the invasion property in Proposition \ref{bcl-prop:inva}. We then describe in broad lines the construction of a subsolution in \Cref{sec:strategy}. The deeper calculations needed for the proof of \Cref{thm:main} are the object of \Cref{sec:estimate1} and \ref{sec:estimate}. Finally, the numerical method employed for the simulations is presented in \Cref{sec:numerics}, with additional experiments.

\section{Tails and flattening estimates}\label{sec:asymp-esti}
\subsection{About the tails of $u$ at $t=1$}\label{subsec:tails}
\noindent In this section, we show that, starting from a Heaviside step function as initial datum, the solution immediately gets polynomial tails of order $2s$, for any positive value of the parameter $s$. For this, we construct a subsolution for short times. 

\noindent Let us introduce the function $v$ defined by
\[
v(t,x)=
\begin{cases}
\frac{1}{\nu} & \text{ for } t>0\text{ and }x\le 0,\\
\frac{\chi t}{x^{2s}+\chi \nu t} & \text{ for } t>0\text{ and }x> 0,
\end{cases}
\]
where $\nu$ and $\chi$ are positive constants to be fixed.\\ Note that $\lim_{t\to0,\ t>0}v(t,\cdot)=\frac{1}{\nu}\mathds{1}_{(-\infty,0]}$.
 
\begin{lem}\label{lem:initialtail}
Assume that the dispersal kernel $J$ satisfies Hypothesis \ref{hyp:J}. For all positive real numbers $\nu$ and $\chi$ verifying $0<\chi\nu\le \frac{1}{2s\mathcal{J}_0}$, one has, for all $t$ $(0,1)$ and $x>R_0+1$,
\[
\partial_tv(t,x) -\opd{v}(t,x) \le \frac{\mathcal{J}_0}{2s} v(t,x).
\]
\end{lem}
\begin{proof}
For $t>0$ and $x>0$, compute
\begin{align*}
&\partial_tv(t,x)=\frac{\chi x^{2s}}{(x^{2s}+\chi\nu t)^2},\quad \partial_xv(t,x)=-2s\frac{\chi t x^{2s-1}}{(x^{2s}+\chi\nu t)^2},\\
&\partial^2_xv(t,x)=2sv^2(t,x)\frac{x^{2s-2}}{\chi t}\left[4s\frac{x^{2s}}{x^{2s}+\nu \chi t} -2s +1 \right].
\end{align*} 
For any $t>0$, note that $x\mapsto v(t,x)$ is a convex and decreasing function. 

Let us estimate $\opd{v}(t,x)$ for all $t>0$ and $x\ge R_0+1$. We first decompose this quantity into three parts as follows: 
\[
\opd{v}(t,x)=\int_{-\infty}^{-1} [v(t,x+z)-v(t,x)]J(z)\,dz 
 +\int_{-1}^{1}[v(t,x+z)-v(t,x)]J(z)\,dz 
 +\int_{1}^{+\infty}[v(t,x+z)-v(t,x)]J(z)\,dz.
\]
Using that $v(t,\cdot)$ is positive for any $t>0$, we have 
\[
\int_{1}^{+\infty}[v(t,x+z)-v(t,x)]J(z)\,dz\ge -v(t,x)\int_{1}^{+\infty}J(z)\,dz.
\]
Similarly, using the definition of $v$ and the fact that $v(t,\cdot)$ is monotone decreasing for any $t>0$ yields
\begin{align*}
\int_{-\infty}^{-1} [v(t,x+z)-v(t,x)]J(z)\,dz &\ge \int_{-\infty}^{-x} [v(t,x+z)-v(t,x)]J(z)\,dz,\\
&\ge \left[\frac{1}{\nu}-v(t,x)\right] \int_{x}^{+\infty}J(z)\,dz,
\end{align*}
so that, using Hypothesis \ref{hyp:J}, we get, for all $t>0$ and $x\ge R_0+1$,
\begin{equation}
\opd{v}(t,x) \ge \int_{-1}^{1} [v(t,x+z)-v(t,x)]J(z)\,dz+\frac{{\mathcal{J}_0}^{-1}}{2s}\left[\frac{1}{\nu}-v(t,x)\right]\frac{1}{ x^{2s}} - \frac{\mathcal{J}_0}{2s} v(t,x).
\end{equation}
The remaining integral is estimated using the regularity of $v$, its convexity with respect to $x$, and the symmetry of $J$. It can be rewritten it as follows
\begin{align*}
\int_{-1}^{1}[v(t,x+z)-v(t,x)]J(z)\,dz&=\int_{0}^1\int_{0}^{1}\int_{-1}^{1}\partial_x^2v(t,x+\tau \sigma z) \tau z^2 J(z)\,dz d\tau d\sigma\\&\geq 0,
\end{align*}
the nonnegativity being infered from the fact that  $\partial_x^2v(t,x+\xi)\ge 0$ for any $\xi$ in $(-1,1)$ since $x\ge R_0+1\ge 2$. We hence conclude that, for all $t>0$ and $x\ge R_0+1$,
\[
\opd{v}(t,x)\ge \frac{{\mathcal{J}_0}^{-1}}{2s}\left[\frac{1}{\nu}-v(t,x)\right] \frac{1}{ x^{2s}} -\frac{\mathcal{J}_0}{2s} v(t,x).
\]
We then have, for all $t$ in $(0,1)$ and $x\ge R_0+1$,
\[
\partial_tv(t,x) -\opd{v}(t,x)\le \frac{\chi x^{2s}}{(x^{2s}+\chi\nu t)^2}-\frac{{\mathcal{J}_0}^{-1}}{2s}\left[\frac{1}{\nu}-v(t,x)\right] \frac{1}{ x^{2s}} +\frac{\mathcal{J}_0}{2s} v(t,x).
\]
Using again the definition of $v$, a quick computation shows that 
\[
\frac{{\mathcal{J}_0}^{-1}}{2s}\left[\frac{1}{\nu}-v(t,x)\right] \frac{1}{ x^{2s}}= \frac{{\mathcal{J}_0}^{-1}}{2s\nu} \frac{1}{x^{2s} +\nu \chi t}.
\]
Since $\frac{x^{2s}}{x^{2s}+\chi \nu t}\le 1$, one gets, for all $t$ in $(0,1)$ and $x\ge R_0+1$,
\[
\partial_tv(t,x) -\opd{v}(t,x)\le  \frac{\chi}{x^{2s}+\chi\nu t} \left(1-\frac{{\mathcal{J}_0}^{-1}}{2s\nu\chi}\right) +\frac{\mathcal{J}_0}{2s} v(t,x) \leq \frac{\mathcal{J}_0}{2s} v(t,x),
\]
as soon as $\chi\nu\le\frac{1}{2s\mathcal{J}_0}$.
\end{proof}

Equipped with the above lemma, we can prove the following result.

\begin{prop}\label{prop:initialtail}
Assume that the dispersal kernel $J$ and the initial datum $u_0$ respectively satisfy Hypotheses \ref{hyp:J}, \ref{hyp:D} and  \ref{hyp:u0}, and let $u$ be a solution to problem \eqref{eq:main-gen}-\eqref{bc:main-gen}. Then, there exists $D > 0$ such that
\[
\lim_{x\to +\infty} x^{2s}u(1,x)\ge 2 D^{2s}.
\]
\end{prop}
\begin{proof}
Observe that, since $u_0$ satisfies Hypothesis \ref{hyp:u0}, it is enough to prove the proposition for a \emph{monotone} initial datum. Indeed, if $u_0$ is monotone nonincreasing then, by a straightforward application of the comparison principle, so is $x\mapsto u(t,x)$ for all $t\geq0$, and we have 
$u(t,x)\ge u(t,R_0+1)$ for all $t >0$ and $x\le R_0+1$. Since $u(t,x)>0$ for all $t>0$ and $x$ in $\mathbb{R}$, we have $\delta:=\inf_{t\in\left[\tfrac12,\tfrac32\right]}u(t,R_0+1)>0$ and thus 
$u(t,x)\ge \delta$ for all $t$ in $\left(\tfrac12,\tfrac32\right)$ and $x = R_0+1$. Consider the function $v$ defined above, with $\nu>\frac{1}{\delta}$ and $\chi$ such that $\chi \nu \leq \frac{1}{2s\mathcal{J}_0}$. For such a choice of constants, the function $\tilde v(t, \cdot):= \left(1 - t \right) e^{-\frac{\mathcal{J}_0}{2s}t}v(t,\cdot)$ satisfies 
\begin{equation}\label{tildez}
\begin{cases}
\partial_t\tilde v(t,x) \le \opd{\tilde v}(t,x) &\quad \text{ for }  t\in (0,1)\text{ and }x>R_0+1,\\
\tilde v(0,\cdot)= \frac{1}{\nu}\mathds{1}_{(-\infty,0]},\\
\tilde v(1,x) = 0, &\quad \text{ for } x \geq R_0+1,\\
\tilde v(t,x)\le \frac{1}{\nu},&\quad \text{ for }  t\in (0,1)\text{ and }x \leq R_0+1.
\end{cases}
\end{equation}
The function $\tilde u := u ( \cdot + \tfrac12 , \cdot)$ satisfies 
\begin{equation}\label{tildeu}
\begin{cases}
\partial_t\tilde u(t,x)\ge \opd{\tilde u} (t,x) &\quad \text{ for }  t\in (0,1)\text{ and } x>R_0+1,\\
\tilde u(0,\cdot) = u ( \tfrac12 , \cdot)> \frac{1}{\nu}\mathds{1}_{(-\infty,0]},\\
\tilde u(1,x) \geq 0, &\quad \text{ for } x \geq R_0+1,\\
\tilde u(t,x) \geq \delta > \frac{1}{\nu},&\quad \text{ for }  t\in (0,1)\text{ and } x \leq R_0+1.
\end{cases}
\end{equation}
Using the parabolic comparison principle, it follows that, for all $t$ in $(0,1)$ and $x\geq R_0+1$, one has $\tilde u(t,x)\ge \tilde v(t,x)$ and thus
\[
\lim_{x\to +\infty} x^{2s}u(1,x)\ge \lim_{x\to +\infty} x^{2s}\tilde v\left(\tfrac12,x\right)=\frac{\chi}{2}\left(1 - \frac{1}{2} \right) e^{-\frac{\mathcal{J}_0}{4s}} \lim_{x\to\infty}\frac{x^{2s} }{x^{2s}+ \tfrac12\chi \nu}, = \frac{\chi}{4} e^{-\frac{\mathcal{J}_0}{4s}} := 2 D^{2s}.
\]
\end{proof}

\subsection{Tail flattening estimates for large times}
We next push further our analysis of the tail of the solution to problem \eqref{eq:main-gen}-\eqref{bc:main-gen} by obtaining the flattening estimate hinted at in the introduction.

\begin{prop}\label{bcl-prop-flatnonlin}
Assume that the dispersal kernel $J$, the initial datum $u_0$ and the nonlinearity $f$ respectively satisfy Hypotheses \ref{hyp:J}, \ref{hyp:D}, \ref{hyp:u0} and \ref{hyp:f},  and let $u$ be a positive solution to problem \eqref{eq:main-gen}-\eqref{bc:main-gen}. Then, for any $C>0$, there exists $t_C>0$ such that
\[
\lim_{x\to +\infty} x^{2s}u(t_C,x)\ge C.
\]
\end{prop}
 
Before proving this result, we shall first establish an invasion property for the solution to \eqref{eq:main-gen}-\eqref{bc:main-gen}. 
 
\begin{prop}\label{bcl-prop:inva}
Assume that the dispersal kernel $J$, the initial datum $u_0$ and the nonlinearity $f$ respectively satisfy Hypotheses \ref{hyp:J}, \ref{hyp:D}, \ref{hyp:u0} and \ref{hyp:f}. Then the solution $u$ to problem \eqref{eq:main-gen}-\eqref{bc:main-gen} satisfies, for all $t\geq0$,
\[
\liminf_{x\to -\infty} u(t,x)\ge \frac{a}{4}.
\]
Moreover, for any positive real number $\mathfrak{R}$, one has
\[
u(t,x)\to 1 \text{ uniformly in $(-\infty,\mathfrak{R}]$ as $t\to +\infty$}.
\]
\end{prop}
\begin{proof}
As above, let us observe that, due to the parabolic comparison principle and since the function $u_0$ satisfies Hypothesis \ref{hyp:u0}, it is enough to prove this proposition for a monotone initial datum. 
Observe also that, when the kernel $J$ is integrable or $\opd{\cdot}$ is the fractional Laplace operator, using the asymptotic global stability of the bistable front (see \cite{Achleitner2015,Chen1997}), the above invasion statement has already been shown in \cite{Alfaro2017,Coville2021}. However, to the best of our knowledge, such an asymptotic stability result has not been established for a generic kernel, and another argument is thus needed in the present context.

Before going into the proof, let us first make a simplification. Since the nonlinearity $f$ satisfies Hypothesis \ref{hyp:f}, we may find a sufficiently small positive number $r_0$ such that $f(z)\ge r_0z^{\beta}(1-z)$ with $\beta>1$ and, using again the parabolic comparison principle, it then is enough to prove the proposition for a nonlinearity of the form $f(z)=r_0z^{\beta}(1-z)$ with $\beta>1$. We thus assume that it is the case here, the idea being to construct a subsolution to equation \eqref{eq:main-gen} that fills the whole space. We first observe that, for any nonnegative nonlinearity $f$, any positive function $v$ and $R\ge R_0$, one has, for all $x$ in $\mathbb{R}$, 

\begin{align*}
\opd{v}(x)&= \int_{-\infty}^{+\infty}[v(x+h)-v(x)] J(h)\,dh \\
&= \int_{|h|<R}[v(x+h)-v(x)] J(h)\,dh +\int_{|h|\ge R} [v(x+h)-v(x)] J(h)\,dh\\
&\ge \int_{|h|<R}[v(x+h)-v(x)] J(h)\,dh -v(x)\int_{|h|\ge R} J(h)\,dh, 
\end{align*}
and so
\[
\opd{v}(x)+f(v(x))\ge \int_{|h|<R}[v(x+h)-v(x)] J(h)\,dh -\frac{\mathcal{J}_0}{R^{2s}}v(x) +f(v(x)).
\]
Set $f_R(s):=-\frac{\mathcal{J}_0}{R^{2s}}s +f(s)$ and denote by $\D_{J_{R}}[\,\cdot\,]$ the dispersal operator with kernel $J\mathds{1}_{B_R(0)}$. From the above computations, for any positive solution $u$ to \eqref{eq:main-gen}, we have, for all $t>0$ and $x$ in $\mathbb{R}$,
\[
 \partial_t u(t,x) -\opdr{u}(t,x)-f_R(u(t,x))\ge \partial_t u(t,x) -\opd{u}(t,x)-f(u(t,x))=0.
\]
Let $0<\theta< a=\liminf_{x\to+\infty} u_0(x)$.\\ For $\theta<\frac{1}{2}$, define $f_\theta(z):=\eps_0 z\left(z-\frac{1}{4}\right)(1-\theta-z)$, with $\eps_0$ positive and small enough so that $f(z)> f_{\theta}(z)>0$ for all $z$ in $\left[\frac{1}{4},1-\theta\right]$. We next choose $\theta<\frac{a}{8}$ small, so that $1-\theta>a$ and $\int_{0}^{1-\theta}f_{\theta}(z)\,dz>0$.\\
Since $f_R$ tends to $f$ as $R$ tends to $+\infty$, we may find $R_\theta$ such that $f_\theta\le f_{R_\theta}$ and so, for $R\ge R_\theta$, we have, for all $t>0$ and $x$ in $\mathbb{R}$,
\begin{equation}\label{bcl:eq:subsol-inva}
 \partial_t u(t,x) -\opdr{u}(t,x)-f_\theta(u(t,x))\ge 0.
\end{equation}
Let us define an extension of class $\mathscr{C}^1$ of $f_\theta$ outside $[0,1-\theta]$ as follows:
$$
\begin{cases} 
{f_\theta}'(0)z \quad &\text{ when } z<0,\\
f_{\theta}(z) \quad &\text{ when } 0\le z \le 1-\theta,\\
{f_{\theta}'}(1-\theta)(z-1+\theta) \quad &\text{ when } 1-\theta<z,
\end{cases}
$$
and, for convenience, let us still denote by $f_\theta$ this extension. Considering the equation, for all $t>0$ and $x$ in $\mathbb{R}$,
\begin{equation}\label{eq-main-approx}
\partial_t v(t,x) -\opdr{v}(t,x)-f_\theta(v(t,x))= 0,
\end{equation}
it follows from \eqref{bcl:eq:subsol-inva} that $u$ is a supersolution to \eqref{eq-main-approx}. We next construct an adequate subsolution to \eqref{eq-main-approx}. From \cite{Alberti1998,Chmaj2013,Fang2015}, we know that equation \eqref{eq-main-approx} admits a unique monotonically decreasing travelling front solution $(\varphi_\theta,c_\theta)$, connecting $1-\theta$ to $0$. Thanks to the results in \cite{Coville2005}, without any loss of generality and by reducing $\theta$ if necessary, we can assume that this front is smooth, that is $(\varphi_\theta,c_\theta)$ is a smooth solution to 
\begin{align*}
&c_\theta {\varphi_\theta}'(x) +\opdr{\varphi_\theta}(x) +f_\theta(\varphi_{\theta}(x)) =0 \quad\text{ for all } x\text{ in }\mathbb{R},\\
&\lim_{x\to -\infty}\varphi_\theta(x) =1-\theta, \qquad \lim_{x\to +\infty}\varphi_\theta(x) =0.
\end{align*}

By the definition of $f_\theta$, we must have $c_\theta>0$, since the sign of the speed in such context is given by the sign of the integral $\int_{0}^{1-\theta}f_{\theta}(z)\,dz$. We then normalise $\varphi_\theta$ to have $\varphi_\theta(0)=\frac{a}{2}$ and set 
\[
\Psi(t,x):=\varphi_{\theta}\left(x-c_\theta t + \varsigma(1-e^{-\eps t})\right) -\left(1-\frac{a}{2}\right)e^{-\eps t},
\]
with $\eps$ and $\varsigma$ some parameters. Observe that at $t=0$, we have, for all $L>0$, 
\begin{align*}
&\Psi(0,x+L)\le \frac{a}{2} -\theta <a \quad \text{ for all } \varsigma,\  \eps,\ x,\\
&\Psi(0,x+L)\le a -1 <0 \quad \text{ for all } \varsigma>0,\ \eps>0,\ x>-L .
\end{align*}
As a consequence, since $\varphi_\theta$ is monotone for any fixed $\eps$ and $\varsigma$, we can always find $L_0$ such that $u_0\ge \Psi(0,\cdot+L_0)$. For an adequate choice of $\eps$ and $\varsigma$, the function $\Psi(\cdot,\cdot+L_0)$ is a subsolution to~\eqref{eq-main-approx}.
\begin{claim}\label{claim subsolution}
There exist values for the parameters $\eps$ and $\varsigma$ such that, for any $L$, $\Psi(\cdot,\cdot+L)$ is a subsolution to equation \eqref{eq-main-approx}.
\end{claim}
Let us postpone the proof of this claim for the moment. Equipped with its conclusion and using the parabolic comparison principle, we deduce that for any $t\ge 0$ and $x$ in $\mathbb{R}$, $u(t,x)\ge \Psi(t,x+L_0)$. Thus, for any $t\ge 0$, we have 
$$
\liminf_{x\to -\infty} u(t,x) \ge \liminf_{x\to -\infty} \Psi(t,x+L_0)=1-\theta-\left(1 -\frac{a}{2}\right)e^{-\eps t}\ge \frac{a}{2}-\theta\ge \frac{a}{4},
$$
since $\theta\le\frac{a}{8}$.\\ In addition, we deduce that, for any real number $\mathfrak{R}$, $\lim_{t\to+\infty}u(t,\cdot)\ge 1-\theta$ in $(-\infty,\mathfrak{R}]$. The parameter $\theta$ being arbitrary small the latter argument then implies that $u(t,\cdot)$ tends to $1$ locally uniformly in $(-\infty,\mathfrak{R}]$ as $t$ tends to $+\infty$. Since $u$ is monotone nonincreasing in $x$, this convergence is uniform. 
\end{proof}

To complete the proof of the above proposition, we establish the claim.

\begin{proof}[Proof of Claim \ref{claim subsolution}]
Computing $\partial_t \Psi$, one has, for all $t\geq0$ and $x$ in~$\mathbb{R}$,
\[
\partial_t \Psi(t,x+L)=(-c_\theta+\eps \varsigma e^{-\eps t}){\varphi_\theta}'\left(x-c_\theta t + \varsigma(1-e^{-\eps t})+L\right)\\+\eps\left(1-\frac{a}{2}\right)e^{-\eps t}.
\]
Set $\xi(t,x):=x-c_\theta t + \varsigma(1-e^{-\eps t})+L$. Using the equation satisfied by $\varphi_\theta$, we have, for all $t\geq0$ and $x$ in $\mathbb{R}$,
\[
\partial_t \Psi(t,x+L)-\opdr{\Psi}(t,x+L)-f_{\theta}(\Psi(t,x+L))=RHS(t,x)\]
with 
\[
RHS(t,x):=\eps \varsigma e^{-\eps t}{\varphi_\theta}'\left(\xi(t,x)\right)+\eps\left(1-\frac{a}{2}\right)e^{-\eps t}+f_\theta(\varphi_{\theta}(\xi(t,x)))-f_{\theta}\left(\varphi_{\theta}(\xi(t,x)) -\left(1-\frac{a}{2}\right)e^{-\eps t}\right).
\]
To conclude, we only have to show that $RHS(t,x)\le 0$, for all $t>0$ and $x$ in $\mathbb{R}$, for the right choice of the parameters $\delta,\eps$ and $\varsigma$.
We next choose $0<\delta_0<\frac{a}{8}$ such that $f_\theta$ satisfies 
\begin{align*}
&f_\theta(z)\le \frac{f'_\theta(0)}{2}z \quad \text{ for } z\text{ in }(0,\delta_0),\\
&\frac{3}{2}{f_\theta}'(1-\theta)\le {f_\theta}'(z)\le \frac{1}{2}{f_\theta}'(1-\theta) \quad \text{ for } z\text{ in }(1-\theta -4\delta_0, 1-\theta).
\end{align*}
Then, taking inspiration in the construction in \cite{Brasseur2021}, let $\delta<\delta_0$ and choose $A(\delta)$ such that $\varphi_{\theta}(z)\le \delta$ if $z\ge A$ and $\varphi_{\theta}(z)\ge 1-\theta -\delta $ for $z\le -A$.
We now distinguish between the three situations $\xi(t,x)>A$, $\xi(t,x)<-A$ and $|\xi(t,x)|\le A$ and treat each of them separately. 
\paragraph*{The case $\xi(t,x)>A$.} In this case, there are two possibilities, either $ \varphi_{\theta}(\xi(t,x)) -\left(1-\frac{a}{2}\right)e^{-\eps t}>0$ or $\varphi_{\theta}(\xi(t,x)) -\left(1-\frac{a}{2}\right)e^{-\eps t}\le 0 $. With the latter one, we have $f_{\theta}(\varphi_{\theta}(\xi(t,x))) \le \frac{f_\theta'(0)}{2} \varphi_{\theta}(\xi(t,x))$ and
\[
f_\theta\left(\varphi_{\theta}(\xi(t,x)) -\left(1-\frac{a}{2}\right)e^{-\eps t}\right)=f'_{\theta}(0) \left[\varphi_{\theta}(\xi(t,x)) -\left(1-\frac{a}{2}\right)e^{-\eps t}\right].
\]
Since ${\varphi_\theta}'<0$, we have
\[
RHS(t,x)\le
\eps\left(1-\frac{a}{2}\right)e^{-\eps t}+\frac{f'_\theta(0)}{2}\varphi_{\theta}(\xi(t,x)) \\-f'_\theta(0)\left[\varphi_{\theta}(\xi(t,x)) -\left(1-\frac{a}{2}\right)e^{-\eps t}\right]
\]
and therefore since $-f'_{\theta}(0)> -\frac{f'_\theta(0)}{2}>0$, we get
\[
RHS(t,x)
\le \left[\eps+\frac{f'_\theta(0)}{2}\right]\left(1-\frac{a}{2}\right)e^{-\eps t}\le 0,
\]
as soon as $\eps \le -\frac{f'_\theta(0)}{2}$. In the other situation, we have $$\delta \ge \varphi_{\theta}(\xi(t,x)) -\left(1-\frac{a}{2}\right)e^{-\eps t}\ge 0$$ and therefore
\[
f_\theta(\varphi_{\theta}(\xi(t,x))) -f_{\theta}\left(\varphi_{\theta}(\xi(t,x)) -\left(1-\frac{a}{2}\right)e^{-\eps t}\right)\le \frac{f'_{\theta}(0)}{2} \left(1-\frac{a}{2}\right)e^{-\eps t}.
\]
As above, we conclude that 
\[
RHS(t,x)\le \left[\eps+\frac{f'_\theta(0)}{2}\right]\left(1-\frac{a}{2}\right)e^{-\eps t}\le 0,
\]
as soon as $\eps \le -\frac{f'_\theta(0)}{2}$.

\paragraph*{The case $\xi(t,x)<-A$.}
In this situation, let us first  assume that
\[
\left(1-\frac{a}{2}\right)e^{-\eps t}\le 3\delta_0.
\]
Then, one has
\[
\varphi_\theta(\xi(t,x)) - \left(1-\frac{a}{2}\right)e^{-\eps t}\ge 1-\theta -4\delta_0,
\]
 and therefore
\[
f_\theta(\varphi_{\theta}(\xi(t,x))) -f_{\theta}\left(\varphi_{\theta}(\xi(t,x)) -\left(1-\frac{a}{2}\right)e^{-\eps t}\right) \le \frac{{f_\theta}'(1-\theta)}{2} \left(1-\frac{a}{2}\right)e^{-\eps t}.
\]
It follows that
\[
RHS(t,x)\le \left[\eps+\frac{{f_\theta}'(1-\theta)}{2}\right]\left(1-\frac{a}{2}\right)e^{-\eps t}\le 0,
\]
provided that $\eps\le -\frac{{f_{\theta}}'(1-\theta)}{2}$. Otherwise, one has $\left(1-\frac{a}{2}\right)e^{-\eps t}>3\delta_0$, hence 
\[
\frac{a}{2}-\theta -\delta =1-\theta -\delta -(1-\frac{a}{2})\le \varphi_{\theta}(\xi(t,x)) -\left(1-\frac{a}{2}\right)e^{-\eps t}<1-\theta -3\delta_0.
\]
As a consequence, since $\delta\le \frac{a}{8}$ and $0\le\theta\le\frac{a}{8}$, one has
\[
1-3\delta_0>\varphi_\theta(\xi(t,x)) -\left(1-\frac{a}{2}\right)e^{-\eps t} \ge \frac{a}{4}>\theta.
\]
By the definition of $f_\theta$, we can ensure that 
\[
f_\theta\left( \varphi_{\theta}(\xi(t,x)) -\left(1-\frac{a}{2}\right)e^{-\eps t}\right)\ge m_0:=\min_{s\in \left[\frac{a}{4}, 1-3\delta_0\right] } f_{\theta}(s).
\]
In addition, using that $\varphi_\theta(\xi(t,x))\ge 1-\theta -\delta$, it follows that
\begin{align*}
f_\theta( \varphi_{\theta}(\xi(t,x)))&=f_\theta( \varphi_{\theta}(\xi(t,x)))-f_\theta(1-\theta)\\&\le -\frac{3}{2}{f_\theta}'(1-\theta)\left(1-\theta-\varphi_{\theta}(\xi(t,x))\right)\\&\le -\frac{3}{2}{f_\theta}'(1-\theta)\delta.
\end{align*}
In conclusion, one gets
\[
RHS(t,x)\le\eps\left(1-\frac{a}{2}\right)-\frac{3{f_\theta}'(1-\theta)}{2}\delta -m_0\le 0,
\]
provided that the values of $\eps$ and $\delta$ are chosen small enough, for instance $\eps \le \frac{m_0}{2-a}$ and $\delta\le \frac{m_0}{3{f_\theta}'(1-\theta)}$.

\paragraph*{The case $|\xi(t,x)|\le A$.}  In that region, one has ${\varphi_\theta}'<0$ and therefore
\[
{\varphi_\theta}'(\xi(t,x))\le-\nu_0:=\sup_{z\in [-A,A]} {\varphi_\theta}'(z)<0.
\]
Recalling that $f_{\theta}$ is a Lipschitz continuous function, one also has 
\[
f_\theta(\varphi_{\theta}(\xi(t,x))) -f_{\theta}\left(\varphi_{\theta}(\xi(t,x)) -\left(1-\frac{a}{2}\right)e^{-\eps t}\right) \le \|{f_{\theta}}'\|_\infty \left(1-\frac{a}{2}\right)e^{-\eps t},
\]
and thus ends up with 
\begin{align*}
RHS(t,x)&\le -\varsigma \eps e^{-\eps t}\nu_0+(\eps+\|{f_{\theta}}'\|_\infty)\left(1-\frac{a}{2}\right)e^{-\eps t}\\
&\le \left(-\varsigma \eps \nu_0+(\eps+\|{f_\theta}'\|_\infty)\left(1-\frac{a}{2}\right)\right)e^{-\eps t}\\
&\le 0,
\end{align*}
provided that the value of $\varsigma$ is chosen large enough, $\varsigma\ge \frac{(\eps+\|{f_\theta}'\|)\left(2-a\right)}{2\eps \nu_0}$ for instance.
\end{proof}

\begin{remark}
The above proof makes no assumption of a specific form for the nonlinearity $f$, only that it is of monostable-type, in the sense that $f(0)=f(1)=0$ and $f>0$ in $(0,1)$. As a consequence, it holds for any monostable nonlinearity. In addition, with some minor adaptations in the manner the bistable function $f_\theta$ is constructed, the proof will also be valid for an ignition-type nonlinearity. In such a situation, we can work with the bistable nonlinearity $f_\theta(z):=\frac{1}{2}f(z) \mathcal{G}(z)+\eps_0 z(\frac{\rho+1}{2}-z)(1-\theta-z)$ where $\rho$ is the threshold of the ignition function $f$ and  the function $\mathcal{G}(\cdot)$ is a smooth regularization of the function $\mathds{1}_{\{z\le\frac{\rho+1}{2}\}}(\cdot)$ that satisfies $\mathcal{G}(\frac{2\rho+3}{4})=0$.
\end{remark}

Let us now prove Proposition \ref{bcl-prop-flatnonlin}.

\begin{proof}[Proof of Proposition \ref{bcl-prop-flatnonlin}]
To construct an adequate subsolution, we consider the parametric function $v$ defined in the beginning of \Cref{subsec:tails}, in which we set $\nu=2$, that is
\[
v(t,x)=
\begin{cases}
\frac{1}{2} &\quad \text{ for } t>0\text{ and }x\le 0,\\
\frac{\chi t}{x^{2s}+2 \chi t} &\quad \text{ for } t>0\text{ and }x> 0.
\end{cases}
\]
and we assume that $\chi=\frac{{\mathcal{J}_0}^{-1}}{8s}$. We now estimate $\opd{v}(t,x)$ for all $t>0$ and $x$ in $\mathbb{R}$.

Let $R$ be a real number greater than $1$, chosen as in the proof of Proposition \ref{prop:initialtail}. By reusing the computation in the proof of Lemma \ref{lem:initialtail}, one can check that, for all $t>0$ and $x\ge R_0+R$,
\[
\opd{v}(t,x)\ge \int_{-R}^{R} [v(t,x+z)-v(t,x)]J(z)\,dz +\frac{{\mathcal{J}_0}^{-1}}{2s}\left[\frac{1}{2}-v(t,x)\right] \frac{1}{ x^{2s}} - \frac{\mathcal{J}_0}{2sR^{2s}} v(t,x).
\]
As in the proof of Lemma \ref{lem:initialtail}, the above integral is shown to be nonnegative using the symmetry of $J$, the regularity of $v$ and its convexity with respect to $x$ for fixed $t$. Thus, one has, for all $t>0$ and $x>R_0+R$, 
\[
\opd{v}(t,x)\ge\frac{{\mathcal{J}_0}^{-1}}{2s}\left[\frac{1}{2}-v(t,x)\right] \frac{1}{ x^{2s}} -\frac{\mathcal{J}_0}{2sR^{2s}} v(t,x).
\]
Altogether, we have, for all $t>0$ and $x\ge R_0+R$,
\[
\partial_tv(t,x) -\opd{v}(t,x)\le \frac{\chi x^{2s}}{(x^{2s}+2\chi t)^2} -\frac{{\mathcal{J}_0}^{-1}}{2s}\left[\frac{1}{2}-v(t,x)\right] \frac{1}{ x^{2s}} +\frac{\mathcal{J}_0}{2s R^{2s}} v(t,x),
\]
which, by observing that $\frac{x^{2s}}{x^{2s}+2\chi t}\le 1$ and using the expressions for $v$ and $\chi$, reduces to 
\begin{align*}
\partial_tv(t,x) -\opd{v}(t,x)&   
\le \frac{1}{x^{2s}+2\chi t} \left(-\frac{{\mathcal{J}_0}^{-1}}{8s} + \frac{\mathcal{J}_0\chi t}{2s R^{2s}}\right).
\end{align*}
For any $C>0$, let us now define $t^*:=\frac{2C}{\chi}$ and choose $R$ large enough, for instance $R\ge R_C:= \left(8C{\mathcal{J}_0}^2\right)^{\frac{1}{2s}}$. From the above computation, we then have, for all $0< t <t^*$ and $x\ge R_0+R_C$,
\begin{equation}\label{bcl2:eq-subsol-flatnonlin}
\partial_tv(t,x) -\opd{v}(t,x)\le 0.
\end{equation}

Equipped with this subsolution, let us now conclude. By using the invasion property stated in Proposition \ref{bcl-prop:inva}, there exists $t_C$ such that, for all $t\ge t_C$ and $x\le R_0+R_C$,
\[
u(t,x)\ge \frac{3}{4}.
\]
The function $\tilde u:=u(\cdot+t_C,\cdot)$ then satisfies 
\begin{align*}
&\partial_t\tilde u(t,x) -\opd{\tilde u}(t,x)\ge 0 \quad \text{ for all }0< t <t^*\text{ and }x\text{ in }\mathbb{R},\\
& \tilde{u}(t,x) \ge v(t,x) \quad \text{ for }0\le t \le t^*\text{ and }x\le R_0+R_C.
\end{align*}
Using the comparison principle, it follows that, for all $0\le t \le t^*$ and $x\geq R_0+R_C$, $\tilde u(t,x)\ge v(t,x)$, from which we deduce that
\[
\lim_{x\to +\infty} x^{2s}u\left(t_C+\tfrac{t^*}{2},x\right)\ge \lim_{x\to +\infty} x^{2s} v\left(\tfrac{t^*}{2},x\right)=C.
\] 
\end{proof}

\section{Strategy for the construction of subsolutions}\label{sec:strategy}
Our goal is to construct a function that mimics some expected asymptotic behaviour and can be compared for large times to the solution of problem \eqref{eq:main-gen}-\eqref{bc:main-gen}. As observed in the previous section, since the nonlinearity satisfies Hypothesis \ref{hyp:f}, it is enough to construct a subsolution for equation \eqref{eq:main-gen} with a nonlinearity having the form $f(z)=r_0z^\beta(1-z)$, with $r_0$ a small enough positive number. One may also notice that, by scaling in both time and space the solution $u$ as well as the kernel $J$, that is considering $v:=u\left(\frac{\cdot}{r_0},\frac{\cdot}{r_0}\right)$ and $J(\frac{\cdot}{r_0})$, we can reduce the construction to finding a subsolution to the following equation
\begin{equation}\label{eq:main}
\partial_t v=\mathcal{D}_{r_0}[v]+v^{\beta}(1-v).
\end{equation}
where $\mathcal{D}_{r_0}[\,\cdot\,]$ denotes the dispersal operator associated to the rescaled measure $J(\frac{z}{r_0})\,dz$. In the sequel, to keep tractable notations, we will drop the subscript of this scaled operator. In addition, we will assume that Hypotheses \ref{hyp:J}, \ref{hyp:f} and \ref{hyp:u0}, and inequality \eqref{s-beta rel} for the parameters $s$ and $\beta$, hold throughout.

\subsection{Form of a subsolution}
We are looking for a function $\undu$ satisfying, for all $t\geq1$ and $x$ in $\mathbb{R}$,
\begin{equation}\label{bcl2-eq:subsol} 
\partial_t\undu(t,x) \leq \opd{\undu}(t,x) + (1-\lambda)\undu^\beta(t,x) \qquad \text{and} \qquad \undu(t,x)\leq\lambda, 
\end{equation}
for some $\lambda$ in $(0,1)$. This function seems to be a natural candidate to compare with the solution $u$ to problem \eqref{eq:main-gen}-\eqref{bc:main-gen}. Indeed, if $\undu(1,\cdot) \leq u(t',\cdot)$ for some $t' > 0$,  since, for all $t>0$ and $x$ in $\mathbb{R}$,
\begin{equation*}
\partial_t\undu(t,x) \leq \opd{\undu}(t,x) + (1-\lambda)\undu^\beta(t,x) \leq \opd{\undu}(t,x) + (1-\undu(t,x))\undu^\beta(t,x)
\end{equation*}
then, by the parabolic comparison principle, it will follow that $u(t+t'-1,\cdot)\ge \undu(t,\cdot)$ for all $t\ge 1$.

In our construction, this subsolution is a piecewise function of class $\mathscr{C}^2$ at least, of the form
\begin{equation}\label{bcl2-def:undu}
\undu(t,x):=\begin{cases}\lambda&\quad\text{for } x\le X(t),\\ 
\phi(t,x) &\quad\text{elsewhere}.
\end{cases}
\end{equation}
in which the function $\phi$ is such that $\phi(t,X(t))=\lambda$, the point $X(t)$ being unknown at this stage. We expect $\phi$ to be related to a solution of an ordinary differential equation of the form $n'=n^\beta$ near $x = X(t)$ and to look like a solution to a standard fractional diffusion-reaction equation with Heaviside step function as initial datum at the far edge. A possible candidate would be given by 
\begin{equation}\label{bcl2-def:w}
w(t,x):= \left[ \left(\frac{\kappa t}{x^{2s}}\right)^{1-\beta} -\gamma(\beta-1)t\right]^{-\frac{1}{\beta-1}},
\end{equation}
where the positive real numbers $\kappa$ and $\gamma$ are free parameters to be chosen later.\\This function is well-defined for $t\ge 1$ and $x>X_{0} :=\kappa^{\frac{1}{2s}}\left(\gamma(\beta-1)t\right)^{\frac{\beta}{2s(\beta-1)}}$ and has the structure of a solution to the ordinary differential equation $n' = n^\beta$. Moreover, the expected decay in space of a solution to the standard fractional Laplace equation with Heaviside step function as initial data \cite{Blumenthal1960,Bogdan2014,Cygan2017} being at least of order $tx^{-2s}$, it possesses the correct asymptotic behaviour. The condition $w(t,X(t))=\lambda$ leading to
\begin{equation}\label{bcl2-def:X}
X(t)=(\kappa t)^{\frac{1}{2s}}\left[ \lambda^{1-\beta} +\gamma(\beta -1)t\right]^{\frac{1}{2s(\beta-1)}}\text{ for all }t\geq1,
\end{equation}
one may observe that the position $X(t)$ moves with the velocity expected from \Cref{thm:main}. However, taking $\phi$ equal to $w$ would not lead to a $\mathscr{C}^2$-function at $x=X(t)$. To remedy this issue, we achieve the construction by setting, for all $t\geq1$ and $x\ge X(t)$,
\begin{equation}\label{bcl2-def:phi}
\phi(t,x):=3\left(1-\frac{w(t,x)}{\lambda}+\frac{w^{2}(t,x)}{3\lambda^2}\right)w(t,x).
\end{equation}

The corresponding $\undu$ then satisfies \eqref{bcl2-eq:subsol} if and only if, for all $t\ge 1$,
\begin{align}
0 &\leq \opd{\undu}(t,x) + \lambda^\beta (1-\lambda), \quad\text{for all }x \leq X(t),\label{eq:subsolleft}\\
\partial_t\undu(t,x) &\leq \opd{\undu}(t,x) + (1-\lambda)\undu^\beta(t,x), \quad\text{elsewhere}.\label{eq:subsolright}
\end{align}
Consequently, the main task is to derive estimates for $\opd{\undu}(t,x)$ in both regions $x \leq X(t)$ and $x \geq X(t)$. The estimate in the first region is rather direct to obtain and relies mostly on the fact that $\undu$ is constant there together with the tails of $J$. Getting the estimate in the second region requires a more intricate approach, for which the region is divided into three zones, as depicted in \Cref{fig:zones} below, each of them being the stage of one specific character of the model and thus demanding a specific way of estimating $\opd{\undu}(t,x)$. 

\begin{figure}[h]
\centering
\resizebox{\textwidth}{!}{
\begin{tikzpicture}[scale=1]
%
\draw[->,line width=0.2, dashed](-8,0)--(8,0) node[below]{$x$};

\draw[line width=1, dashed](-3,-1.5)--(-3,6.5) node[right]{$x=X(t)$};
\draw[line width=1, dashed](0,-1.5)--(0,6.5) node[right]{$x=Y(t)$};
\draw[line width=1, dashed](4,-1.5)--(4,6.5) node[right]{$x= 2^{\frac{1}{2s(\beta-1)}}X(t)$};
%
\fill[pattern=north west lines, pattern color=blue, opacity=0.5] 
(-7.5,-1)--(-3,-1)--(-3,6)
-- (-7.5,6)
-- cycle;

\fill[pattern=north west lines, pattern color=orange, opacity=0.5] 
(-3,-1)--(-3,6) -- (0,6)--(0,-1)
-- cycle;

\fill[pattern=north west lines, pattern color=black!50!green,opacity=0.5] 
(0,-1)--(0,6)--(4,6)--(4,-1)-- cycle;


\fill[pattern=north west lines,pattern color=brown, opacity=0.5] 
(4,-1)--(4,6)--(7.5,6)--(7.5,-1)-- cycle;

%
%
\draw [domain=-3:7.5,line width=3, samples=150, red] plot (\x, {3/(1+0.1*(\x+3)^2)}) ;
\draw [line width=3,red] (-3,3)--(-7.5,3) node[anchor=south west]{$\lambda $};
\end{tikzpicture}
}
\caption{Schematic view of the subsolution $\undu$ at a given time $t$. For estimations, several zones are considered. The blue zone is where $\undu$ is constant, making computations easier. In the orange zone, we use crucially the fact that $\undu$ looks like a solution to an ordinary differential equation of the form $n'=n^\beta$. The exact expression of $Y(t)$ will appear naturally later. In the brown (far-field) zone, a decay imitating that of the solution to a fractional Laplace equation provides the correct asymptotic behaviour. Finally, the construction in the intermediate green zone is more subtle and based on a mixture of the treatments in both surrounding zones.}
\label{fig:zones}
\end{figure}

\subsection{Facts and formulas on $X$ and $w$}
First, from direct computations, we have:
\begin{align}
&\partial_t\undu(t,x) = \partial_x\undu(t,x) = \partial_x^2\undu(t,x) =0 &\text{for all } t\ge 1\text{ and }x\le X(t),\label{bcl2-def:u_tx}\\
&\partial_t\undu(t,x) = 3 \partial_tw(t,x) \left(1- \frac{w(t,x)}{\lambda} \right)^2&\text{for all } t\ge 1\text{ and }x>X(t),\label{bcl2-def:u_t}\\
&\partial_x\undu(t,x) = 3 \partial_xw(t,x) \left(1- \frac{w(t,x)}{\lambda} \right)^2&\text{for all } t\ge 1\text{ and }x>X(t),\label{bcl2-def:u_x}
\end{align}
and for all $ t\ge 1\text{ and }x>X(t)$:
\begin{align}
&\partial_x^2\undu(t,x)=3\left(1-\frac{w(t,x)}{\lambda}\right)\left[\partial_x^2w(t,x)\left(1-\frac{w(t,x)}{\lambda}\right)-\frac{2(\partial_xw(t,x))^2}{\lambda}\right]\label{bcl2-def:u_xx}
\end{align}
Note crucially that $\undu$ is a function of class $\mathscr{C}^1$ in $t$ and of class $\mathscr{C}^{2}$ in $x$. For convenience, let us denote 
\begin{equation}\label{eq:phiU}
\Phi(t,x) :=\frac{\kappa t}{x^{2s}},\qquad U:=\frac{w}{\Phi}.
\end{equation}
We will repeatedly need the following information on derivatives of $w$ at any point $(t,x)$ where $w$ is defined,
\begin{align}
\partial_tw(t,x)&=w^\beta(t,x) \left( \gamma + \frac{\partial_t\Phi(t,x)}{\Phi^\beta(t,x)} \right),\label{bcl2-def:w_t}\\
\partial_xw(t,x)&=w^\beta(t,x)\frac{\partial_x\Phi(t,x)}{\Phi^\beta(t,x)}\nonumber\\
&=U^\beta(t,x)\partial_x\Phi(t,x)=-2sw^\beta(t,x)\frac{x^{2s(\beta -1)-1}}{(\kappa t)^{\beta-1}},\label{bcl2-def:w_x}\\
\partial_x^2w(t,x)&=\beta w^{\beta-1}(t,x)\partial_xw(t,x)\frac{\partial_x\Phi(t,x)}{\Phi^\beta(t,x)}\nonumber\\
&\quad + w^\beta(t,x) \frac{\partial_x^2\Phi(t,x)\Phi^\beta(t,x)-\vert\partial_x\Phi(t,x)\vert^2 \beta \Phi^{\beta-1}(t,x) }{\Phi^{2\beta}(t,x)},\label{bcl2-def:w_xx}\\ 
&=\left[\partial_x^2\Phi(t,x)+ \beta\vert\partial_x\Phi(t,x)\vert^2 \Phi^{-1}(t,x) \left( U^{\beta-1}(t,x)- 1\right)\right]U^{\beta}(t,x).\nonumber
\end{align}
Since $U\geq 1$, it follows from the latter identity that $w$ is convex with respect to $x$. In addition, by rewriting the second-order partial derivative $\partial_x^2\undu$ in terms of $U$ and $\Phi$, we observe that 
\[\frac{\partial_x^2\undu(t,x)}{ 3U^{2\beta}(t,x) \vert\partial_x\Phi(t,x)\vert^2\left(1-\frac{w(t,x)}{\lambda}\right)} =\frac{1}{\lambda}\left(\mathscr{G}(t,x) -2\right),\]
with 
\[\mathscr{G}(t,x):=\frac{ \lambda\left[ \left(1+\frac{1}{2s} \right) + \beta \left( U^{\beta-1}(t,x)- 1\right)\right]}{\Phi(t,x) U^{\beta}(t,x)} \left(1- \frac{w(t,x)}{\lambda} \right). 
\]
By using that $U\Phi=w$,  $\mathscr{G}(t,x)$ can be rewritten as
\[\mathscr{G}(t,x)= 
\left[ \left(1+\frac{1}{2s} -\beta\right)U^{1-\beta}(t,x) + \beta \right] \left(\lambda w^{-1}(t,x)- 1 \right).
\]
 
From the above expression, we can observe that $\undu(t,x)$ is convex with respect to $x$ for $t\geq1$ and as soon as $x$ is such that
\[\left[ \left(1+\frac{1}{2s} - \beta \right) U^{1-\beta}(t,x) + \beta \right] \left(\lambda w^{-1}(t,x) - 1 \right) \ge 2.
\]
\begin{lem}
We have $\partial_x^2\undu(t,x)\ge 0$ for all $t\ge 1$  and $x$ in $\mathbb{R}$ such that $w(t,x)\le\frac{\lambda}{\delta_c}$, where
\[
\delta_c := 1+ \frac{2}{\min \left(\beta , 1+\frac{1}{2s} \right)}.
\]
\end{lem}
\begin{proof}
First, recall that, for all $t\ge 1$ and $x$ in $\mathbb{R}$, $U(t,x) \geq 1$, so that $0 \leq (U(t,x))^{1-\beta}\leq 1$. Assume that $1+\frac{1}{2s}-\beta\ge 0$, then, if $w(t,x) \leq \frac{\lambda}{1+ \frac{2}{\beta}}$, one has from the above inequality
\[
\left[ \left(1+\frac{1}{2s} - \beta \right) U^{1-\beta}(t,x) + \beta \right] \left(\lambda w^{-1}(t,x) - 1 \right) \geq \beta \left(\lambda w^{-1}(t,x) - 1 \right) \geq 2,
\]
and so, for all $t\ge 1,$ $\partial_x^2\undu(t,x)\ge 0$ if $w(t,x) \leq \frac{\lambda}{1+\frac{2}{\beta}}$.\\
Assume next that $1+\frac{1}{2s} - \beta<0$. If $w(t,x) \leq \frac{\lambda}{1+ \frac{2}{1+\frac{1}{2s}}}$, one has
\begin{align*}
\left[ \left(1+\tfrac{1}{2s} - \beta \right) U^{1-\beta}(t,x) + \beta \right] \left(\lambda w^{-1}(t,x) - 1 \right)&\geq \left(1+\tfrac{1}{2s}\right) \left(\lambda w^{-1}(t,x) - 1 \right)\\&\geq 2.
\end{align*}
\end{proof}

\begin{prop}\label{bcl2-prop:w2}
For all $t\ge 1$ and $x\ge 2^{\frac{1}{2s(\beta -1)}}X(t)$, one has
\[
w(t,x)\le \frac{2^{\frac{1}{\beta -1}}\kappa t}{x^{2s}}.
\]
\end{prop}
\begin{proof}
Using definition \eqref{bcl2-def:w}, we have, for all $t\ge 1$ and $x$ in $\mathbb{R}$,
\[
w(t,x)= \frac{\kappa t }{x^{2s}} \left(1-\frac{\gamma (\beta-1) t^{\beta} \kappa^{\beta -1}}{x^{2s(\beta -1)}} \right)^{-\frac{1}{\beta -1}}.
\]
By using the expression \eqref{bcl2-def:X} for $X(t)$, it follows that, for all $t\ge 1$ and $x\ge 2^{\frac{1}{2s(\beta -1)}}X(t)$,
\begin{align*}
w(t,x)&\le \frac{\kappa t }{x^{2s}} \left(1-\frac{\gamma (\beta-1) t^{\beta} \kappa^{\beta -1}}{2X(t)^{2s(\beta -1)}} \right)^{-\frac{1}{\beta -1}},\\
&\le \frac{\kappa t }{x^{2s}} \left(1-\frac{\gamma (\beta-1) t^{\beta} \kappa^{\beta -1}}{2(\kappa t)^{\beta -1}\left[ \lambda^{1-\beta} +\gamma(\beta -1)t\right]} \right)^{-\frac{1}{\beta -1}}, \\
&\le \frac{\kappa t }{x^{2s}} \left(1-\frac{\gamma (\beta-1) t}{2\left[ \lambda^{1-\beta} +\gamma(\beta -1)t\right]} \right)^{-\frac{1}{\beta -1}} \le \frac{2^{\frac{1}{\beta -1}} \kappa t }{x^{2s}}.
\end{align*}
\end{proof}

Finally, let us observe that, for all $t\ge 1$, $X(t)$ satisfies
\begin{align}
&\frac{t}{X(t)}\le \left(\frac{1}{\kappa^{\frac{1}{2s}}(\gamma (\beta-1))^{\frac{1}{2s(\beta -1)}}}\right) t^{1 -\frac{\beta}{2s(\beta -1)}} \Longrightarrow \lim_{t\to+\infty} \frac{t}{X(t)} = 0,\label{bcl2-eq:esti-t/X(t)}\\
&\frac{\kappa t}{X^{2s}(t)}=\frac{\lambda}{ (1+ \lambda^{\beta -1} \gamma (\beta-1)t )^{\frac{1}{\beta -1}}} \Longrightarrow \lim_{t\to+\infty} \frac{\kappa t}{X^{2s}(t)}=0. \label{bcl2-eq:esti-t/X2s(t)}
\end{align}
From the above estimates, we can also derive the following useful limits
\begin{align}
&\lim_{t\to+\infty} \frac{t \ln{t}}{X(t)} = 0,\label{bcl2-eq:esti-tlnt/X(t)}\\
&\lim_{t\to+\infty}\partial_xw(t,X(t)) = 0.\label{bcl2-eq:limw_x(t,X(t))}
\end{align} 
The second assertion is based on the fact that, using the definition of $X$, one can deduce that, for all $t\geq1$,
\begin{align}
\partial_xw(t,X(t))&=\frac{-2s\lambda^\beta}{(\kappa t)^{\beta-1}}(X(t))^{2s(\beta-1)-1}\nonumber\\&=-2s\lambda \left(\frac{\lambda}{\kappa}\right)^{\frac{1}{2s}}  \left(\frac{1}{t}\ +\lambda^{\beta-1}\gamma(\beta-1)\right)^{1-\frac{1}{2s(\beta-1)}} t^{1-\frac{\beta}{2s(\beta -1)}}\label{expression for w_x(X)}.
\end{align}

\section{Proof of \Cref{thm:main} for $s \geq 1$}\label{sec:estimate1}
\subsection{Choice of parameters and consequences}
Let us define $t_\lambda:=\frac{\sigma}{\lambda}$ for some $\sigma>1$ and let us show that for a right choice of the previously introduced parameters $\lambda$, $\kappa$, $\sigma$ and $\gamma$, the function $\undu(\cdot+t_\lambda-1,\cdot)$, with $\undu$ defined in \eqref{bcl2-def:undu}, satisfies \eqref{bcl2-eq:subsol} for all $t\ge 1$. Note that, by definition, we have $t_\lambda>1$ for all $\lambda$ in $(0,1]$.

In the remainder of the present section, let us set 
\[
\kappa:=\frac{D^{2s}}{2} \frac{\lambda}{\sigma}, \qquad \gamma:= \frac{\lambda^{2-\beta}}{\beta-1},
\]
the positive constant $D$ being given in Proposition \ref{prop:initialtail}. Let us also respectively define the functions $X_c$ and $Y$ by
\begin{equation}\label{def:Xc}
X_c(t) := (\kappa t)^{\frac{1}{2s}}\left[ \delta_c^{\beta -1}\lambda^{1-\beta} +\gamma(\beta -1)t\right]^{\frac{1}{2s(\beta-1)}},
\end{equation}
\begin{equation}\label{def:Y}
Y(t) :=(\kappa t)^{\frac{1}{2s}}\left[ \left(2\delta_c \right)^{\beta -1}\lambda^{1-\beta} +\gamma(\beta -1)t\right]^{\frac{1}{2s(\beta-1)}}.
\end{equation}
These are such that $X_{c}(t) < Y(t)$, $w(t,X_c(t))=\frac{\lambda}{\delta_c}$ and $w(t,Y(t))=\frac{\lambda}{2\delta_c}$. As a consequence, one has, for all $t\geq t_\lambda$,
\[
\frac{Y(t) - X_c(t)}{(\kappa t)^{\frac{1}{2s}}} =  \left[ \left(2\delta_c \right)^{\beta -1}\lambda^{1-\beta} +\gamma(\beta -1)t\right]^{\frac{1}{2s(\beta-1)}}- \left[ \delta_c^{\beta -1}\lambda^{1-\beta} +\gamma(\beta -1)t\right]^{\frac{1}{2s(\beta-1)}},
\]
therefore 
\[Y(t) - X_c(t) \geq (\kappa t)^{\frac{1}{2s}} \frac{ 2^{\beta -1} - 1}{2s(\beta-1)} \delta_c^{\beta -1} \lambda^{1-\beta} \left[ \zeta_{s,\beta} \delta_c^{\beta -1}\lambda^{1-\beta} +\gamma(\beta -1)t\right]^{\frac{1}{2s(\beta-1)}-1}
\]
where 
\begin{equation*}
\zeta_{s,\beta} =
\begin{cases}
1 &\text{if } 2s(\beta -1) < 1,\\
2^{\beta -1} &\text{if } 2s(\beta -1) > 1.
\end{cases}
\end{equation*}
The latter being increasing in $t$ in both configurations, we obtain, using the values of $\kappa$ and $\gamma$, for all $t\geq t_\lambda$,
\begin{align*}
Y(t) - X_c(t)&\geq (\kappa t_\lambda)^{\frac{1}{2s}} \tfrac{ 2^{\beta -1} - 1}{2s(\beta-1)} \delta_c^{\beta -1} \lambda^{1-\beta} \left[ \zeta_{s,\beta} \delta_c^{\beta -1}\lambda^{1-\beta} +\gamma(\beta -1)t_\lambda\right]^{\frac{1}{2s(\beta-1)}-1} \\
&\geq\tfrac{D}{2^{\frac{1}{2s}}}\tfrac{ 2^{\beta -1} - 1}{2s(\beta-1)} \delta_c^{\beta -1} \lambda^{1-\beta} \lambda^{(1-\beta)\left(\frac{1}{2s(\beta-1)}-1\right)} \left[ \zeta_{s,\beta} \delta_c^{\beta -1}+\sigma\right]^{\frac{1}{2s(\beta-1)}-1}\\
&\geq\tfrac{D}{2^{\frac{1}{2s}}} \tfrac{ 2^{\beta -1} - 1}{2s(\beta-1)} \delta_c^{\beta -1} \lambda^{- \frac{1}{2s}} \left[ \zeta_{s,\beta} \delta_c^{\beta -1}+ \sigma \right]^{\frac{1}{2s(\beta-1)}-1}\\
&=: \mathcal{C}_1 \lambda^{-\frac{1}{2s}}.
\end{align*}
 
We end this section with a useful computation for further use. Since $w$ is decreasing and convex with respect to $x$ for $x>0$, one has, for all $t\geq t_\lambda$,
\begin{equation}\label{bcl:eq-esti-w_x-carre}
(\partial_xw(t,X(t)))^2\le (\partial_xw(t_\lambda,X(t_\lambda)))^2=\frac{4s^2\left(1 + \sigma\right)^{2-\frac{1}{s(\beta-1)}}}{D^2}  \lambda^{2\left( 1+\frac{1}{2s}\right)}.
\end{equation}

\subsection{Estimating $\mathcal{D}[\underline{u}](t,x)$ where $x \leq X(t)$}
In this region, by the definition of $\undu$, we have, for all $t\geq1$,
\begin{equation*}
\opd{\undu}(t,x) = \int_{y \geq X(t)} [\undu(t,y) - \lambda] J(x-y) \, dy,
\end{equation*}
and we aim at showing \eqref{eq:subsolleft}. For the convenience of the reader, we shall state the following result.

\begin{prop}\label{bcl2-cla:esti-frac1}
For any $\sigma>1$, there exists $\lambda_0(\sigma)$ such that, for any $\lambda \le \lambda_0(\sigma)$, we have, for all $t\geq t_\lambda$ and $x\leq X(t)$,
\[
\opd{\undu}(t,x) + \frac{\lambda^{\beta}}{2} (1-\lambda) \ge 0.
\]
\end{prop}
\begin{proof}
Let us subdivide the interval $(-\infty,X(t)]$ into $(-\infty,X(t)-B]$ and $(X(t)-B,X(t)]$, with $B>1$ to be chosen later, and estimate $\opd{\undu}(t,x)$ on these two subsets.

\paragraph*{When $ x \leq X(t)-B$:} in this subset, Hypothesis \ref{hyp:J} and a short computation gives, for all $t\geq1$,
\begin{align*}
\opd{\undu}(t,x)&= \int_{X(t)}^{+\infty} \frac{\undu(t,y) - \lambda}{ \vert x-y\vert^{1+2s}} J(x-y)\vert x-y\vert^{1+2s} \, dy\\&\geq - \lambda \mathcal{J}_0\int_{X(t)}^{+\infty} \frac{dy}{(y-x)^{1+2s}}\\
&\geq - \frac{\lambda\mathcal{J}_0}{2s} \frac{1}{(X(t)-x)^{2s}}\\
&\geq - \frac{\lambda\mathcal{J}_0}{2s} \frac{1}{B^{2s}}.
\end{align*}

\paragraph*{When $X(t)-B< x \le X(t)$:} in this subset, by making the change of variable $z=y-x$, a short computation gives, for all $t\geq1$,
\begin{align*}
\opd{\undu}(t,x)&= \int_{X(t)-x}^{X(t)-x + B} [\undu(t,x+z) - \lambda] J(z) \, dz+\int_{X(t)-x + B}^{+\infty} [\undu(t,x+z) - \lambda]J(z)\, dz\\
 &\geq \int_{X(t)-x}^{X(t)-x + B}[\undu(t,x+z) - \undu(t,x)]J(z)\,dz- \lambda\mathcal{J}_0 \int_{X(t)-x + B}^{+\infty} \frac{dz}{z^{1+2s}}\\
&= \int_{X(t)-x}^{X(t)-x + B} [\undu(t,x+z) - \undu(t,x)]J(z) \, dy- \frac{\lambda\mathcal{J}_0}{2s} \frac{1}{(X(t) +B -x)^{2s}}\\
&\geq  \int_{X(t)-x}^{X(t)-x + B} [\undu(t,x+z) - \undu(t,x)]J(z) \, dy- \frac{\lambda\mathcal{J}_0}{2s} \frac{1}{B^{2s}},
\end{align*}
since $B>1$ and $\undu(t,x)=\lambda$. By using Taylor's theorem with integral form of the remainder, we have, for all $t\geq1$ and $z$ in $(X(t)-x,X(t)-x + B)$,
\[
\undu(t,x+z) - \undu(t,x)=z\int_{0}^1\partial_x\undu(t,x+\tau z)\,d\tau,
\]
and thus we can estimate the remaining integral by 
\begin{align*}
I&:=\int^{X(t)-x + B}_{X(t)-x} [\undu(t,x+z) - \undu(t,x)]J(z) \, dz\\&= \int^{X(t)-x + B}_{X(t)-x}\int_{0}^1 \partial_x\undu(t,x+\tau z)zJ(z)\,d\tau\,dz.
\end{align*} 
Since $\partial_x\undu$ is a $\mathscr{C}^1$-function with respect to $x$, and $\partial_x\undu(t,x)=0$ for all $x\le X(t)$, we can again apply Taylor's theorem in order to rewrite the last integral as 
\[
I=\int_{0}^{1}\int_{0}^1\int_{X(t)-x}^{B+X(t)-x}\partial_x^2\undu(t,x+\tau \omega z) J(z)\tau z^2 \, dzd\tau d\omega.
\]
Since $\partial_x^2\undu(t,x)=0$ for all $x\le X(t)$ (see \eqref{bcl2-def:u_tx}), the integral can be rewritten as
\[
I=\int_{0}^{1}\int_{0}^1\int_{0}^{B+X(t)-x}\partial_x^2\undu(t,x+\tau \omega z) J(z)\tau z^2 \, dzd\tau d\omega.
\]
From expression \eqref{bcl2-def:u_xx} and the convexity of $w$, we get
\begin{align*}
I&\ge \min_{X(t)\le \xi\le B+X(t)-x}\partial_x^2\undu(t,\xi)\int_{0}^1\int_{0}^{1}\int_{0}^{B+X(t)-x} J(z)\tau z^2 \, dzd\tau d\omega,\\
I&\ge -\frac{6}{\lambda}(\partial_xw(t,X(t)))^2 \int_{0}^1\int_{0}^{1}\int_{0}^{2B} J(z)\tau z^2 \, dzd\tau d\omega\\
&\ge -\frac{6}{\lambda}(\partial_xw(t,X(t)))^2 \left(\int_{0}^1\int_{0}^1\int_{0}^1J(z)\tau z^2\,d\tau d\omega dz +\int_{1}^{2B}\int_{0}^1\int_{0}^1 J(z)\tau z^2 \,d\tau d\omega dz \right)\\
&\ge -\frac{3}{\lambda}(\partial_xw(t,X(t)))^2 \left(\mathcal{J}_1 + \mathcal{J}_0 \int_{1}^{2B} z^{1-2s} \, dz \right),
\end{align*}
using again Hypothesis \ref{hyp:J}. As a consequence, we obtain the following estimate, for all $t\geq1$,
\begin{equation}\label{bcl-eq:esti-1}
\opd{\undu}(t,x)\geq - \frac{\lambda\mathcal{J}_0}{2s} \frac{1}{B^{2s}} - \frac{3}{\lambda} \left(\mathcal{J}_1+\mathcal{J}_0 \int_{1}^{2B} z^{1-2s} \, dz\right)(\partial_xw(t,X(t)))^2.
\end{equation}
Setting $B:=\left(\frac{2\mathcal{J}_0}{s\lambda^{\beta-1}(1-\lambda)}+1\right)^{\frac{1}{2s}}$ then implies that one has, in both cases,
\begin{align*}
\opd{\undu}(t,x) &\ge - \frac{\lambda^{\beta}(1-\lambda)}{4} - \frac{3}{\lambda}\left(\mathcal{J}_1 + \mathcal{J}_0 \int_{1}^{2B} z^{1-2s} \, dz\right)(\partial_xw(t,X(t)))^2\\
&\ge - \frac{\lambda^{\beta}(1-\lambda)}{4} - \frac{3}{\lambda}\left(\mathcal{J}_1 + \mathcal{J}_0 \ln(2B)\right)(\partial_xw(t,X(t)))^2,
\end{align*}
which, using estimate \eqref{bcl:eq-esti-w_x-carre}, leads to, for all $t\geq t_\lambda$, 
\[
\opd{\undu}(t,x)\ge - \left[ \frac{1-\lambda}{4} + \frac{12 s^2 \left[ 1 + \sigma\right]^{2-\frac{2}{(\beta-1)2s}}}{D^2} \left(\mathcal{J}_1 + \mathcal{J}_0 \ln(2B)\right)\lambda^{1-\beta+\frac{1}{s}}\right] \lambda^{\beta}.
\]

Recall that the relation \eqref{s-beta rel} between $s$ and $\beta$ holds, so that we have 
$\beta <1+\frac{1}{2s-1}$ and $2s(\beta -1)<\beta$. 
Then  the fact that $s\ge 1$ implies that $\beta<2$ and that $1-\beta+\frac{1}{s}=\frac{2 - 2s(\beta-1)}{2s}\ge \frac{2-\beta}{2s} >0$. As a consequence,  $\lambda^{1-\beta+\frac{1}{s}}\left( \mathcal{J}_1 + \mathcal{J}_0\ln{(2B)}\right)\le \lambda^{\frac{2 -\beta}{2s}}\left( \mathcal{J}_1 + \mathcal{J}_0\ln{(2B)}\right)$ tends to $0$ with $\lambda$, and there exists a positive real number $\lambda_0$ (explicitly depending on $\sigma$) such that, for all $\lambda\le \lambda_0$ and for all $t\geq t_\lambda$ and $x\leq X(t)$,
\[
\opd{\undu}(t,x) + \frac{\lambda^{\beta}}{2}(1-\lambda)\ge 0. 
\]
\end{proof}

\subsection{Estimating $\mathcal{D}[\underline{u}](t,x)$ where $x>X(t)$}
As explained earlier and shown in \Cref{fig:zones}, we shall estimate $\opd{\undu}(t,x)$ differently in the three distinct intervals $[X(t), Y(t)]$, $[Y(t), 2^{\frac{1}{2s(\beta-1)}}X(t)]$ and $[2^{\frac{1}{2s(\beta-1)}}X(t),+\infty)$. Recall that the exact expression of $Y(t)$ is explicit and is such that $w(t,Y(t))=\frac{\lambda}{2\delta_c}$. Note also that, by definition, $Y>X_{c}$ and that $Y(t)\geq X(t) + R_0$ for $t\geq t_\lambda$, if $\lambda$ is small enough.

\subsubsection{The region $X(t)< x\le Y(t)$}
Let us begin with a technical estimate.

\begin{lem}\label{bcl-lem:reg-2}
For any $B>1$,  one has, for all $t\ge 1$ and $x> X(t)$,
\[
\opd{\undu} (t,x) \ge -\frac{\mathcal{J}_0\lambda}{s B^{2s}}-\frac{6}{\lambda}\left(\mathcal{J}_1+\mathcal{J}_0 \int_{1}^{B} z^{1-2s}\, dz\right) \left(\partial_xw(t,X(t))\right)^2.
\]
\end{lem}

\begin{proof}
Let $\delta>0$ be a parameter with value be fixed later. Since $J$ is nonnegative and $\undu$ is monotone nonincreasing in $x$, we have, for any $\delta \ge R_0$ and for all $t\ge 1$ and $x$ in $\mathbb{R}$, 
\begin{align*}
\opd{\undu}(t,x)&= \int_{-\infty}^{X(t)-x-\delta}[\lambda - \undu(t,x)]J(z) \, dz + \int_{X(t)-x-\delta}^{+\infty}[\undu(t,x+z) - \undu(t,x)] J(z) \, dz\\
&\geq \int_{X(t)-x-\delta}^{+\infty}[\undu(t,x+z) - \undu(t,x)]J(z) \, dz,
\end{align*}
which leads to
\begin{equation}\label{bcl-def:frac}
\opd{\undu}(t,x)\geq \int_{x+z \geq X(t)-\delta, \vert z\vert \leq B} [\undu(t,x+z) - \undu(t,x)]J(z) \, dz + \int_{x+z \geq X(t)-\delta,\vert z\vert \geq B} [\undu(t,x+z) - \undu(t,x)] J(z) \, dz.
\end{equation}

The second integral in the right-hand side of the above equality is the easiest to deal with. Since $\undu $ is positive and $J$ satisfies Hypothesis \ref{hyp:J}, we have, for $B>1$, 
\[
\int_{x+z \geq X(t)-\delta, \vert z\vert \geq B} \hspace{-1cm}[\undu(t,x+z) - \undu(t,x)]J(z) \, dz \geq - \undu(t,x)\mathcal{J}_0 \int_{x+z \geq X(t)-\delta, \vert z\vert \geq B} \frac{dz}{ \vert z\vert^{1+2s}}.
\]
When $X(t)-\delta\le x-B$ and $x > X(t)$, a short computation shows that 
\begin{align*}
\int_{x+z \geq X(t)-\delta, \vert z\vert \geq B} \frac{dz}{ \vert z\vert^{1+2s}} &= \int_{X(t)-x-\delta}^{-B} \frac{dz}{ \vert z\vert^{1+2s}} + \int_{B}^{+\infty} \frac{dz}{ \vert z\vert^{1+2s}}\\
&= \int_{X(t)-x-\delta}^{-B} \frac{dz}{ z^{1+2s}} + \int_{B}^{+\infty} \frac{dz}{ z^{1+2s}}\\
&= \frac{1}{2sB^{2s}} - \frac{1}{2s(x+\delta -X(t))^{2s}} +\frac{1}{2sB^{2s}}.
\end{align*}
On the other hand, if $ X(t)-x-\delta\ge -B$, one has 
\[
\int_{x+z \geq X(t)-\delta, \vert z\vert \geq B} \frac{dz}{ \vert z\vert^{1+2s}} = \int_{B}^{+\infty}\frac{dz}{ z^{1+2s}}= \frac{1}{2sB^{2s}}.
\]
In each situation, we have, for all $t\ge 1$ and $x>X(t)$,
\[
\int_{x+z \geq X(t)-\delta, \vert z\vert \geq B} [\undu(t,x+z) - \undu(t,x)] J(z) \, dz \geq - \frac{\undu(t,x)\mathcal{J}_0}{sB^{2s}} \geq - \frac{\mathcal{J}_0 \lambda}{sB^{2s}}.
\]

Let us now estimate the first integral of the right-hand side of \eqref{bcl-def:frac}, that is
\[
I:=\int_{x+z \geq X(t)-\delta, \vert z\vert \leq B} [\undu(t,x+z) - \undu(t,x)]J(z) \, dz.
\]
As in the  proof of Proposition \ref{bcl2-cla:esti-frac1}, since the function $\undu(t,\cdot)$ is of class $\mathscr{C}^{1}$ for all $t \geq 1$, by using Taylor's theorem with integral form of the remainder, we have, for all $t\ge 1$ and $x$ in $\mathbb{R}$,
\[
I=\int_{x+z \geq X(t)-\delta, \vert z\vert \leq B}\int_{0}^1\partial_x\undu(t,x+\tau z)z J(z)\,d\tau dz.
\]
Now observe that for $\delta\ge B$, we have $-B \geq X(t)-x-\delta$ and therefore 
\[
\int_{x+z \geq X(t)-\delta, \vert z\vert \leq B}\int_{0}^1z J(z)\,d\tau dz= \int_{\vert z\vert \leq B}\int_{0}^1 z J(z)\,d\tau dz=0.
\]
As a consequence, since the partial derivative $\partial_x\undu$ is of class $\mathscr{C}^1$ with respect to $x$,  we can again apply Taylor's theorem in order to rewrite the integral $I$ as
\[
I=\int_{x+z \geq X(t)-\delta, \vert z\vert \leq B}\int_{0}^1\int_{0}^1\partial_x^2\undu(t,x+\omega\tau z)\tau z^2 J(z)\,d\tau d\omega dz,
\]
which leads to the estimate
\begin{align*}
I&\ge \min_{-B<\xi<B}\partial_x^2\undu(t,x+\xi)\left(\int_{\vert z\vert \leq B}\int_{0}^1\int_{0}^1 \tau z^2 J(z)\,d\tau d\omega dz \right),\\
&\ge \min_{-B<\xi<B}\partial_x^2\undu(t,x+\xi)\left(\int_{\vert z\vert\le 1}\int_{0}^1\int_{0}^1 \tau z^2 J(z)\,d\tau d\omega dz +\int_{1\le \vert z\vert \leq B}\int_{0}^1\int_{0}^1 \tau z^2 J(z)\,d\tau d\omega dz\right).
\end{align*}
By using the convexity of $w$ and properties \eqref{bcl2-def:u_tx} and \eqref{bcl2-def:u_xx}, we deduce that 
\begin{align*}
I &\ge -\frac{6}{\lambda} \left(\mathcal{J}_1 + \mathcal{J}_0 \int_{1}^{B} z^{1-2s}\, dz\right) \sup_{\stackrel{-B<\xi<B,}{x+\xi>X(t)}} \partial_xw(t,x+\xi)^2\\
&\ge -\frac{6}{\lambda} \left(\mathcal{J}_1 + \mathcal{J}_0 \int_{1}^{B} z^{1-2s}\, dz\right)\partial_xw(t,X(t))^2.\label{bcl2-eq:esti4} 
\end{align*}
By taking $\delta\ge \sup\{R_0,B\}$, gathering the previous results yields the expected estimate.
\end{proof}

With this lemma at hand, we show the following result. 

\begin{prop}\label{bcl2-cla:esti-frac2}
For any $\sigma>1$, there exists $\lambda_1(\sigma)$ such that, for all $\lambda \le \lambda_1(\sigma)$, we have, for all $t\geq t_\lambda$ and $X(t)<x\leq Y(t)$,
\[
\opd{\undu}(t,x)+ \frac12 (1-\lambda) \undu^{\beta}(t,x)\ge 0.
\]
\end{prop}
\begin{proof}
Set $B:=\nu\lambda^{\frac{1-\beta}{2s}}$, with $\nu>1$ to be chosen later. Note that $B>1$ since $\nu>1$, $\beta>1$ and $\lambda\le 1$. With this choice, the inequality from Lemma \ref{bcl-lem:reg-2} reads, for all $t\geq1$ and $x>X(t)$,
\[
\opd{\undu} (t,x) \ge -\frac{\mathcal{J}_0\lambda^{\beta}}{s \nu^{2s}}-\frac{6}{\lambda} \left(\mathcal{J}_1 +\mathcal{J}_0 \int_{1}^{\nu\lambda^{\frac{1-\beta}{2s}}} z^{1-2s}\, dz \right) \left(\partial_xw(t,X(t))\right)^2.
\]

Next, observe that it follows from the explicit form of $\undu$ (see \eqref{bcl2-def:undu}) that $3w(t,x)\ge \undu(t,x)\ge w(t,x)$ for all $t\ge 1$ and $x>X(t)$. Since we also have $w(t,x)\ge \frac{\lambda}{2\delta_c}$ for all $t\geq1$ and $x\le Y(t)$, we get $\undu(t,x)\ge \frac{\lambda}{2\delta_c}$. As a consequence, one has for any $t\ge 1$ and $X(t)<x\le Y(t)$,
\[
\opd{\undu}(t,x)+\frac{(1-\lambda)}{2}\undu^{\beta}(t,x) \ge 
-\frac{ {6}}{\lambda}\left(\mathcal{J}_1 +\mathcal{J}_0 \int_{1}^{\nu\lambda^{\frac{1-\beta}{2s}}}\hspace{-0.5cm}  z^{1-2s}\, dz\right) \left(\partial_xw(t,X(t))\right)^2+\lambda^\beta\frac{1-\lambda}{2 (2\delta_c)^\beta} -\frac{\mathcal{J}_0\lambda^{\beta}}{s \nu^{2s}},
\]
and thus,
\[
\opd{\undu}(t,x)+\frac{(1-\lambda)}{2}\undu^{\beta}(t,x) \ge 
 - \frac{6}{\lambda}\left(\mathcal{J}_1 +\mathcal{J}_0 \int_{1}^{\nu\lambda^{\frac{1-\beta}{2s}}}\hspace{-0.5cm} z^{1-2s}\, dz\right)(\partial_xw(t,X(t)))^2+\lambda^\beta\frac{1-\lambda}{4(2\delta_c)^\beta} \\,
\]
where we have set $\nu:=\max\left\{1,\left(\frac{4 (2\delta_c)^\beta\mathcal{J}_0}{s(1-\lambda)}\right)^{\frac{1}{2s}}\right\}$. We may now reproduce the argument used in the proof of Proposition \ref{bcl2-cla:esti-frac1} to find an adequate positive real number $\lambda_1$.
\end{proof}

\subsubsection{A preliminary estimate in the range $x>Y(t)$}

\noindent In this zone, the function $\undu$ is convex with respect to $x$, since $Y(t)\geq X_c(t)$. 
 
\begin{lem}\label{lem:estDfarbis}
There exists a positive constant $\mathcal{C}_0$ such that, for any $\lambda$ in $(0,1)$, one has, for all $t \ge 1$ and $x \ge Y(t)$, and for any $B > R_0$ such that $x - B \ge X_c(t)$,
\begin{equation} \label{bcl2-eq:esti-frac4-ter}
\opd{\undu}(t,x)\ge \frac{\lambda - \undu(t,x)}{2s\mathcal{J}_0x^{2s}} - \frac{\mathcal{C}_0}{B^{2s-1}}\,\frac{x^{2s(\beta -1) -1}}{(\kappa t)^{\beta-1}} w^{\beta}(t,x).
\end{equation}
\end{lem}
\begin{proof}
For $t\ge 1$, let us consider the expression for $\opd{\undu}(t,x)$, that we split into three parts:
\begin{multline*}
\opd{\undu}(t,x) = \int_{-\infty}^{-B} [{\undu}(t,x+z) - \undu(t,x)]J(z) \, dz\\+ \int_{-B}^{B} [\undu(t,x+z) - \undu(t,x)]J(z) \, dz + \int^{+\infty}_{B} [{\undu}(t,x+z) - \undu(t,x)] J(z) \, dz.
\end{multline*}
To obtain an estimate of the second integral, we follow the same steps as previously to obtain
\begin{align*}
\int_{-B}^{B} [\undu(t,x+z) - \undu(t,x)]J(z) \, dz&=\int_{-B}^{B}\int_{0}^1\int_{0}^1\hspace{-0.1cm}\partial_x^2\undu(t,x+\tau \omega z) \tau z^2J(z) \, d\tau d\omega dz\\&\geq 0,
\end{align*}
since $\undu$ is convex with respect to $x$ in the domain of integration. Next, using again Taylor's theorem, the last integral in the decomposition may be rewritten as
\[
\int_{B}^{+\infty} [\undu(t,x+z) - \undu(t,x)]J(z) \, dz=\int_{B}^{+\infty}\int_{0}^1\partial_x\undu(t,x+\tau z) z J(z) \,d\tau dz.
\]
Observe that since $x \geq Y(t)$ and $w$ is decreasing and convex with respect to $x$, identity \eqref{bcl2-def:u_x} implies that, for all positive $\tau$ and $z$, one has
\begin{align*}
\partial_x\undu(t,x+\tau z)&= 3\partial_xw (t,x+\tau z) \left(1- \frac{w(t,x+\tau z)}{\lambda} \right)^2\\&\geq 3 \left(1- \frac{1}{2 \delta_c} \right)^2\partial_xw(t,x). 
\end{align*}
It then follows from Hypothesis \ref{hyp:J} that 
\begin{align*}
\int_{B}^{\infty} [\undu(t,x+z) - \undu(t,x)]J(z) \, dz &\geq 3 \left(1- \frac{1}{2 \delta_c} \right)^2 \left(\int_{B}^{\infty} z J(z) \, dz\right) \partial_xw(t,x)\\
&\geq 3 \left(1- \frac{1}{2 \delta_c} \right)^2 \frac{\mathcal{J}_0}{(2s-1)B^{2s-1}} \partial_xw(t,x).
\end{align*}
Finally, since $X(t)-x \leq X(t) - X_c(t) -B \leq -B$, the first integral can be estimated as follows
\[
\int_{-\infty}^{-B} [\undu(t,x+z) - \undu(t,x)] J(z) \, dz
\ge {\mathcal{J}_0}^{-1}\int_{-\infty}^{X(t)-x} \frac{\undu(t,x+z) - \undu(t,x)}{ \vert z\vert^{1+2s}} \, dz+\int_{X(t)-x}^{-B} [\undu(t,x+z) - \undu(t,x)]J(z) \, dz,
\]
which by taking advantage of the fact that $\undu$ is decreasing with respect to $x$ leads to
\begin{align*}
\int_{-\infty}^{-B} [\undu(t,x+z) - \undu(t,x)] J(z) \, dz&\geq \frac{{\mathcal{J}_0}^{-1}}{2s} \frac{\lambda - \undu(t,x)}{(x-X(t))^{2s}} {\geq\frac{{\mathcal{J}_0}^{-1}}{2s} \frac{\lambda - \undu(t,x)}{x^{2s}}}.
\end{align*}

Finally, collecting these estimates and recalling the expression of $\partial_xw$ in \eqref{bcl2-def:w_x} give the result by setting $\mathcal{C}_0=\frac{6s\mathcal{J}_0}{(2s-1)}\left(1- \frac{1}{2 \delta_c} \right)^2$.
\end{proof}

\subsubsection{The region $Y(t) <  x \le 2^{\frac{1}{2s(\beta -1)}}X(t)$}
Let us now estimate $\opd{\undu}(t,x)$ when $t\ge 1$ and $x > Y(t)$. 

\begin{prop}\label{bcl2-prop:esti-frac1}
For any $\sigma>1$, there exists $\lambda_2(\sigma)$ such that for all $\lambda \le \lambda_2(\sigma)$, we have, for all $t\geq t_\lambda$ and $Y(t)<x\leq 2^{\frac{1}{2s(\beta -1)}}X(t)$,
\[
\opd{\undu}(t,x)+ \frac12 (1-\lambda) \undu^{\beta}(t,x)\ge 0.
\]
\end{prop}
\begin{proof}
Let us recall that $Y(t)$ is such that $w(t,Y(t))=\frac{\lambda}{2\delta_c}$. For $t\ge 1$, consider $x> Y(t)$. As long as $B$ is chosen such that $x - B \geq X_c(t)$, it follows from Lemma \ref{lem:estDfarbis} that
\[
\opd{\undu}(t,x)\ge -\frac{\mathcal{C}_0}{B^{2s-1}}\frac{x^{2s(\beta -1) -1}}{(\kappa t)^{\beta-1}} w^{\beta}(t,x).
\]
The rest of the proof deals with the choice of $B$.\\Since $X(t) < x \le 2^{\frac{1}{2s(\beta -1)}}X(t)$, we have directly 
\begin{align*}
\frac{x^{2s(\beta -1) -1}}{(\kappa t)^{\beta-1}} &\leq \frac{2 X(t)^{2s(\beta -1) -1}}{(\kappa t)^{\beta-1}} = 2 (\kappa t)^{-\frac{1}{2s}}\left[ \lambda^{1-\beta} +\gamma(\beta -1)t\right]^{1 - \frac{1}{2s(\beta-1)}} \\
&\leq {\frac{2^{1+\frac{1}{2s}}}{D}}\left[1+\sigma \right]^{1 - \frac{1}{2s(\beta-1)}} \lambda^{(\beta -1)(\frac{1}{2s(\beta-1)}-1)}.
\end{align*}
Consequently, one has, for all $t\ge 1$ and $Y(t)<x\le 2^{\frac{1}{2s(\beta -1)}}X(t),$
\[
\opd{\undu}(t,x) +\frac{1}{2}(1-\lambda)\undu^\beta(t,x)\ge \left[ \frac{1}{2}(1-\lambda) -2 \frac{\mathcal{C}_0}{D} \left[1+\sigma \right]^{1 - \frac{1}{2s(\beta-1)}} \frac{\lambda^{-(\beta -1)\left(1-\frac{1}{2s(\beta-1)}\right)} }{B^{2s-1}} \right] w^{\beta}(t,x).
\]
Observe that the proof is ended if the quantity in the brackets above is positive. We thus set 
\[
B:=\sup\left\{R_0,\left(\tfrac{4 \mathcal{C}_0}{(1 -\lambda)D}\right)^{\frac{1}{2s -1}}(1+\sigma)^{\frac{1}{2s-1}-\frac{1}{2s(2s-1)(\beta -1)}} \lambda^{-\frac{\beta -1}{2s-1}(1 - \frac{1}{2s(\beta-1)})_+}\right\},
\]
which is an adequate value. Indeed, let us point out that the limitation on the choice of $B$ is due to the fact that we need to ensure that $x - X_c(t) - B \geq 0$ for all $t \geq t_\lambda$ and $x> Y(t)$. Since $Y(t) - X_c(t) \geq \mathcal{C}_1 \lambda^{-\frac{1}{2s}}$, this is satisfied as long as $B \leq \mathcal{C}_1 \lambda^{-\frac{1}{2s}}$. Since $\frac{1}{2s}-\frac{\beta -1}{2s-1}(1 - \frac{1}{2s(\beta-1)}) = \frac{2-\beta}{2s-1} > 0$, one may observe that the condition is satisfied by taking $\lambda$ sufficiently small for any choice of $\sigma$.
\end{proof}

\subsubsection{The region $x>2^{\frac{1}{2s(\beta -1)}}X(t)$}
In this region, we claim

\begin{prop}\label{bcl2-prop:esti-frac3}
There exists $\sigma_{LR}>1$ such that, for all $\sigma\ge \sigma_{LR}$ and $\lambda$ in $(0,1)$, we have, for all $t\ge t_\lambda$ and $x> 2^{\frac{1}{2s(\beta -1)}}X(t)$,
\[
\opd{\undu}(t,x) \ge \frac{\lambda(1-\tau_0)}{4\mathcal{J}_0sx^{2s}}.
\]
with $\tau_0:=\frac{3}{2\delta_c} (1-\frac{1}{2\delta_c}+\frac{1}{3}\frac{1}{(2\delta_c)^2})$.  
\end{prop}

\begin{proof}
For a given $x$, let us define $B=\frac{x}{K}$  with $K=\frac{2^{\frac{1}{2s(\beta -1)}+1}}{2^{\frac{1}{2s(\beta -1)}}-1}$. Note that by the definition of $K$, we have $2^{\frac{1}{2s(\beta -1)}}\left(1-\frac{1}{K}\right)=\frac{2^{\frac{1}{2s(\beta -1)}}+1}{2}:=\mathcal{C}_2$. As a consequence, for all $t\ge 1$, one has $x-B\ge \mathcal{C}_2X(t)$. A straightforward computation then shows that, for all $\lambda$ in $(0,1)$ and $t\ge t_\lambda$, one has $\mathcal{C}_2X(t)\ge X_{c}(t)$ as soon as $\mathcal{C}_2X(t_\lambda)\ge X_{c}(t_\lambda)$, that is
\[
{\mathcal{C}_2}^{2s(\beta -1)}\ge \frac{\delta_c^{\beta-1}+ \sigma}{1+\sigma}.
\] 
Since $\mathcal{C}_2>1$ and $\lim_{z\to \infty}\frac{\delta_c^{\beta-1}+z}{1+z}=1$, the above inequality is always true for a large $\sigma$. We may apply Lemma \ref{lem:estDfarbis} to get, for all $\lambda$ in $(0,1)$, $t\ge t_\lambda$ and $x>2^{\frac{1}{2s(\beta -1)}}X(t)$,
\[
\opd{\undu}(t,x)\ge \frac{\lambda - \undu(t,x)}{2s\mathcal{J}_0x^{2s}} -\frac{\mathcal{C}_0}{B^{2s-1}}\frac{x^{2s(\beta -1) -1}}{(\kappa t)^{\beta-1}}w^{\beta}(t,x).
\]
Let us recall that $Y(t)$ is such that $w(t,Y(t))=\frac{\lambda}{2\delta_c}$, and thus $\undu(t,Y(t))= \frac{3\lambda}{2\delta_c} (1-\frac{1}{2\delta_c}+\frac{1}{3}\frac{1}{(2\delta_c)^2}) := \lambda\tau_0$. We can easily check that since $\frac{1}{2\delta_c}<1$, we have $\tau_0<1$ (one may see it as a level set for the function $\undu$ with $\lambda=1$). Therefore, using Proposition \ref{bcl2-prop:w2}, it follows that, for all $t\ge t_\lambda $ and $x> 2^{\frac{1}{2s(\beta -1)}}X(t)$,
\begin{align*}
\opd{\undu}(t,x)&\ge \frac{\lambda(1-\tau_0)}{2\mathcal{J}_0sx^{2s}} -\frac{\mathcal{C}_0}{B^{2s-1}}w^{\beta}(t,x)\frac{x^{2s(\beta -1) -1}}{(\kappa t)^{\beta-1}}\\
&\ge \frac{\lambda(1-\tau_0)}{2\mathcal{J}_0sx^{2s}} -\frac{\mathcal{C}_0 2^{\frac{\beta}{\beta-1}} }{B^{2s-1}}\frac{1}{x^{2s}}\frac{\kappa t}{x}\\
&= \left[ \frac{1-\tau_0}{2\mathcal{J}_0s} - \mathcal{C}_0 2^{\frac{\beta}{\beta-1}} K^{2s -1}\frac{\kappa t}{\lambda x^{2s}} \right] \frac{\lambda}{x^{2s}}\\
&\geq \left[ \frac{1-\tau_0}{2\mathcal{J}_0s} - 2\mathcal{C}_0K^{2s -1}\frac{\kappa t}{\lambda X(t)^{2s}} \right] \frac{\lambda}{x^{2s}}\\
&= \left[ \frac{1-\tau_0}{2\mathcal{J}_0s} - \frac{2\mathcal{C}_0K^{2s -1}}{\left[ 1 +\lambda t\right]^{\frac{1}{\beta-1}}} \right] \frac{\lambda}{x^{2s}}.
\end{align*} 
We finally get, for all $t\ge t_\lambda$ and $x> 2^{\frac{1}{2s(\beta -1)}}X(t)$,
\[
\opd{\undu}(t,x)\ge \left[ \frac{1-\tau_0}{2\mathcal{J}_0s} - \frac{2\mathcal{C}_0K^{2s -1}}{\left[ 1 + \sigma\right]^{\frac{1}{\beta-1}}} \right] \frac{\lambda}{x^{2s}} \ge \frac{\lambda(1-\tau_0)}{4\mathcal{J}_0sx^{2s}},
\]
by choosing $\sigma>1$ large enough. 
\end{proof}

\subsection{Tuning the parameters $\sigma$ and $\lambda$ 
}
In the last part of the proof, we have to choose the parameters $\sigma$ and $\lambda$ in order that, for $t\ge 1,$ the function $\undu(\cdot+t_\lambda-1,\cdot)$ satisfies inequation \eqref{bcl2-eq:subsol}. Recall that it is case if and only if inequations \eqref{eq:subsolleft} and \eqref{eq:subsolright} hold simultaneously for $t\ge t_\lambda$. Since the former holds unconditionally for any $t\ge t_\lambda$ and $\lambda\le \lambda_0$, one only needs to check that the latter holds for a suitable choice of $\sigma>1$.

By using \eqref{bcl2-def:u_t} and \eqref{bcl2-def:w_t} and since $\gamma=\frac{\lambda^{2-\beta}}{\beta -1}$, inequality \eqref{eq:subsolright} holds, in particular, for $t\ge t_\lambda$ and $x>X(t)$, so that one can show that
\[
3\frac{\partial_t\Phi(t,x)}{\Phi^\beta(t,x)}w^{\beta}(t,x)\leq \opd{\undu}(t,x)+ (1-\lambda)\undu^{\beta}(t,x) -\frac{\lambda^{2-\beta}}{\beta -1} w^{\beta}(t,x),
\]
Proving the above inequality is the purpose of the present section. 

For $\sigma>1$, set  $\lambda^*(\sigma):=\inf\{\lambda_0(\sigma),\lambda_1(\sigma),\lambda_2(\sigma)\}$, where $\lambda_0$, $\lambda_1$ and $\lambda_2$ were respectively introduced in Propositions \ref{bcl2-cla:esti-frac1}, \ref{bcl2-cla:esti-frac2}, and \ref{bcl2-prop:esti-frac1}. Let us next decompose the set $[X(t),+\infty)$ into the two subsets
\[
I_1:= [X(t),2^{\frac{1}{2s(\beta -1)}}X(t)]\text{ and }I_2:=(2^{\frac{1}{2s(\beta -1)}}X(t),+\infty).
\]
On the first interval, we have the following result.
\begin{lem}
For all $\sigma>1$, there exists $\lambda_4(\sigma)$ such that, $\forall\,\lambda\le \lambda_4(\sigma)$, one has, for all $t \ge t_\lambda$ and $x$ in $I_1$, 
\[
3\frac{\partial_t\Phi(t,x)}{\Phi^\beta(t,x)}w^{\beta}(t,x)\leq \opd{\undu}(t,x) + (1-\lambda)\undu^{\beta}(t,x) -\frac{\lambda^{2-\beta}}{\beta -1} w^{\beta}(t,x).
\]
\end{lem}
\begin{proof}
By the definition \eqref{eq:phiU} of $\Phi$, one has, for all $t\ge 1$ and $x$ in $\mathbb{R}$,
\[
3\frac{\partial_t\Phi(t,x)}{\Phi^\beta(t,x)}w^{\beta}(t,x) = \frac{3}{t}\frac{x^{2s(\beta -1)}}{(\kappa t)^{\beta-1}}w^{\beta}(t,x).
\]
By exploiting \eqref{bcl2-def:X}, it follows that, for $t\ge 1$ and $x\le 2^{\frac{1}{2s(\beta -1)}}X(t)$, 
\[
3\frac{\partial_t\Phi(t,x)}{\Phi^\beta(t,x)}w^{\beta}(t,x)\le \left(\frac{6}{t}\lambda^{1-\beta}+6\lambda^{2-\beta}\right) w^{\beta}(t,x).
\]
So for $t\ge t_\lambda$ and $x\le 2^{\frac{1}{2s(\beta -1)}}X(t)$ we have, since $\sigma>1$, 
\[
3\frac{\partial_t\Phi(t,x)}{\Phi^\beta(t,x)}w^{\beta}(t,x)\le 6\lambda^{2-\beta}\left[1+ \sigma^{-1}\right] w^{\beta}(t,x)\le 12\lambda^{2-\beta} w^{\beta}(t,x).
\]
From the above, we see that we have, for all $\lambda\le \lambda':=\left(\frac{1}{96}\right)^{\frac{1}{2-\beta}}$, $t\ge t_\lambda$ and $x\le 2^{\frac{1}{2s(\beta -1)}}X(t)$,
\[
3\frac{\partial_t\Phi(t,x)}{\Phi^\beta(t,x)}w^{\beta}(t,x)\le\frac{1}{4} w^{\beta}(t,x).
\]
Recall that by Propositions \ref{bcl2-cla:esti-frac2} and \ref{bcl2-prop:esti-frac1}, we have, for all $\lambda \le \lambda^*$ and for all $t\ge t_\lambda$ and $x$ in $I_1$, 
\[
\opd{\undu}(t,x) + (1-\lambda)\undu^{\beta}(t,x)-\frac{\lambda^{2-\beta}}{\beta -1} w^{\beta}(t,x)\ge \left(\frac{1-\lambda}{2}-\frac{\lambda^{2-\beta}}{\beta-1}\right)w^{\beta}(t,x),
\]
since  $\undu\ge w$ for all $x\ge X(t)$.\\ We end the proof by taking $\lambda\le \lambda_4:=\inf\{\lambda^*(\sigma),\lambda',\lambda''\}$, where $\lambda''$ is such that $\frac{1-\lambda}{2}-\frac{\lambda^{2-\beta}}{\beta-1}\geq\frac{1}{4}$ for all $\lambda\le \lambda''$, which is possible since $\beta < 2$ in that case (see \eqref{s-beta rel}).
\end{proof}

Finally, let us check what happens on $I_2$.
\begin{lem}
There exists  $\sigma^*>1$ such that, for all $\lambda\le \lambda_4(\sigma^*)$, one has, for all $t\ge t_\lambda$ and $x$ in $I_2$, 
\[
3\frac{\partial_t\Phi(t,x) }{\Phi^\beta(t,x) }w^{\beta}(t,x) \leq \opd{\undu}(t,x) + (1-\lambda)\undu^{\beta}(t,x) -\frac{\lambda^{2-\beta}}{\beta -1} w^{\beta}(t,x).
\]
\end{lem}
\begin{proof}
As in the preceding proof, by the definition of $\Phi$, we have, for all $t\ge 1$ and $x$ in $\mathbb{R}$,
\[
3\frac{\partial_t\Phi(t,x)}{\Phi^\beta(t,x)}w^{\beta}(t,x)= 3\kappa\frac{x^{2s(\beta -1)}}{(\kappa t)^{\beta}}w^{\beta}(t,x).
\]
By Proposition \ref{bcl2-prop:w2}, we have, for $t\ge t_\lambda$ and $x$ in $I_2$,
\[
w^{\beta}(t,x)\le 2^{\frac{\beta}{\beta -1}}\frac{(\kappa t)^{\beta}}{x^{2s\beta}}.
\]
Therefore, since we have set $\kappa=\frac{D^{2s}}{2} \frac{\lambda}{\sigma}$, we get, for $t\ge t_\lambda$ and $x$ in $I_2$,
\[
3\frac{\partial_t\Phi(t,x)}{\Phi^\beta(t,x)}w^{\beta}(t,x)\le 3\kappa\frac{x^{2s(\beta -1)}}{(\kappa t)^{\beta}} 2^{\frac{\beta}{\beta-1}}\frac{(\kappa t)^{\beta}}{x^{2s\beta}}= \frac{2^{\frac{\beta}{\beta-1} {-1}}3D^{2s} \lambda}{\sigma x^{2s}}.
\]
 Observe that for all $\sigma \ge \sigma':= \sup\left\{2,D^{2s}2^{\frac{\beta}{\beta-1}} \frac{ {6}\mathcal{J}_0s}{1-\tau_0}\right\}$, so that, for $t\ge t_\lambda$ and $x$ in $I_2$,
\[
3\frac{\partial_t\Phi(t,x)}{\Phi^\beta(t,x)}w^{\beta}(t,x)\le \frac{\lambda (1-\tau_0)}{4\mathcal{J}_0 sx^{2s}}.
\]
Now recall that by Proposition \ref{bcl2-prop:esti-frac3}, we have, for all $\sigma\ge \sigma_{LR}$ and $\lambda$ in $(0,1)$, and for all $t\ge t_\lambda$ and $x$ in $I_2$, 
\begin{multline*}
\opd{\undu}(t,x) + (1-\lambda)\undu^{\beta}(t,x)-\frac{\lambda^{2-\beta}}{\beta -1} w^{\beta}(t,x)\\\ge\left(1-\lambda-\frac{\lambda^{2-\beta}}{\beta -1}\right)w^{\beta}(t,x)+\frac{\lambda(1-\tau_0)}{4\mathcal{J}_0s x^{2s}},
\end{multline*}
since $\undu\ge w$ for all $x\ge X(t)$. The claim is then proved by taking $\sigma^*:=\sup\{\sigma',\sigma_{LR}\}$ and $\lambda\le\lambda_4(\sigma^*)$.
\end{proof}

\subsection{Final argument}
\noindent From the above section, one may find a sufficiently small $\lambda$, say $\lambda\le \lambda_4(\sigma^*)$, such that the function $\undu(\cdot+t_\lambda-1,\cdot)$ satisfies inequation \eqref{bcl2-eq:subsol} for all $t\ge 1$. Having this function at hand, it is enough to check that, for some $R^*$ and $T$, one has $u(T,x+R^*)\ge \undu(t_\lambda,x)$ to conclude the proof. Indeed, if so, by the parabolic comparison principle, one would then have $u(t+T,x+R^*)\ge \undu(t_\lambda+t-1,x)$ for all $t\ge 1$, and the level set
\[
E_\lambda(t):=\{x\in \mathbb{R}|u(t,x)\ge \lambda\}\supset (  - \infty,X(t-T+t_\lambda-1)+R^*],
\]
which implies $x_\lambda(t) \ge X(t+t_\lambda -1-T)+R^*$, and thus $x_\lambda(t) \gtrsim   t^{\frac{\beta}{2s(\beta -1)}}$ as $t$ tends to $+\infty$.

Let us find the adequate $R^*$ and $T$. For $T$, an adequate value is $T=1$ since, by Proposition \ref{prop:initialtail},
\[
\lim_{x\to +\infty}x^{2s}u({1},x) \geq 2 D^{2s}.
\]
On the other hand, by the definition of $\undu$, a quick computation shows that
\[
\lim_{x\to +\infty}x^{2s}\undu(t_\lambda,x)=\frac{3D^{2s}}{2}.
\] 
Therefore, there exists $R_1>0$ such that $u( {1},x)\ge \undu(t_\lambda,x)$ for all $x\ge R_1$. In particular, we have $u(1,x-R_1)\ge \undu(t_\lambda,x)$ for all $x\ge R_1$, since $u(1,\cdot)$ is monotone nonincreasing. To end the proof, we just need to ensure that $\liminf_{x\to -\infty}u( {1},x)>\lambda$. Indeed, if so, then there would exist $R_2>0$ such that $u(1,x)>\undu(t_\lambda,x)$ for all $x<-R_2$. This allows to conclude that $u(1,x-R_1-R_2)\ge\undu(t_\lambda,x)$ since, by monotonicity of $u(1,\cdot)$, we have $u(1,x-R_1-R_2)\ge\lambda\ge \undu(t_\lambda,x)$ for all $x\le R_1$ and $u(1,x-R_1-R_2)\ge u(1,x-R_1)\ge\undu(t_\lambda,x)$ for all $x\ge R_1$.

To prove that $\liminf_{x\to-\infty}u( {1},x)>\lambda$, we just need to observe that, by Proposition \ref{bcl-prop:inva}, we have $\liminf_{x\to -\infty} u(t,x)\ge \frac{a}{4}$ for all $t\ge 0$, and, as a consequence, $\liminf_{x\to-\infty}u( {1},x)>\lambda$  as soon as $\lambda\le \frac{a}{4}$.
As a consequence, $\undu(t+t_\lambda -1,x) \le u(t+1,x-R_1 -R_2)$ for all $t\ge 1$ as soon as $\lambda \le \inf\{\frac{a}{4},\lambda_4(\sigma^*)\}$.

Having a lower bound at hand for $x_\lambda(t)$ when $\lambda\le \inf\{\frac{a}{4},\lambda_4(\sigma^*)\}$, we can obtain one for $x_\lambda(t)$ when $\lambda\ge\inf\{\frac{a}{4},\lambda_4(\sigma^*)\}$ by arguing as in the proof in \cite{Alfaro2017,Coville2021}, using the adequate invasion property, namely Proposition \ref{bcl-prop:inva}, we obtain the desired lower bound for the speed of any level set with value $\lambda$ in $(0,1)$.

\section{Proof of \Cref{thm:main} for $s<1$}\label{sec:estimate}
In this section, we prove \Cref{thm:main} for $s<1$. In such a situation, the above construction based on a fine control of the time $t_\lambda=\frac{\sigma}{\lambda}$, is inadequate for a large set of parameters $(s,\beta)$, especially when $\beta \geq 2$. In this case, the constraint imposed on the form of $t_\lambda$ would cause the proof to fail. To cover all possible situations, new ideas must be developed. When $s<1$, the diffusion process plays an important role in inducing a strong flattening of the solution. So, with this in mind, we exploit the flattening properties of the solution to \eqref{eq:main} to remove the constraint imposed on $t_\lambda$ in the above construction, hoping that we can find a time $t^*$ after which $\undu(\cdot+t^*-1,\cdot)$ satisfies \eqref{bcl2-eq:subsol}. By doing so, we get more flexibility in the construction, but at the expense of a clear understanding of the time after which the true acceleration regime starts. 

We shall show that, for the right choice of $\lambda$, $\kappa$ and $\gamma$, the function $\undu(t+t^*-1,\cdot)$ indeed satisfies \eqref{bcl2-eq:subsol} for all $t\ge 1$ for some $t^*>1$.

\subsection{Estimating $\mathcal{D}[\underline{u}](t,x)$ where $x \leq X(t)$}
In this region, by the definition of $\undu$, we have 
\[
\opd{\undu}(t,x) = \int_{X(t)}^{+\infty} [\undu(t,y) - \lambda] J(x-y) \, dy.
\]
This section aims at showing \eqref{eq:subsolleft}, stated as the following result.

\begin{prop}\label{bcl2-cla:esti-frac1-1}
For all $\lambda\le \frac{1}{2}, \gamma$ and $\kappa$, $ \exists\, t_0(\lambda,\kappa,\gamma,\beta,s) >1$ such that, for all $t\ge t_0$ and $x\le X(t)$,
\[
\opd{\undu}(t,x) + \frac{\lambda^{\beta}}{2} (1-\lambda) \ge 0.
\]
\end{prop}

\begin{proof}
The starting point being the same as in Proposition \ref{bcl2-cla:esti-frac1}, we do not reproduce the beginning of the proof and follow from \eqref{bcl-eq:esti-1}, that is, for all $t\ge 1$ and $X(t)\ge x$,
\[
\opd{\undu(t,x)}\geq - \frac{\lambda\mathcal{J}_0}{2s} \frac{1}{B^{2s}} - \frac{3}{\lambda} \left(\mathcal{J}_1+\mathcal{J}_0 \int_{1}^{2B} z^{1-2s} \, dz\right)\partial_xw(t,X(t))^2.
\]
Choosing $B:=\left(\frac{2\mathcal{J}_0}{s\lambda^{\beta-1}(1-\lambda)}+1\right)^{\frac{1}{2s}}$, we get, for all $t\ge 1$ and $X(t)\ge x$,
\[
\opd{\undu}(t,x)\ge - \frac{\lambda^{\beta}(1-\lambda)}{4} - \frac{3}{\lambda}\left(\mathcal{J}_1 + \mathcal{J}_0 \int_{1}^{2B} z^{1-2s} \, dz\right)(\partial_xw(t,X(t)))^2.
\]
It then follows from \eqref{bcl2-eq:limw_x(t,X(t))} that we can find $t_0>1$ such that, for all $t\ge t_0$,
\[
\frac{3}{\lambda}\left(\mathcal{J}_1 + \mathcal{J}_0 \int_{1}^{2B} z^{1-2s} \, dz\right)(\partial_xw(t,X(t)))^2 \le \frac{\lambda^{\beta}(1-\lambda)}{4},
\]
thus ending the proof.
\end{proof}

\subsection{Estimating $\mathcal{D}[\underline{u}](t,x)$ where $x>X(t)$}
Here again, we estimate $\opd{\undu}(t,x)$ in the three separate intervals 
\begin{equation*}
[X(t), Y(t)], \qquad [Y(t), 2^{\frac{1}{2s(\beta-1)}}X(t)]\quad\text{and }\quad [2^{\frac{1}{2s(\beta-1)}}X(t),+\infty ),
\end{equation*}
recalling that $Y(t)>X_c(t)$ for all $t>0$ such that $w(t,Y(t))=\frac{\lambda}{2\delta_c}$. Note that for all $\lambda$, $\gamma$, $\kappa$, $s$, $\beta$, we may find $t^{\#}>1$ such that $Y(t) \geq X_c(t) + R_0$ for all $t \geq t^\#$.

\subsubsection{The region $X(t)< x\le Y(t)$}
In this region, owing to Lemma \ref{bcl-lem:reg-2}, we claim the following.
\begin{prop}\label{bcl-prop:esti-reg2}
For all $\lambda<\frac{1}{2}$, $\kappa$ and $\gamma$, there exists $t_1>1$ such that, for all $t\ge t_1$ and $X(t)<x\le Y(t)$,
\[
\opd{\undu} \ge -\frac12 (1-\lambda) \undu^{\beta}.
\]
\end{prop}

\begin{proof}
We follow essentially the same steps as in the proof of Proposition \ref{bcl2-cla:esti-frac2}.\\ 
First, using Lemma \ref{bcl-lem:reg-2} with $B:=\eta\lambda^{\frac{1-\beta}{2s}}$, where $\eta:=\sup\left\{\left(\frac{4(2 \delta_c)^{ \beta} \mathcal{J}_0}{s(1-\lambda)}\right)^{\frac{1}{2s}},1\right\}$, and the fact that, for all $x\le Y(t)$, it holds $$\undu(t,x) \geq \undu(t,Y(t)) \geq w(t,Y(t))=\frac{\lambda}{2\delta_c},$$ one has, for $t\ge 1$ and $X(t)<x\le Y(t)$,
\[
\opd{\undu}(t,x)+\frac12 (1-\lambda)\undu^{\beta}(t,x)\ge -\frac{6}{\lambda}\left(\partial_xw(t,X(t))\right)^2\left(\mathcal{J}_1+ \mathcal{J}_0 \int_{1}^{\eta\lambda^{\frac{1-\beta}{2s}}}\hspace{-0.5cm} z^{1-2s}\, dz\right)+\frac14 (1-\lambda) \frac{\lambda^\beta}{(2 \delta_c)^{\beta}}.
\]
From there, we can argue as in the proof of Proposition \ref{bcl2-cla:esti-frac1-1}, using that $\lim_{t\to +\infty}\left(\partial_xw(t,X(t))\right)^2=0$ to find $t_1>1$ such that, for all $t\ge t_1$,
\[
\frac{6}{\lambda}\left(\partial_xw(t,X(t))\right)^2\left(\mathcal{J}_1+ \mathcal{J}_0 \int_{1}^{\eta\lambda^{\frac{1-\beta}{2s}}} z^{1-2s}\, dz\right)\le \frac14 (1-\lambda) \frac{\lambda^\beta}{(2 \delta_c)^{\beta}},
\]
which concludes the proof.
\end{proof}

\subsubsection{A preliminary estimate in the range $x> Y(t)$}
The estimate obtained in Lemma \ref{lem:estDfarbis} is not valid for $s<1$. We thus derive an estimate of $\opd{\undu}(t,x)$ in such case, but only valid in the range $x>Y(t)$.

\begin{lem}\label{lem:estDfar}
For any $B>1$ and for all $t >t^{\#}$ and $x> Y(t)$, one has 
\begin{equation} \label{bcl2-eq:esti-frac4}
\opd{\undu}(t,x)\ge \frac{\lambda(1-\tau_0)}{2s\mathcal{J}_0x^{2s}} -\frac{\mathcal{J}_0\lambda\tau_0}{2sB^{2s}}+3\mathcal{J}_0 \left(\int_{1}^{B} z^{-2s} \, dz\right)\partial_xw(t,x),
\end{equation}
where $\tau_0:=\frac{3}{2\delta_c}\left(1-\frac{1}{2\delta_c}+\frac{1}{3(2\delta_c)^2}\right)$.
\end{lem}
\begin{proof}
For $t\ge 1$, let us split into three parts the integral defining $\opd{\undu}(t,x)$:
\begin{multline*}
\opd{\undu}(t,x) = \int_{-\infty}^{-1}\hspace{-0.1cm} [u(t,x+z) - \undu(t,x)]J(z) \, dz + \int_{-1}^{1}\hspace{-0.1cm} [\undu(t,x+z) - \undu(t,x)]J(z) \, dz \\+ \int^{\infty}_{1} \hspace{-0.1cm}[u(t,x+z) - \undu(t,x)] J(z) \, dz.
\end{multline*}
Since $\undu$ is decreasing for $t\ge t^{\#}$ and $x>Y(t)\ge X_c(t)+R_0\ge X(t)+R_0$, the first integral can be estimated in the following way, using Hypothesis \ref{hyp:J}, 
\begin{align}
\int_{-\infty}^{-1} [\undu(t,x+z) - \undu(t,x)] J(z) \, dz \nonumber&\ge {\mathcal{J}_0}^{-1}\int_{-\infty}^{X(t)-x} \frac{\undu(t,x+z) - \undu(t,x)}{ \vert z\vert^{1+2s}} \, dz+\int_{X(t)-x}^{-1} [\undu(t,x+z) - \undu(t,x)]J(z) \, dz\nonumber\\
&\ge \frac{{\mathcal{J}_0}^{-1}}{2s} \frac{\lambda - \undu(t,x)}{(x-X(t))^{2s}}.\label{bcl2-eq:esti5}
\end{align}
To estimate  the second integral, we proceed as previously, using Taylor's theorem and Hypothesis \ref{hyp:J} to obtain
\begin{align}
\int_{-1}^{1} [\undu(t,x+z) - \undu(t,x)]J(z) \, dz &=\int_{-1}^{1}\int_{0}^1\int_{0}^1\partial_x^2\undu(t,x+\tau \omega z) \tau z^2J(z) \, d\tau d\omega dz \nonumber\\
&\ge \mathcal{J}_1\min_{-1<\xi<1}\partial_x^2\undu(t,x+\xi) \ge 0,
\label{bcl2-eq:esti6}
\end{align}
since $x-1 > Y(t)-1 \geq X_c(t)+R_0 -1\ge X_c(t)$, so that the function $\undu$ is convex with respect to $x$ there.

Finally, the last integral is estimated by splitting it into two parts, that is
\begin{align*}
I&:=\int_{1}^{+\infty} [\undu(t,x+z) - \undu(t,x)]J(z)\,dz\\
&= \int_{1}^{B} [\undu(t,x+z) - \undu(t,x)] J(z) \, dz + \int_{B}^{+\infty} [\undu(t,x+z) - \undu(t,x)] J(z) \, dz.
\end{align*}
with $B>1$. Using Taylor's theorem, the first of these is rewritten as
\[
\int_{1}^{B} [\undu(t,x+z) - \undu(t,x)]J(z) \, dz=\int_{1}^{B}\int_{0}^1\partial_x\undu(t,x+\tau z) z J(z) \,d\tau dz.
\]
Using the expression \eqref{bcl2-def:u_x} for $\partial_x\undu$, observe that we have, for all $\tau z\ge 0$, 
\[
\partial_x\undu(t,x+\tau z) \geq 3\partial_xw(t,x+\tau z) \geq 3\partial_xw(t,x),
\]
by convexity of $w$ is convex with respect to $x$. It then follows that 
\begin{align}
\int_{1}^{B} [\undu(t,x+z) - \undu(t,x)]J(z) \, dz&\geq 3 \left(\int_{1}^{B} z J(z) \, dz\right)\partial_xw(t,x)\nonumber\\& \geq 3\mathcal{J}_0 \left(\int_{1}^{B} z^{-2s} \, dz\right)\partial_xw(t,x)\label{bcl2-eq:esti8}
\end{align}
using again Hypothesis \ref{hyp:J}. For the second integral, we have, since $\undu $ is positive, 
\begin{equation}\label{bcl2-eq:esti7}
\int_{B}^{\infty} [\undu(t,x+z) - \undu(t,x)]J(z)\,dz \geq - \frac{\mathcal{J}_0\undu(t,x)}{2sB^{2s}}.
\end{equation}
Collecting \eqref{bcl2-eq:esti5}, \eqref{bcl2-eq:esti6}, \eqref{bcl2-eq:esti7} and \eqref{bcl2-eq:esti8}, we find that, for all $t>t^\#$ and $x>Y(t)$,
\[
\opd{\undu}(t,x)\ge \frac{\lambda - \undu(t,x)}{2s\mathcal{J}_0x^{2s}} -\frac{\mathcal{J}_0\undu(t,x)}{2sB^{2s}}+3\mathcal{J}_0 \left(\int_{1}^{B} z^{-2s} \, dz\right)\partial_xw(t,x).
\]
The proof is ended by observing that $\undu(t,x)\le \undu(t,Y(t))=\lambda\tau_0$ for all $t>t^\#$ and $x\ge Y(t)$.
\end{proof}

\subsubsection{The region $Y(t) < x< 2^{\frac{1}{2s(\beta -1)}}X(t)$}
With the previous result at hand, let us next estimate $\opd{\undu}(t,x)$ when $x>Y(t)$. 

\begin{prop}\label{bcl2-prop:esti-frac1-1}
For any $0<\lambda\le \frac{1}{2}$ and any $\gamma$, $\kappa >0$, there exists $t_2 > 0$ such that, for all $t\ge t_2$ and $Y(t)< x< 2^{\frac{1}{2s(\beta -1)}}X(t)$,
\[
\opd{\undu}(t,x)+\frac{1}{2}(1-\lambda)\undu^\beta(t,x) \ge 0.
\]
\end{prop}
\begin{proof}
First let us observe that since $Y(t)$ tends to $+\infty$ as $t$ tends to $+\infty$, we may find $t'>t^{\#}$ such that, for all $t\ge t'$,
\[
\left(\frac{2\tau_0{\mathcal{J}_0}^2}{1-\tau_0}\right)^{\frac{1}{2s}}Y(t)>1.
\]
Set $B:= \left(\frac{2\tau_0{\mathcal{J}_0}^2}{1-\tau}\right)^{\frac{1}{2s}}x$, with, again, $\tau_0:=\frac{3}{2\delta_c}\left(1-\frac{1}{2\delta_c}+\frac{1}{3(2\delta_c)^2}\right)$. From Lemma \ref{lem:estDfar}, using the expression \eqref{bcl2-def:w_x} for $\partial_xw(t,x)$, we deduce that, for all $t\ge t'$ and $x> Y(t)$,
\[
\opd{\undu}(t,x)\ge \frac{\lambda(1 - \tau_0)}{4s\mathcal{J}_0x^{2s}} -6s\mathcal{J}_0 w^{\beta}(t,x)\frac{x^{2s(\beta-1)-1}}{(\kappa t)^{\beta-1}}\left(\int_{1}^{B} z^{-2s} \, dz\right).
\]
Therefore, since $\undu(t,x)\ge w(t,x)$, we get, for all $t\ge t'$ and $x> Y(t)$,
\[
\opd{\undu}(t,x)+\frac{1}{2}(1-\lambda)\undu^\beta(t,x)\ge w^\beta(t,x)\left[\frac{1}{2}(1-\lambda) -6s\mathcal{J}_0 \frac{x^{2s(\beta-1)-1}}{(\kappa t)^{\beta-1}}\left(\int_{1}^{B} z^{-2s} \, dz\right)\right].
\]
Set $\mathcal{C}_3:=\left(\frac{2\tau_0{\mathcal{J}_0}^2}{1-\tau_0}\right)^{\frac{1}{2s}}$. We now treat separately the three cases $\frac{1}{2}<s<1$, $s=\frac{1}{2}$ and $s<\frac{1}{2}$.

\paragraph*{Case $\frac{1}{2}<s<1$:}
In this situation, the above integral is bounded from above by $\frac{1}{2s-1}$ and we have, for all $t\ge t'$,
\[
\opd{\undu}(t,x)+\frac{1}{2}(1-\lambda)\undu^\beta(t,x)\ge w^\beta(t,x)\left[\frac{1}{2}(1-\lambda) -\frac{6s\mathcal{J}_0}{2s-1} \frac{x^{2s(\beta-1)-1}}{(\kappa t)^{\beta-1}}\right].
\]
Since $x< 2^{\frac{1}{2s(\beta -1)}}X(t)$ for all $t\ge t'$,
\[
\frac{x^{2s(\beta -1)-1}}{\kappa^{\beta -1} t^{\beta-1}}=\frac{t x^{2s(\beta -1)-1}}{\kappa^{\beta -1} t^{\beta}} \leq \frac{2(X(t))^{2s(\beta -1)-1}}{\kappa^{\beta -1} t^{\beta-1}}.
\]
which, using the expression \eqref{bcl2-def:X} for $X(t)$, enforces, for all $t\ge t'$ and $Y(t)<x<2^{\frac{1}{2s(\beta -1)}}X(t)$,
\[
\opd{\undu}(t,x)+\frac{1}{2}(1-\lambda)\undu^\beta(t,x)\ge w^\beta(t,x)\left[\frac{1}{2}(1-\lambda) -\frac{12s\mathcal{J}_0}{2s-1} \frac{\left[\lambda^{1-\beta}+\gamma(\beta-1)t\right]}{X(t)}\right].
\]
Using the fact that both $\frac{1}{X(t)}$ and $\frac{t}{X(t)}$ tend to $0$ as $t$ tends to $+\infty$ according to \eqref{bcl2-eq:esti-t/X(t)}, we may find $t_2>t'$ such that, for all $t\ge t_2$,
\[
\frac{1}{2}(1-\lambda) \ge \frac{12s\mathcal{J}_0}{2s-1} \frac{\left[\lambda^{1-\beta}+\gamma(\beta-1)t\right]}{X(t)}.
\]
\paragraph*{Case $s=\frac{1}{2}$:}
In this situation, the above integral is bounded from above by $\ln(B)$ and we thus get, for all $t\ge t'$ and $Y(t)<x<2^{\frac{1}{2s(\beta -1)}}X(t)$,
\[
\opd{\undu}(t,x)+\frac{1}{2}(1-\lambda)\undu^\beta(t,x)\ge w^\beta(t,x)\left[\frac{1}{2}(1-\lambda) -12s\mathcal{J}_0 \frac{\left[\lambda^{1-\beta}+\gamma(\beta-1)t\right]\ln(2^{\frac{1}{2s(\beta-1)}}\mathcal{C}_3X(t))}{X(t)}\right].
\]
It follows from \eqref{bcl2-def:X} that $\ln(X(t))\lesssim\ln(t)$, and from \eqref{bcl2-eq:esti-tlnt/X(t)} that we can find $t_2>t'$ such that, for all $t\ge t_2$,
\[
\frac{1}{2}(1-\lambda) \ge -12s\mathcal{J}_0 \frac{\left[\lambda^{1-\beta}+\gamma(\beta-1)t\right]\ln(2^{\frac{1}{2s(\beta-1)}}\mathcal{C}_3X(t))}{X(t)}.
\]

\paragraph*{Case $0<s<\frac{1}{2}$:}
In this last situation, the integral is bounded from above by $\frac{{\mathcal{C}_3}^{1-2s}x^{1-2s}}{1-2s}$ and therefore, for all $t\ge t'$ and $Y(t)<x<2^{\frac{1}{2s(\beta -1)}}X(t)$,
\[
\opd{\undu}(t,x)+\frac{1}{2}(1-\lambda)\undu^\beta(t,x)\\\ge w^\beta(t,x)\left[\frac{1}{2}(1-\lambda) -\frac{6s\mathcal{J}_0{\mathcal{C}_{3}}^{1-2s}}{1-2s} \frac{x^{2s(\beta-1)-2s}}{(\kappa t)^{\beta-1}}\right].
\]
which, using that $X(t)\le x\le 2^{\frac{1}{2s(\beta-1)}}X(t)$, yields
\[
\opd{\undu}(t,x)+\frac{1}{2}(1-\lambda)\undu^\beta(t,x)\\\ge w^\beta(t,x)\left[\frac{1}{2}(1-\lambda) -\frac{12s\mathcal{J}_0\mathcal{C}_3^{1-2s}}{1-2s} \frac{\left[\lambda^{1-\beta}+\gamma(\beta-1)t\right]}{X^{2s}(t)}\right].
\]
Using again \eqref{bcl2-def:X} and since $\frac{t}{X^{2s}(t)}$ tends to $0$ as $t$ tends to $+\infty$ according to \eqref{bcl2-eq:esti-t/X2s(t)}, we may find $t_2>t'$ such that, for all $t\ge t_2$
\[
\frac{1}{2}(1-\lambda) \ge \frac{12s\mathcal{J}_0\mathcal{C}_3^{1-2s}}{1-2s} \frac{\left[\lambda^{1-\beta}+\gamma(\beta-1)t\right]}{X^{2s}(t)}.
\]

In conclusion, for any $s$ in $(0,1)$, we have found $t_2>t'>1$ such that, for all $t\ge t_2$ and $Y(t)< x< 2^{\frac{1}{2s(\beta -1)}}X(t)$,
\[
\opd{\undu}(t,x)+\frac{1}{2}(1-\lambda)\undu^\beta(t,x)\ge 0.
\]
\end{proof}

\subsubsection{The region $x\ge 2^{\frac{1}{2s(\beta -1)}}X(t)$}
In this region, we claim the following.

\begin{prop}\label{bcl2-cla:esti-frac3}
For all $\lambda\le \frac{1}{2}$, $\gamma$ and $\kappa$, there exists $t_3>1$ such that, for all $t\ge t_3$ and $x\ge 2^{\frac{1}{2s(\beta -1)}}X(t)$,
\[
\opd{\undu}(t,x) \ge \frac{\lambda(1-\tau_0)}{8s\mathcal{J}_0x^{2s}}, 
\]
where $\tau_0=\frac{3}{2\delta_c}\left(1-\frac{1}{2\delta_c}+\frac{1}{3(2\delta_c)^2}\right)$.
\end{prop}
\begin{proof}
The proof follows the same steps as for Proposition \ref{bcl2-prop:esti-frac1-1}, but with some adaptations. We first set the same values for $B$ and $t'>t^{\#}$, and observe that, by a straightforward computation from the expression \eqref{bcl2-def:X} for $X(t)$, there exists $t''>0$ such that $2^{\frac{1}{2s(\beta-1)}}X(t)>Y(t)$ for all $t\ge t''$.

From Lemma \ref{lem:estDfar}, we then have, for $t\ge\sup\{t',t''\}$ and $x\ge 2^{\frac{1}{2s(\beta -1)}}X(t)$,
\[
\opd{\undu}(t,x)\ge \frac{\lambda(1 - \tau_0)}{4s\mathcal{J}_0x^{2s}} -6s\mathcal{J}_0 w^{\beta}(t,x)\frac{x^{2s(\beta-1)-1}}{(\kappa t)^{\beta-1}}\left(\int_{1}^{B} z^{-2s} \, dz\right).
\]
Using Proposition \ref{bcl2-prop:w2}, we get 
$$w^{\beta}(t,x)\le 2^{\frac{\beta}{\beta -1}}\frac{(\kappa t)^{\beta}}{x^{2s\beta}}$$
and therefore 
\begin{align*}
\opd{\undu}(t,x)&\ge \frac{\lambda(1 - \tau_0)}{4s\mathcal{J}_0x^{2s}} -\frac{6s\mathcal{J}_02^{\frac{\beta}{\beta-1}}}{x^{2s}} \frac{\kappa t}{x}\left(\int_{1}^{B} z^{-2s} \, dz\right)\\
&\ge \frac{1}{x^{2s}}\left[\frac{\lambda(1 - \tau_0)}{4s\mathcal{J}_0} -6s\mathcal{J}_02^{\frac{\beta}{\beta-1}} \frac{\kappa t}{x}\left(\int_{1}^{B} z^{-2s} \, dz\right)\right].
\end{align*}
Finally, by considering separately the three cases $\frac12<s<1$, $s=\frac12$, $0<s<\frac12$ and reproducing the argument previously used in the proof of Proposition \ref{bcl2-prop:esti-frac1-1}, we may find $t_3>1$ such that, for all $t\ge t_3$ and $x\ge 2^{\frac{1}{2s(\beta-1)}}X(t)$,
\[
\opd{\undu}(t,x)\ge \frac{\lambda(1 - \tau_0)}{8s\mathcal{J}_0x^{2s}}.
\]
\end{proof}

\subsection{Tuning the parameters $\kappa$ and $\gamma$}
We must now choose the values of the parameters $\gamma$ and $\kappa$ in such a way that, for some $t^*>1$, the function  $\undu(\cdot+t^*-1,\cdot)$ satisfies \eqref{bcl2-eq:subsol} for all $t\ge 1$. We first recall that $\undu(\cdot+t^*-1,\cdot)$ satisfies \eqref{bcl2-eq:subsol} for all $t\ge 1$ if and only if inequations \eqref{eq:subsolleft} and \eqref{eq:subsolright} hold simultaneously for all $t\ge t^*$. Since \eqref{eq:subsolleft} holds unconditionally for $t$ sufficiently large, the only thing left to check is that \eqref{eq:subsolright} holds for a suitable choice of the values of $\gamma$ and $\kappa$.

By using \eqref{bcl2-def:u_t} and \eqref{bcl2-def:w_t}, inequation \eqref{eq:subsolright} holds, in particular, for all $t\ge 1$ and $x > X(t)$,
\[
3\frac{\partial_t\Phi(t,x)}{\Phi^\beta(t,x)}w^{\beta}(t,x)\leq \opd{\undu}(t,x) + (1-\lambda)\undu^{\beta}(t,x) -\gamma w^{\beta}.
\]
We next set $t^*:=\sup\{t_0,t_1,t_2,t_3\}$, where $t_0$, $t_1$, $t_2$ and $t_3$ are respectively determined in Propositions \ref{bcl2-cla:esti-frac1-1}, \ref{bcl-prop:esti-reg2}, \ref{bcl2-prop:esti-frac1-1} and \ref{bcl2-cla:esti-frac3}, and decompose the set $[X(t),+\infty)$ into the two subsets
\[
I_1:= [X(t),2^{\frac{1}{2s(\beta -1)}}X(t))\text{ and }I_2:=[2^{\frac{1}{2s(\beta -1)}}X(t),+\infty).
\]

In the first interval, we have the following result.
\begin{lem}
For any $\lambda<\frac12$, there exists $\gamma^*$ such that, for all $\kappa$ and $\gamma\le \gamma^*$, one has, for $t\ge \sup\{\frac{48}{\lambda^{\beta-1}(1-\lambda)},t^*\}$ and $x$ in $I_1$, 
\[
3\frac{\partial_t\Phi(t,x)}{\Phi^\beta(t,x)}w^{\beta}(t,x)\leq \opd{\undu}(t,x) + (1-\lambda)\undu^{\beta}(t,x) -\gamma w^{\beta}(t,x).
\]
\end{lem}
\begin{proof}
By the definition of $\Phi$, we have, for all $t\ge 1$ and $X(t)\le x$, 
\[
3\frac{\partial_t\Phi(t,x)}{\Phi^\beta(t,x)}w^{\beta}(t,x) = \frac{3}{t}\frac{x^{2s(\beta -1)}}{(\kappa t)^{\beta-1}}w^{\beta}(t,x).
\]
By exploiting the definition of $X(t)$, it follows that, for all $t\ge 1$ and $x< 2^{\frac{1}{2s(\beta -1)}}X(t)$,
\[
3\frac{\partial_t\Phi(t,x)}{\Phi^\beta(t,x)}w^{\beta}(t,x)= \frac{3}{t}\frac{x^{2s(\beta -1)}}{(\kappa t)^{\beta-1}}w^{\beta}(t,x)\le \left[\frac{6}{t}\left(\frac{1}{\lambda}\right)^{\beta-1}+6\gamma(\beta-1)\right] w^{\beta}(t,x).
\]
Let $\gamma_0:=\frac{1-\lambda}{48(\beta -1)}$, then, for all $\gamma \le \gamma_0$, we have, for all $t\ge 1$ and $x$ in $I_1$, 
\[
3\frac{\partial_t\Phi(t,x)}{\Phi^\beta(t,x)}w^{\beta}(t,x)\le \left[\frac{6}{t}\left(\frac{1}{\lambda}\right)^{\beta-1}+\frac{1-\lambda}{8}\right] w^{\beta}(t,x),
\]
which, for all $t$ large enough, for instance $t\ge\frac{48\lambda^{1-\beta}}{1-\lambda}$, and $x$ in $I_1$, gives
\[
3\frac{\partial_t\Phi(t,x)}{\Phi^\beta(t,x)}w^{\beta}(t,x)\le \frac{1-\lambda}{4} w^{\beta}(t,x).
\]
From Propositions \ref{bcl-prop:esti-reg2} and \ref{bcl2-prop:esti-frac1-1}, we have, for all $t\ge t^*$ and $x$ in $I_1$, 
\[
\opd{\undu}(t,x) + (1-\lambda)\undu^{\beta}(t,x)-\gamma w^{\beta}(t,x)\ge \left(\frac{1-\lambda}{2}-\gamma\right)w^{\beta}(t,x),
\]
since $\undu(t,x)\ge w(t,x)$ for all $x\ge X(t)$. We end the proof by taking $\gamma^*:=\inf\{\gamma_0,\frac{1-\lambda}{4}\}$ and $t\ge \sup\{\frac{48 \lambda^{1-\beta}}{1-\lambda},t^*\}$.
\end{proof}

We next turn to the interval $I_2$.

\begin{lem}
For any $\lambda\le \frac12$, there exists $\kappa^*$ such that, for all $\gamma\le \gamma^*$ and $\kappa\le \kappa^*$, one has for all $t\ge t^*$ and $x$ in $I_2$, 
\[
3\frac{\partial_t\Phi(t,x) }{\Phi^\beta(t,x) }w^{\beta}(t,x) \leq \opd{\undu}(t,x) + (1-\lambda)\undu^{\beta}(t,x) -\gamma w^{\beta}(t,x).
\]
\end{lem}

\begin{proof}
As in the above proof, by the definition of $\Phi$, we have, for all $t\ge 1$ and $x>X(t)$,
\[
3\frac{\partial_t\Phi(t,x)}{\Phi^\beta(t,x)}w^{\beta}(t,x)= 3\kappa\frac{x^{2s(\beta -1)}}{(\kappa t)^{\beta}}w^{\beta}(t,x).
\]
By Proposition \ref{bcl2-prop:w2}, we have, for all $t\ge 1$ and $x$ in $I_2$,
\[
w^{\beta}(t,x)\le 2^{\frac{\beta}{\beta-1}}\frac{(\kappa t)^{\beta}}{x^{2s\beta}},
\]
therefore, for all $t\ge 1$ and $x$ in $I_2$,
 \[
 3\frac{\partial_t\Phi(t,x)}{\Phi^\beta(t,x)}w^{\beta}(t,x)\le 3\kappa 2^{\frac{\beta}{\beta-1}}\frac{1}{x^{2s}}.
 \]
Finally, by Proposition \ref{bcl2-cla:esti-frac3}, we have, for all $t\ge t^*$ and $x$ in $I_2$, 
\[
\opd{\undu}(t,x) + (1-\lambda)\undu^{\beta}(t,x)-\gamma w^{\beta}(t,x)\ge \left(\frac{1-\lambda}{2}-\gamma\right)w^{\beta}(t,x)+\frac{\lambda(1-\tau_0)}{8\mathcal{J}_0s x^{2s}},
\]
since $\undu(t,x)\ge w(t,x)$ for all $x\ge X(t)$. The lemma is then proved by taking $\gamma\le \gamma_0$ and $\kappa \le \kappa^*:=\frac{\lambda(1-\tau_0)}{24(2^{\frac{\beta}{\beta-1}}\mathcal{J}_0 s)}$.
\end{proof}

\subsection{Conclusion}
From the preceding subsections, for any $\lambda\le \frac{1}{2}$, there exist $\kappa^*$, $\gamma^*$ and $t^*$ such that the function $\undu(\cdot+t^*-1,\cdot)$ satisfies inequation \eqref{bcl2-eq:subsol} for all $t\ge 1$ and $x$ in $\mathbb{R}$. To conclude the proof, we need, as in \Cref{sec:estimate1}, to check that, for some time $T$, we have $u(T,x)\ge \undu(t^*,x)$ for all $x$ in $\mathbb{R}$. 

To do so, let us observe that $\undu(t^*,\cdot) \le 3w(t^*,\cdot)$ and, by using the definition of $w$, we have 
\[
\lim_{x\to +\infty}x^{2s}\undu(t^*,x)\le 3\kappa^*t^*.
\]
It also follows from Proposition \ref{bcl-prop-flatnonlin} that there exists $t_{3\kappa^*t^*}>0$ such that, for all $t\ge t_{3\kappa^*t^*}$,
\[
\lim_{x\to +\infty}x^{2s}u(t,x)\ge 3\kappa^*t^*.
\]
Therefore, we have $\undu(t^*,x)\le u(t_{3\kappa^*t^*},x)$ for $x$ large enough, for instance $x>x_0$. Moreover, due to the monotonic behaviour of $u$, $\undu(t^*,x)\le u(t,x)$ for all $t\ge t_{3\kappa^*t^*}$ and $x>x_0$. On the other hand, by Proposition \ref{bcl-prop:inva}, $u(t,\cdot)$ tends to $1$ uniformly in $(-\infty,x_0]$ as $t$ tends to infinity, and, since $\undu \le \frac{1}{2}<1$, we can find $\hat t>0$ such that $\undu(t^*,x)\le \frac{1}{2}< u(\hat t,x)$ for all $x\le x_0$. Thus, by taking $T\ge\sup\{t_{3\kappa^*t^*},\hat t\}$, we obtain that $\undu(t^*,x)\le u(T,x) $ for all $x$ in $\mathbb{R}$. 
 
The proof of \Cref{thm:main} is then complete for all $\lambda\le \frac{1}{2}$. To obtain the speed of the level line for $\lambda \ge \frac{1}{2}$, we can reproduce the proof used in \cite{Alfaro2017,Coville2021} using the adequate invasion property, namely Proposition \ref{bcl-prop:inva}.

\section{Numerical exploration}\label{sec:numerics}
In this Section, we provide, in the particular case of the fractional Laplace operator, numerical experiments illustrating the theoretical findings reported in the present work.

To compute approximations to the solution of the Cauchy problem \eqref{eq:main}, the integro-differential equation is first discretised in space using a quadrature rule-based finite difference method on a uniform Cartesian grid, and then integrated in time using an implicit-explicit (IMEX) scheme. To do so, one needs to set the problem on a bounded domain, which is achieved by truncating the real line to a bounded interval and imposing an \emph{exterior} boundary condition.

The integral representation of the fractional Laplacian involves a singular integrand, and proper care is needed when discretising this operator. A common approach to deal with this difficulty is to split the singular integral into a sum of an isolated contribution from the singular part with another having a smooth integrand and on which standard quadrature rules can be employed. Such a strategy has been used to solve both non-local (see \cite{Tian2013}) and fractional (see \cite{Huang2014,Duo2018,Minden2020}) diffusion models. In the present work, we followed the splitting approach introduced in \cite{Duo2018}. It consists in writing the singular integral representation of the fractional Laplacian as a weighted integral of a weaker singular function by introducing a splitting parameter, namely, for any $s$ in $(0,1)$ and for all $x$ in $\mathbb{R}$,
\[
(-\Delta)^su(x)=C_{1,s}\,\text{P.V.}\int_{\mathbb{R}}\frac{u(x)-u(y)}{\vert{x-y}\vert^{\gamma}}\vert{x-y}\vert^{\gamma-1-2s}\,\mathrm{d}y\text,
\]
where $\gamma$ is a real number appropriately chosen in $(2s,2)$. The discretisation of the fractional Laplacian in a bounded interval $\Omega=(a,b)$, such that $b-a=L>0$, with the extended Dirichlet boundary condition $u=g$ in $\mathbb{R}\setminus\Omega$, works as follows. Using a uniform Cartesian grid $\{x_j=a+j(\Delta x)\,|\,j\in\mathbb{Z}\}$, with $\Delta x=\frac{L}{M}$ for some nonzero naturel integer $M$, the fractional operator, evaluated at a given gridpoint $x_j$ in $\Omega$ (that is, for $j$ in $\{0,\dots,M\}$), is then decomposed into two parts
\begin{multline}\label{integral splitting}
(-\Delta)^su(x_j)=-C_{1,s}\int_0^L\frac{u(x_j-z)-2\,u(x_j)+u(x_j+z)}{z^{1+2s}}\,\mathrm{d}z\\-C_{1,s}\int_L^{+\infty}\frac{u(x_j-z)-2\,u(x_j)+u(x_j+z)}{z^{1+2s}}\,\mathrm{d}z.
\end{multline}
The first integral in the decomposition being singular, the splitting is used. Denoting $z_k=k(\Delta x)$, for any integer $k$ in $\{0,\dots,M\}$, one first writes
\[
\int_0^L\tfrac{u(x_j-z)-2\,u(x_j)+u(x_j+z)}{z^{1+2s}}\,\mathrm{d}z=\sum_{k=1}^M\int_{z_{k-1}}^{z_k}\tfrac{u(x_j-z)-2\,u(x_j)+u(x_j+z)}{z^{\gamma}} z^{\gamma-1-2s}\,\mathrm{d}z.
\]
For any index $k$ in $\{2,\dots,M\}$, the integral
\[
\int_{z_{k-1}}^{z_k}\tfrac{u(x_j-z)-2\,u(x_j)+u(x_j+z)}{z^{\gamma}} z^{\gamma-1-2s}\,\mathrm{d}z
\]
in the above sum is regular and approximated using the weighted trapezoidal rule, that is, it is replaced by the approximate value
\[
\tfrac{1}{2(\gamma-2s)}\left(\tfrac{u(x_j-z_{k-1})-2\,u(x_j)+u(x_j+z_{k-1})}{{z_{k-1}}^{\gamma}}+\tfrac{u(x_j-z_k)-2\,u(x_j)+u(x_j+z_k)}{{z_k}^{\gamma}}\right)\left({z_k}^{\gamma-2s}-{z_{k-1}}^{\gamma-2s}\right).
\]
For $k=1$, assuming that the solution $u$ is smooth enough (of class $\mathscr{C}^2$ for instance), the corresponding integral can also be formally approximated by the weighted trapezoidal rule, that is
\[
\int_{z_0}^{z_1}\frac{u(x_j-z)-2\,u(x_j)+u(x_j+z)}{z^{\gamma}} z^{\gamma-1-2s}\,\mathrm{d}z\approx\frac{(\Delta x)^{\gamma-2s}}{2(\gamma-2s)}\frac{u(x_{j-1})-2\,u(x_j)+u(x_{j+1})}{(\Delta x)^{\gamma}}.
\]
Note that an optimal convergence rate for this scheme is obtained for $\gamma=1+s$ (see the discussion in \cite{Duo2018}).

Next, observe that, for any $z$ larger than $L$, $x_j\pm z$ belongs to $\mathbb{R}\setminus\Omega$ and thus the value of $u(x_j\pm z)$ is given by the extended Dirichlet boundary condition. As a consequence, the second integral in \eqref{integral splitting} reduces to
\[
\int_L^{+\infty}\frac{u(x_j-z)-2\,u(x_j)+u(x_j+z)}{z^{1+2s}}\,\mathrm{d}z=-\frac{1}{sL^{2s}}\,u(x_j)+\int_L^{+\infty}\frac{g(x_j-z)+g(x_j-z)}{z^{1+2s}}\,\mathrm{d}z,
\]
and may be computed explicitly depending on the extended boundary datum $g$. For the problem at hand, it is known that the solution tends to $1$ at $-\infty$ and $0$ at $+\infty$ and we used boundary datum with constant value $1$ or $0$ where appropriate.

A forward-backward Euler $(1,1,1)$ IMEX scheme (see \cite{Asher1997}), with stepsize $\Delta t$, is then applied to the semi-discretized equation, the diffusion term in the equation being treated implicitly (by the backward Euler method) and the nonlinear reaction term being dealt with explicitly (by the forward Euler method). Due to the use of a uniform grid, the resulting linear system to be solved at each step possesses a Toeplitz-type square matrix of order $M-1$. Its solution can be advantageously tackled by the Levinson recursion, at a cost of $O(M^2)$ arithmetic operations.

To cope with the algebraic decay of solutions and their spreading over a given period of time, which is necessary to observe the setting of a travelling or accelerated front, we implemented a very crude adaptation mechanism of the domain size along the iteration. At each time step, a criterion decides whether the discretisation grid is to be expanded on each side or not, according to the measured spreading of the numerical approximation at the current time and a given tolerance. This allows for discretisation points to be added to the grid (the space step being fixed one and for all at the beginning) over the course of the computation, which results in an ever increasing cost for each new iteration. The maximum number of added points at each step is a fixed parameter in the code, and, to complete the values of the approximation at these points, the boundary conditions are used, that is, the value $1$ on the left side of the grid and the value $0$ on the right one. This results in using extremely large computational domains as the simulation progresses, and thus an ever increasing computational effort\footnote{In practice, the stepsize in space $\Delta x$ is fixed and the integer $M$ is modified at each time step.}. Such a basic approach nevertheless allowed to qualitatively confirm a number of theoretical results established in the present paper, but showed its limitations in experiments in which smaller values of the fractional exponent were used, the computational domain being too small (with the parameters chosen for the computations) to correctly account for the spreading of the solution. As a consequence, the influence of the homogeneous Dirichlet boundary conditions was felt, affecting the asymptotic behaviour of the approximation.

The numerical scheme was implemented with Python using standard NumPy and SciPy libraries, notably the \texttt{scipy.linalg.solve\_toeplitz} routine to solve the Toeplitz linear system. In all the computations presented, a stepsize $\Delta t$ equal to $0.01$ was used and the starting computational domain was the interval to $(-1000,1000)$, discretised with 10001 points, that is, a stepsize $\Delta x$ equal to $0.2$ in space. The maximum number of points that could be added to each side of the domain at each iteration was $150$.

The first important feature that we were able to recover numerically is the expected dynamics of the invasion with respect to the Allee effect. Namely, after a transition period, stabilisation to a regime in which the level set of the solution evolves with a speed of order $t^{\frac{1}{2s(\beta -1)}+\frac{1}{2s}}$ occurs, as seen in \Cref{semilogy_s_comparison} below. By plotting the evolution of the position of a given level set $x_\lambda(t)$ using a semi-logarithmic scale for different values of the parameter $\beta$ and a fixed value of the parameter $s$, we observe that, except for $\beta=1$ for which the dynamics differs, the shapes of resulting curves are somehow identical, meaning that $\log(x_\lambda(t))\sim C(s,\beta) \log(t)$.

\begin{figure}[h]
\centering
\includegraphics[width=0.4\textwidth]{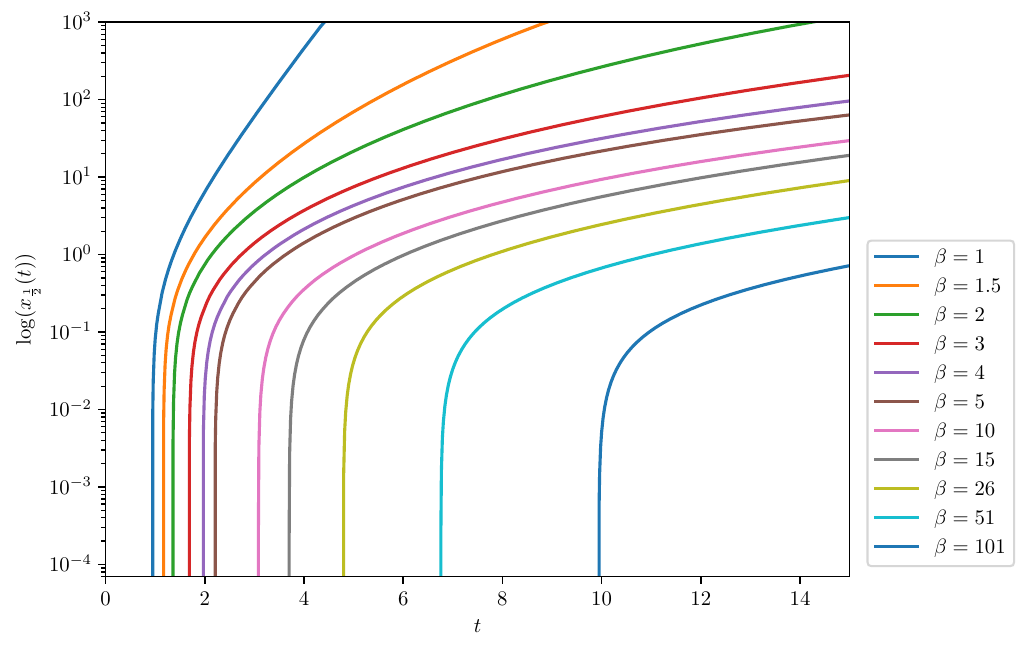}
\caption{Logarithm of the position of the level line of height $\frac{1}{2}$ of numerical approximations of the solution to the problem with fractional diffusion, plotted as a function of time, for different values of $\beta$ and $s$ equal to $\frac{1}{2}$.}\label{semilogy_s_comparison}
\end{figure}

\begin{figure}[h]
	\centering
	\begin{subfigure}[b]{0.3\textwidth}
    	\centering
    	\includegraphics[width=\textwidth]{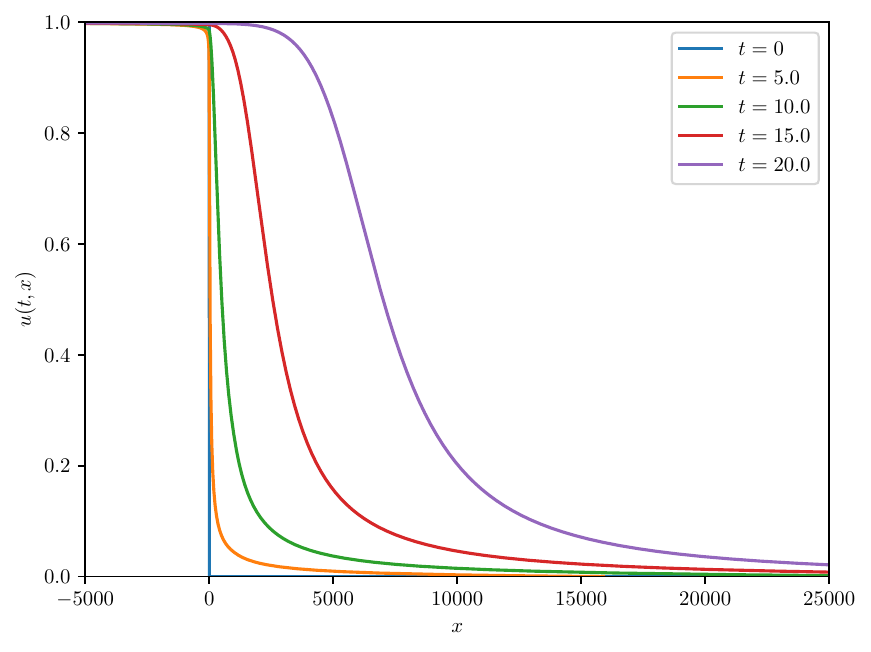}
    	\caption{$s=0.3$.}
    \end{subfigure} %
	\begin{subfigure}[b]{0.3\textwidth}
		\centering
		\includegraphics[width=\textwidth]{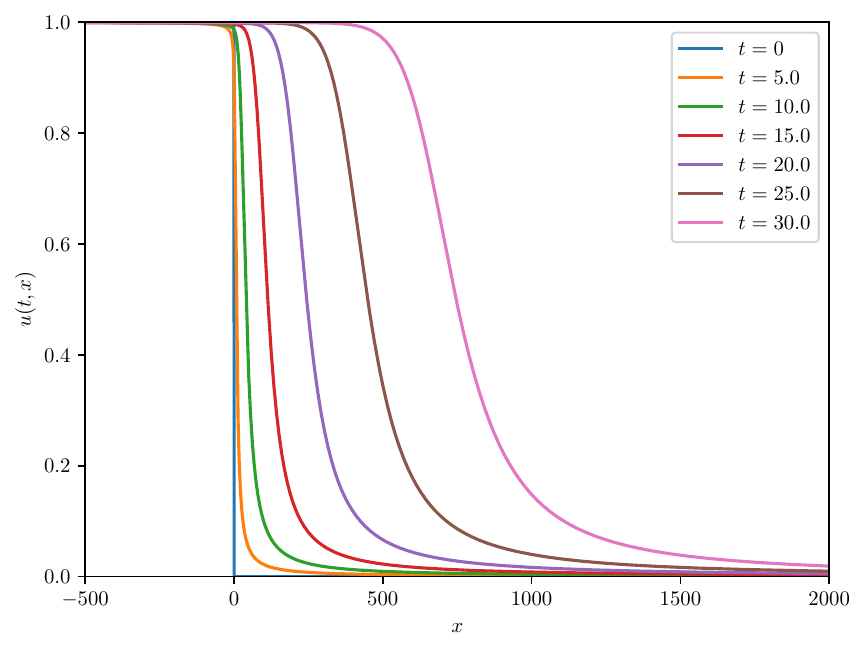}
		\caption{$s=0.5$.}
	\end{subfigure} %
	\begin{subfigure}[b]{0.3\textwidth}
		\centering
		\includegraphics[width=\textwidth]{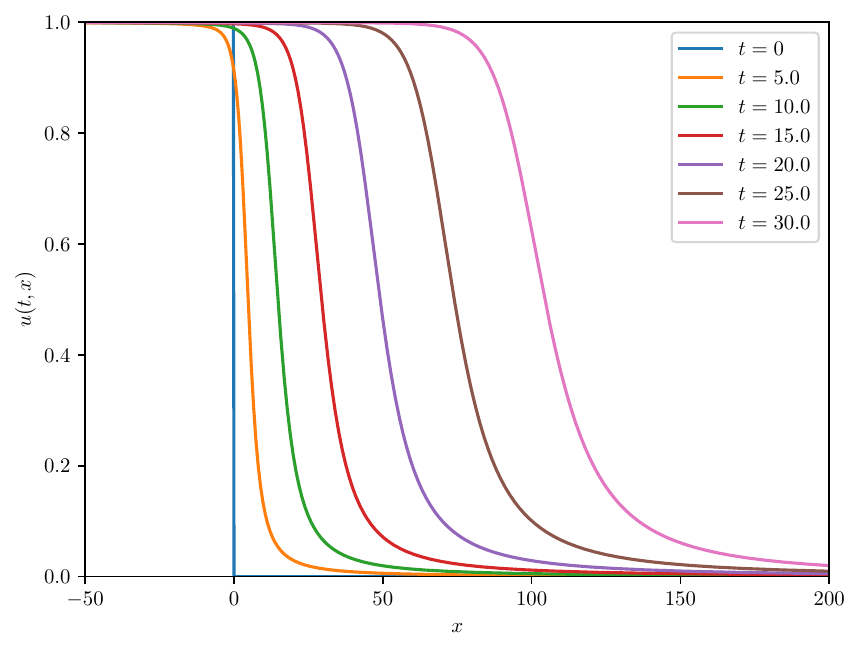}
		\caption{$s=0.7$.}
	\end{subfigure} %

	\begin{subfigure}[b]{0.3\textwidth}
    	\centering
        \includegraphics[width=\textwidth]{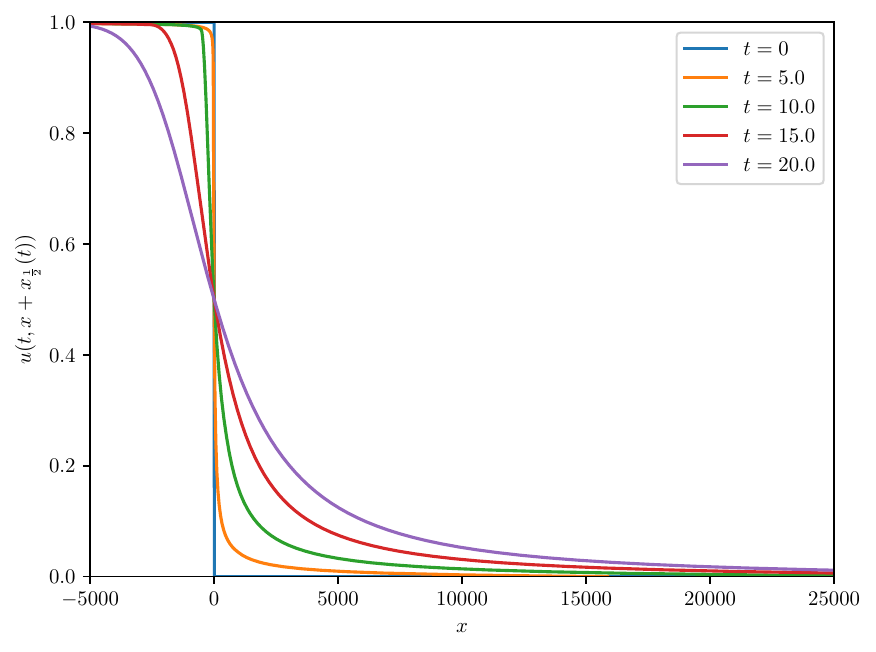}
    	\caption{$s=0.3$.}
    \end{subfigure} %
	\begin{subfigure}[b]{0.3\textwidth}
		\centering
        \includegraphics[width=\textwidth]{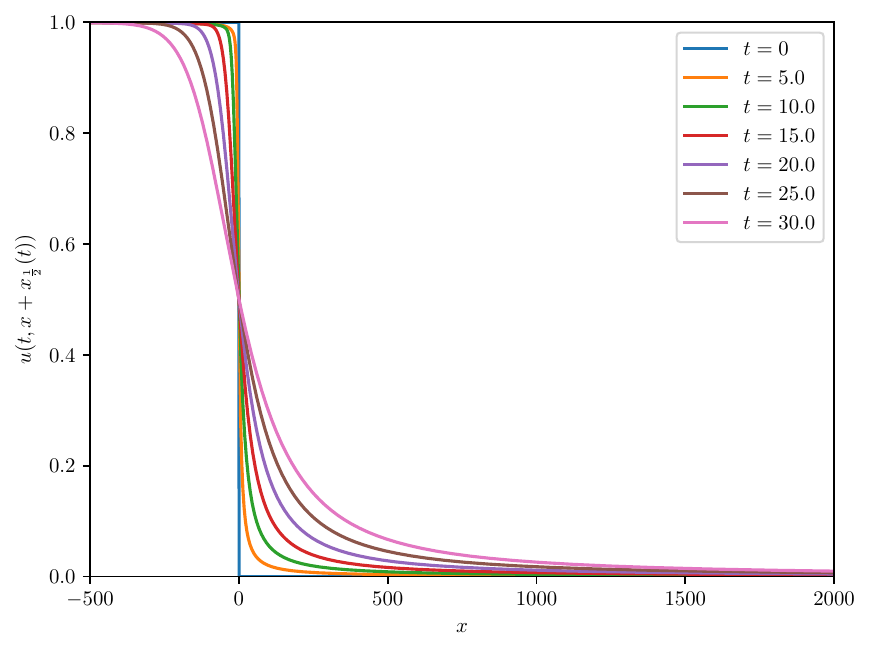}
		\caption{$s=0.5$.}
	\end{subfigure} %
	\begin{subfigure}[b]{0.3\textwidth}
		\centering
        \includegraphics[width=\textwidth]{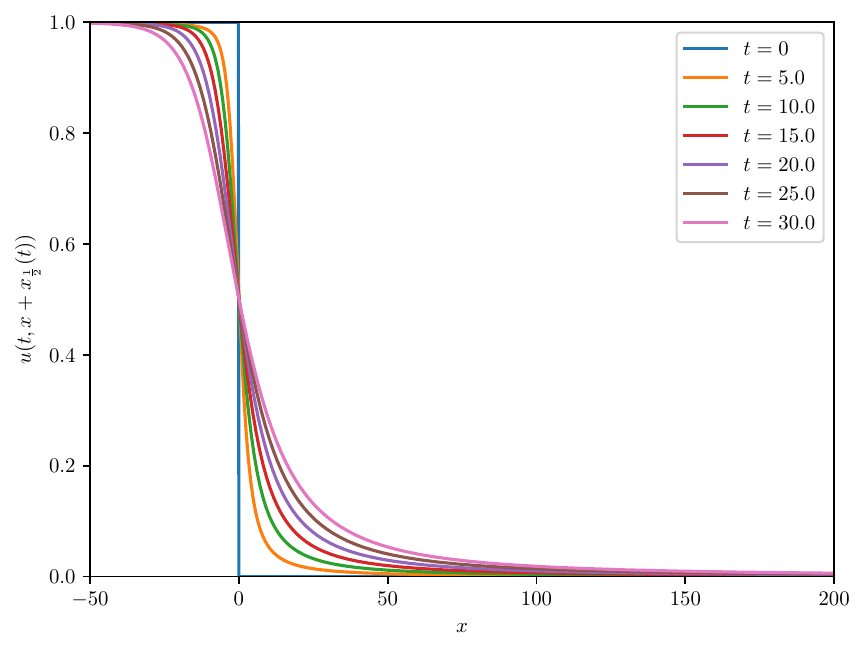}
		\caption{$s=0.7$.}
	\end{subfigure} %
\caption{Numerical approximations of the solution to the problem with fractional diffusion at different times for $\beta=1.5$ and different values of $s$. In the lower row, the graphs have been shifted by setting the position of the level line of value $\frac{1}{2}$ at $x=0$, for comparison purposes.}\label{solutions beta=1.5}
\end{figure}

\begin{figure}[h]
	\centering
	\begin{subfigure}[b]{0.3\textwidth}
    	\centering
    	\includegraphics[width=\textwidth]{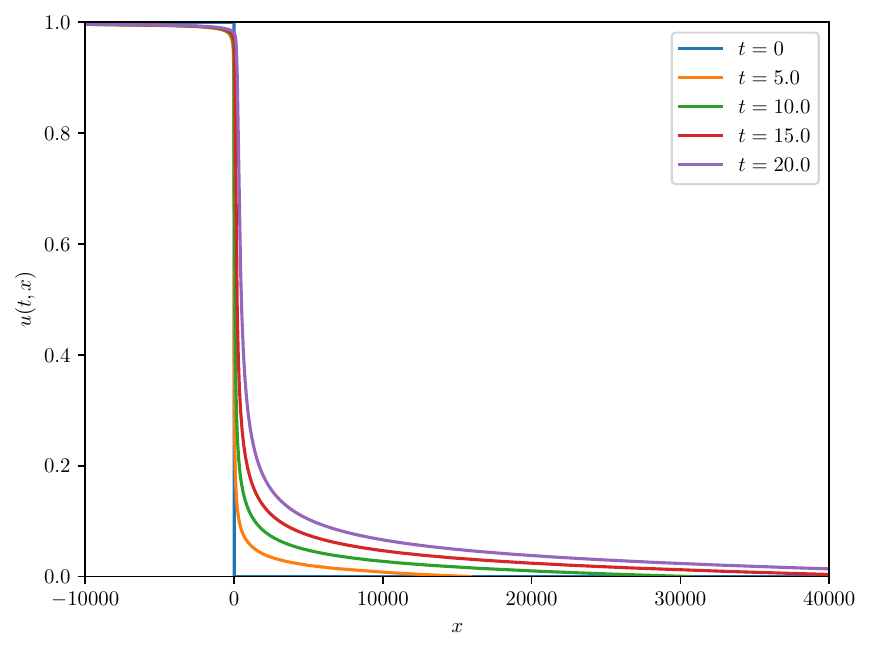}
    	\caption{$s=0.25$.}
    \end{subfigure} %
	\begin{subfigure}[b]{0.3\textwidth}
    	\centering
    	\includegraphics[width=\textwidth]{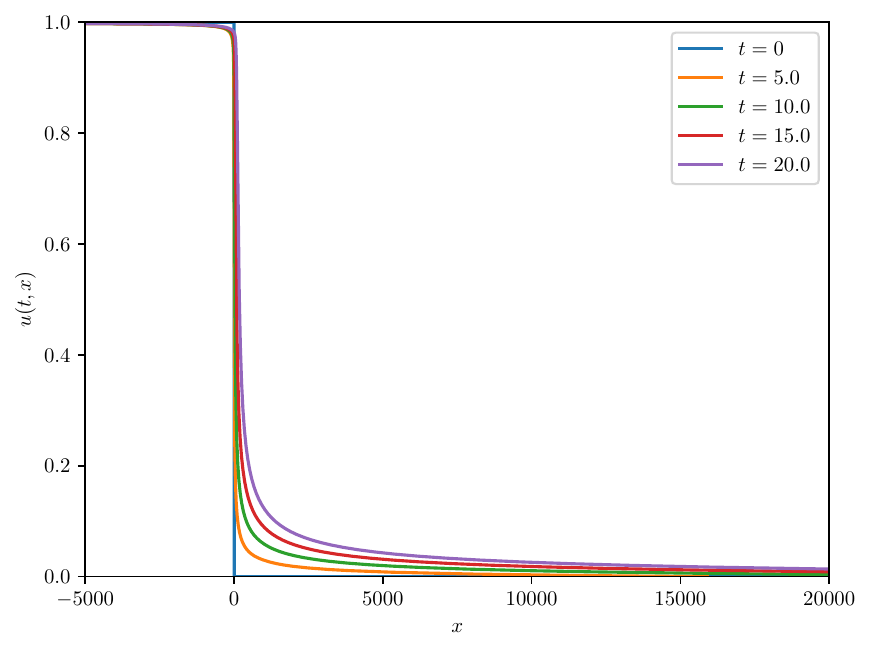}
    	\caption{$s=0.3$.}
    \end{subfigure} %
	\begin{subfigure}[b]{0.3\textwidth}
		\centering
		\includegraphics[width=\textwidth]{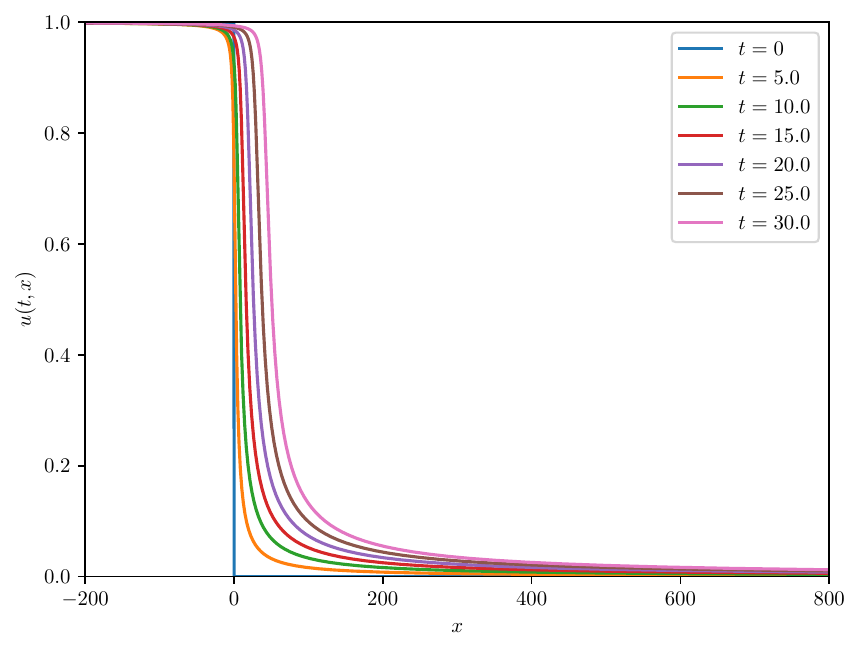}
		\caption{$s=0.5$.}
	\end{subfigure} %

	\begin{subfigure}[b]{0.3\textwidth}
    	\centering
        \includegraphics[width=\textwidth]{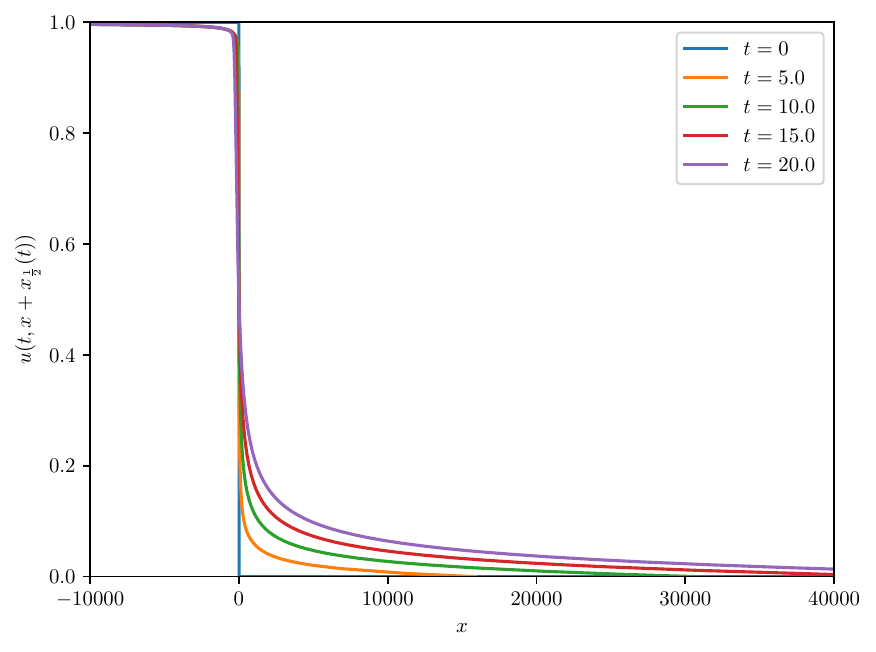}
    	\caption{$s=0.25$.}
    \end{subfigure} %
	\begin{subfigure}[b]{0.3\textwidth}
    	\centering
        \includegraphics[width=\textwidth]{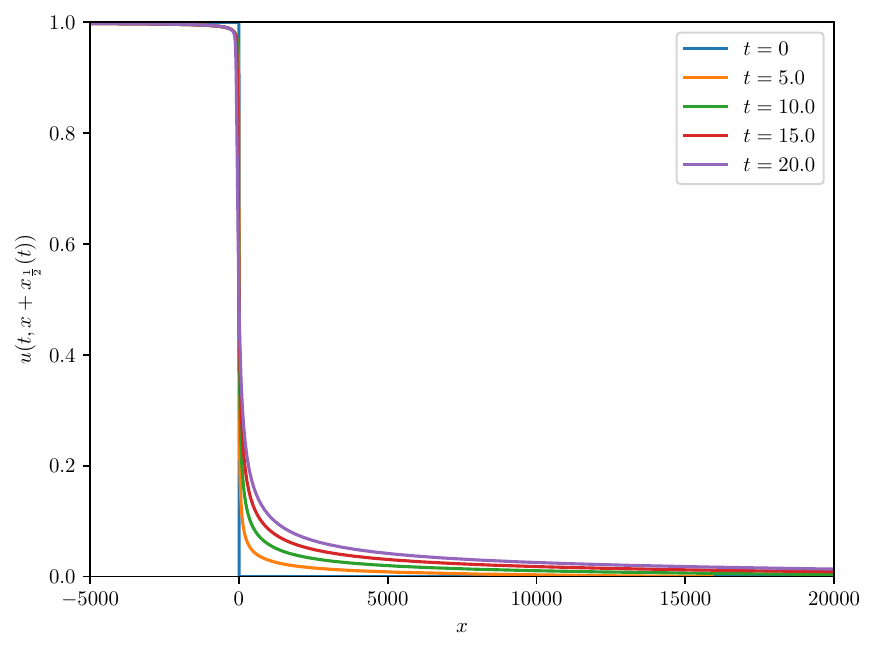}
    	\caption{$s=0.3$.}
    \end{subfigure} %
	\begin{subfigure}[b]{0.3\textwidth}
		\centering
        \includegraphics[width=\textwidth]{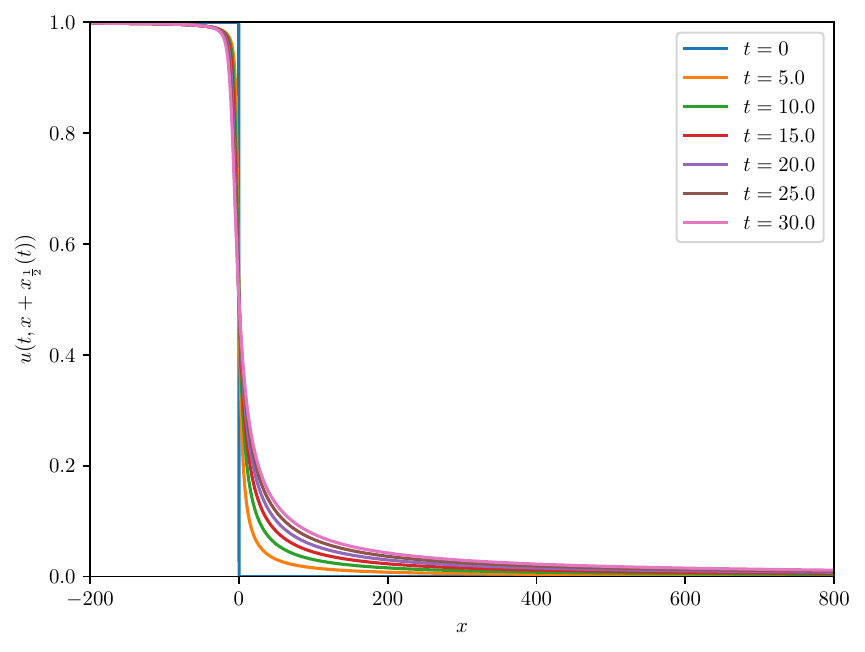}
		\caption{$s=0.5$.}
	\end{subfigure}
\caption{Numerical approximations of the solution to the problem with fractional diffusion at different times for $\beta=3$ and different values of the fractional Laplacian exponent $s$. In the lower row, the graphs have been shifted by setting the position of the level line of value $\frac{1}{2}$ at $x=0$, for comparison purposes.}\label{solutions beta=3. (a)}
\end{figure}

\begin{figure}[h]
	\centering
	\begin{subfigure}[b]{0.3\textwidth}
		\centering
		\includegraphics[width=\textwidth]{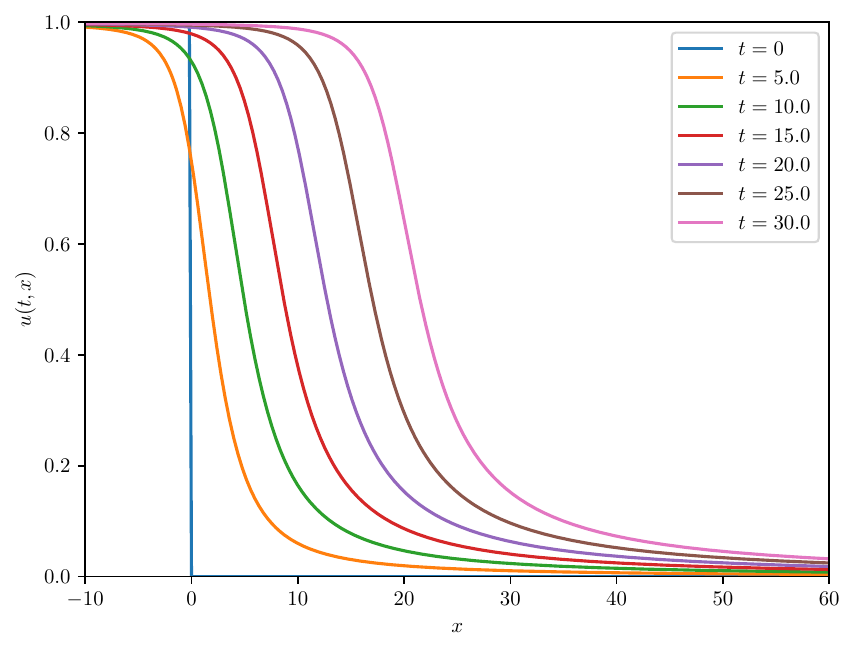}
		\caption{$s=0.7$.}
	\end{subfigure} %
	\begin{subfigure}[b]{0.3\textwidth}
    	\centering
    	\includegraphics[width=\textwidth]{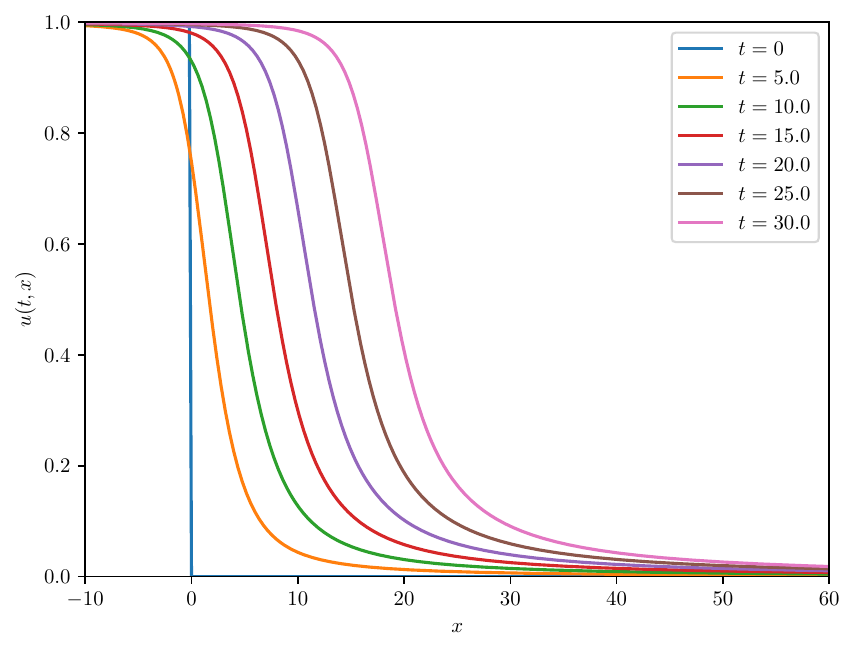}
    	\caption{$s=0.75$.}
    \end{subfigure} %
	\begin{subfigure}[b]{0.3\textwidth}
		\centering
		\includegraphics[width=\textwidth]{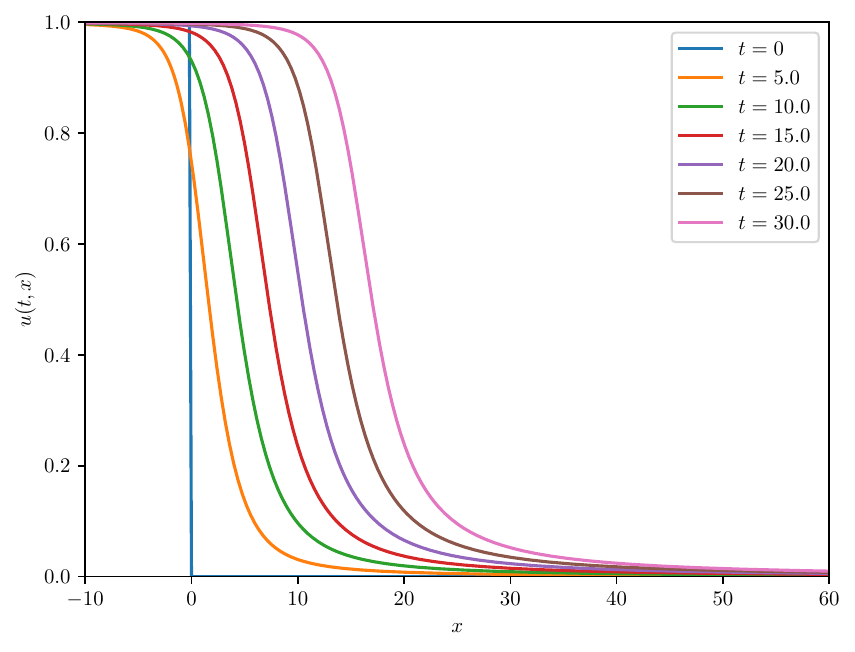}
		\caption{$s=0.8$.}
	\end{subfigure} %

	\begin{subfigure}[b]{0.3\textwidth}
		\centering
        \includegraphics[width=\textwidth]{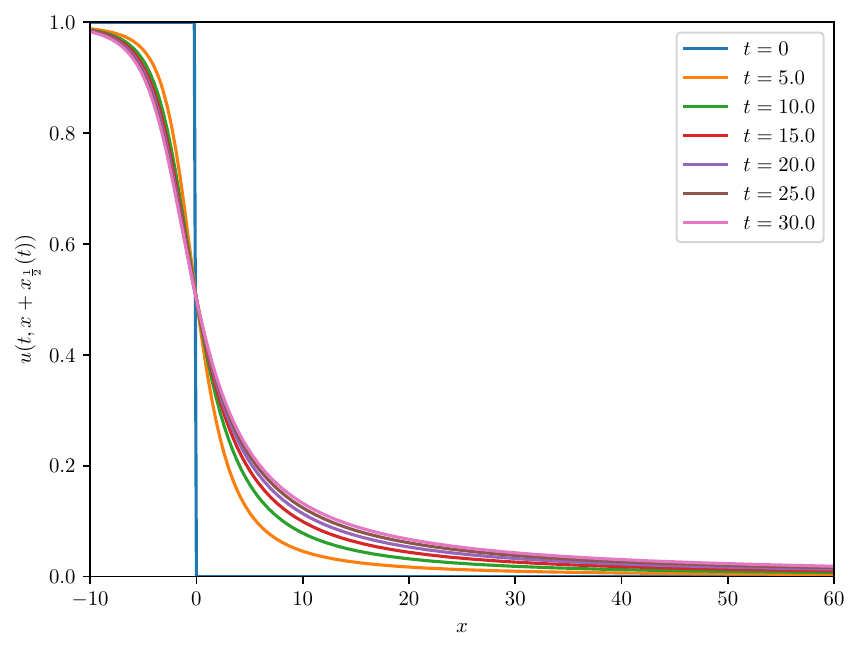}
		\caption{$s=0.7$.}
	\end{subfigure} %
	\begin{subfigure}[b]{0.3\textwidth}
    	\centering
        \includegraphics[width=\textwidth]{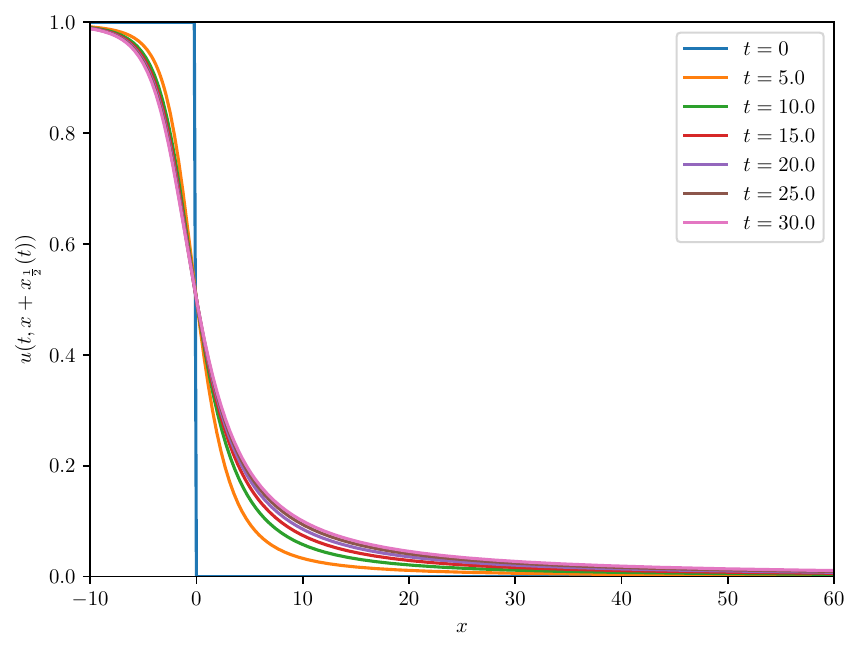}
    	\caption{$s=0.75$.}
    \end{subfigure} %
	\begin{subfigure}[b]{0.3\textwidth}
		\centering
        \includegraphics[width=\textwidth]{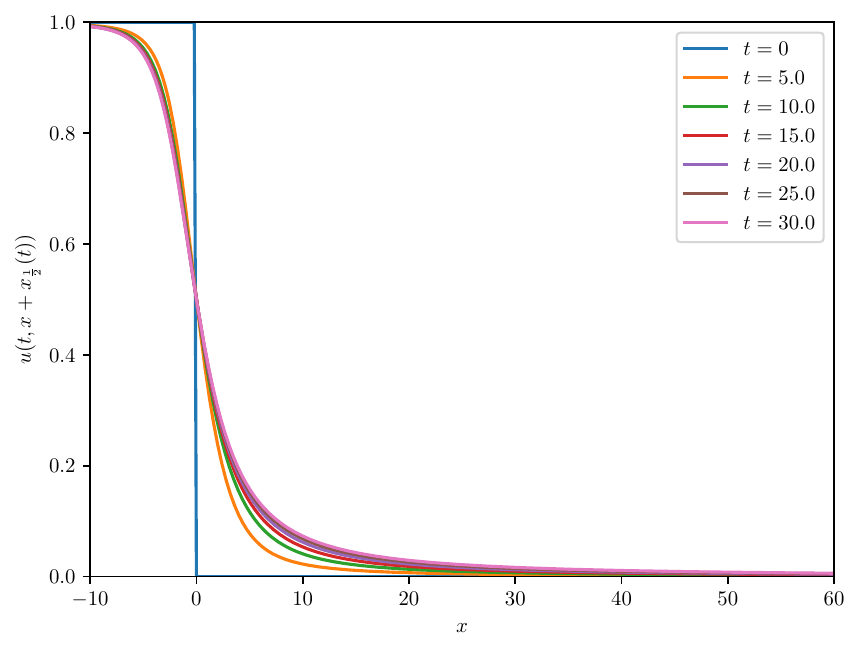}
		\caption{$s=0.8$.}
	\end{subfigure} %
\caption{Numerical approximations of the solution to the problem with fractional diffusion at different times for $\beta=3$ and different values of $s$. In the lower row, the graphs have been shifted by setting the position of the level line of value $\frac{1}{2}$ at $x=0$, for comparison purposes.}\label{solutions beta=3. (b)}
\end{figure}

In contrast, Figures \ref{solutions beta=1.5}, \ref{solutions beta=3. (a)} and \ref{solutions beta=3. (b)} illustrate the different behaviours observed when the value of the parameter $s$ varies while the value of the parameter $\beta$ is fixed. For $\beta=1.5$, \Cref{solutions beta=1.5} shows that acceleration occurs for any of the values of $s$ we considered, that is $s=0.3$, $0.5$, and $0.7$. For $\beta=3$, Figures \ref{solutions beta=3. (a)} and \ref{solutions beta=3. (b)} offer a more complex picture. According to theoretical prediction, one can observe a transition from an accelerated invasion for values of $s$ lower than $0.7$ to an invasion at constant speed for $s$ equal to $0.8$, the transition being captured for the value $s=0.75$. In both cases, it is observed that acceleration always occurs when $s<\frac{1}{2}$.

The acceleration being more pronounced for small values of the parameter $s$, one may notice that the ranges used to plot the profiles of the solution vary drastically from case to case, which may lead to some possible misinterpretations of the numerical results. There is, for instance, a ratio of about $20$ between the range used for the case $s=0.3$ and the one for the case $s=0.5$. To properly compare the profile deformation, we have plotted in Figure \ref{solution t=20, beta=3.} the shifted profile of the level set at a given time and for several values of~$s$. By doing so, we are able to observe more easily the transition that occurs at $s=0.75$, the profiles associated with values of $s$ greater than $0.75$ being very similar whereas they exhibit a noticeable deformation for lower values.

\begin{figure}[h]
\centering
\includegraphics[width=0.3\textwidth]{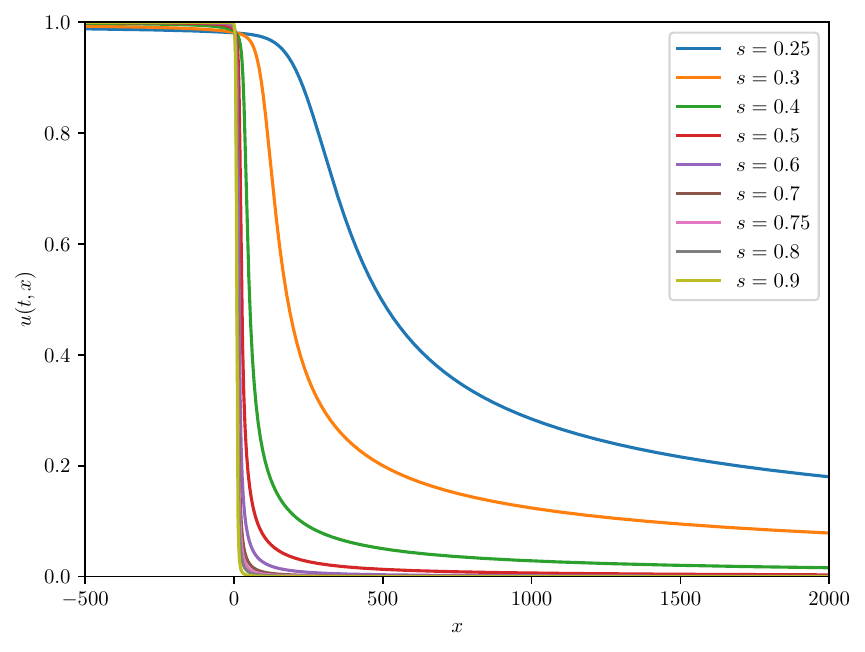}
\quad
\includegraphics[width=0.3\textwidth]{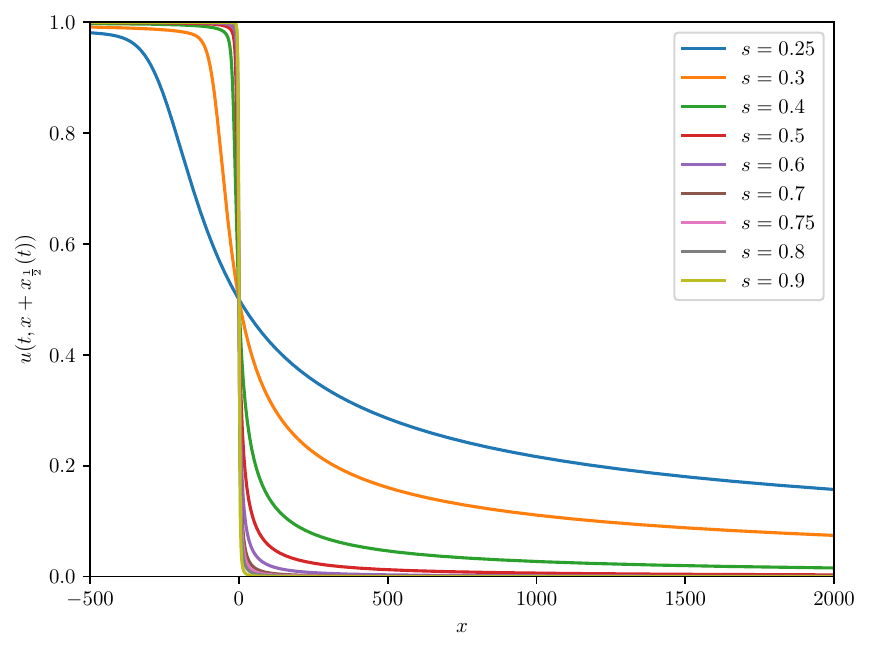}
\caption{Numerical approximations of the solution to the problem with fractional diffusion at time $t=20$ for $\beta=3$ and different values of $s$. On the right, the graphs have been shifted by setting the position of the level line of value $\frac{1}{2}$ at $x=0$, for comparison purposes.}\label{solution t=20, beta=3.}
\end{figure}

\begin{figure}[h]
	\centering
	\begin{subfigure}[b]{0.24\textwidth}
		\centering
		\includegraphics[width=\textwidth]{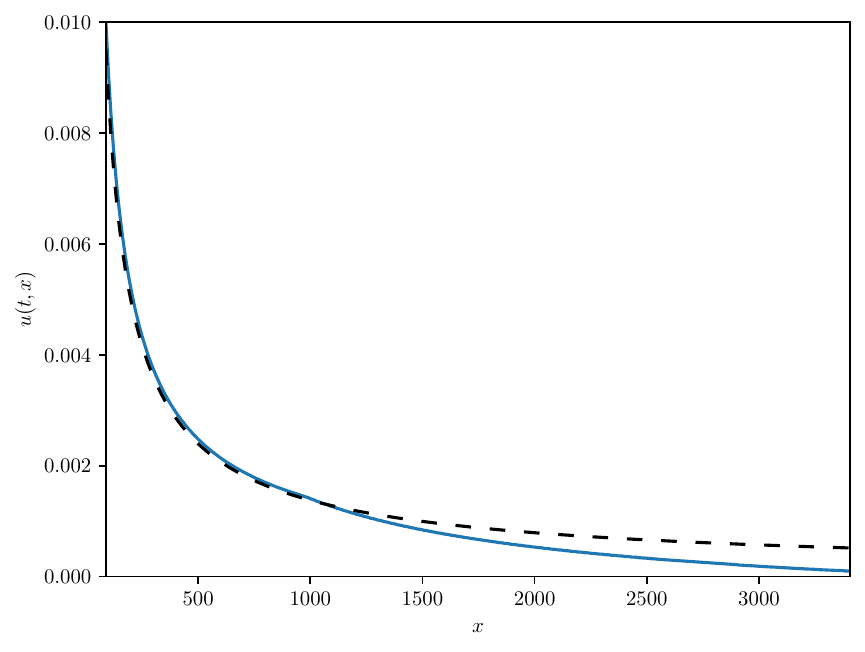}
		\caption{$s=0.4$.}
	\end{subfigure} %
	\begin{subfigure}[b]{0.24\textwidth}
        \includegraphics[width=\textwidth]{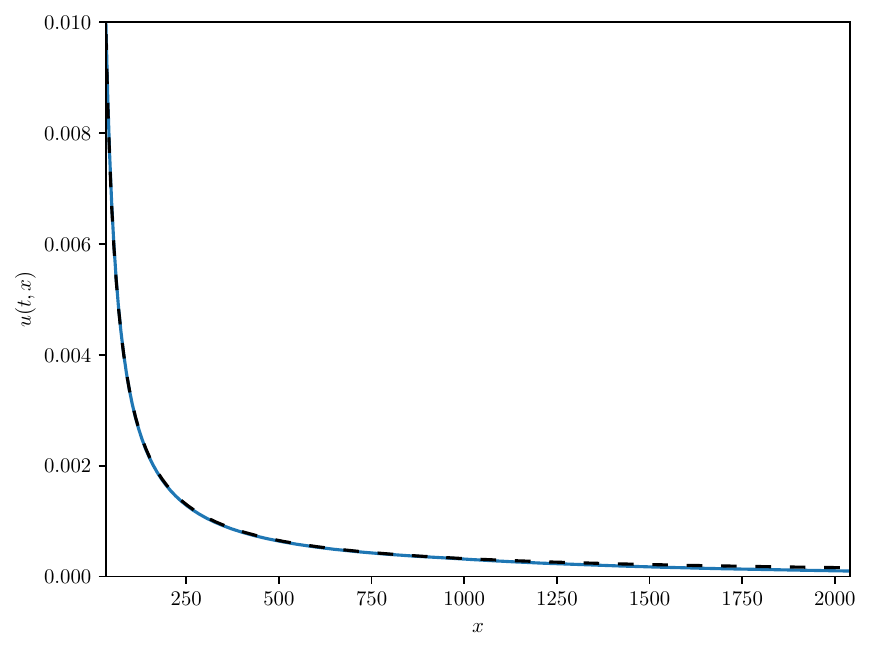}
		\caption{$s=0.5$.}
	\end{subfigure} %
	\begin{subfigure}[b]{0.24\textwidth}
    	\centering
    	\includegraphics[width=\textwidth]{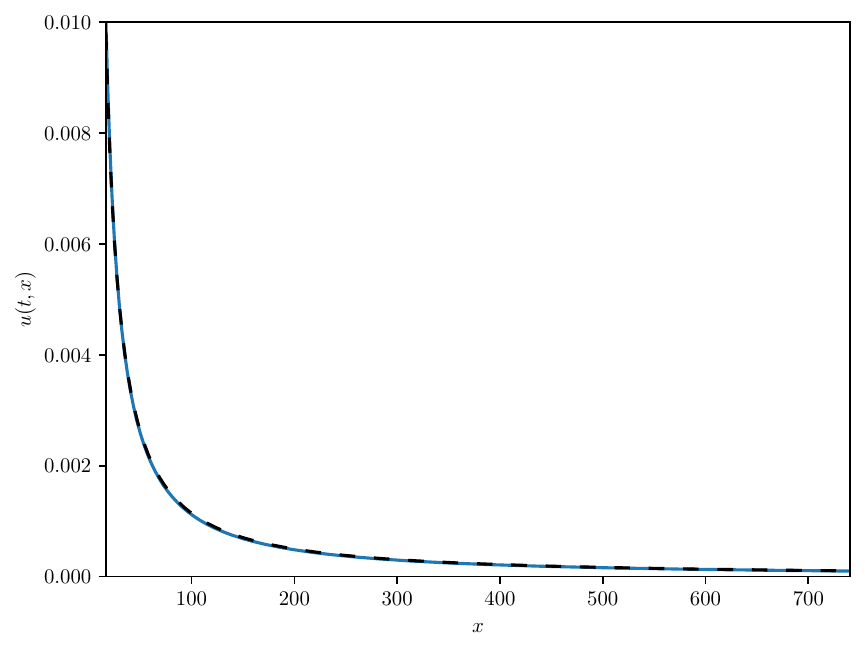}
		\caption{$s=0.6$.}
	\end{subfigure} %
	\begin{subfigure}[b]{0.24\textwidth}
        \includegraphics[width=\textwidth]{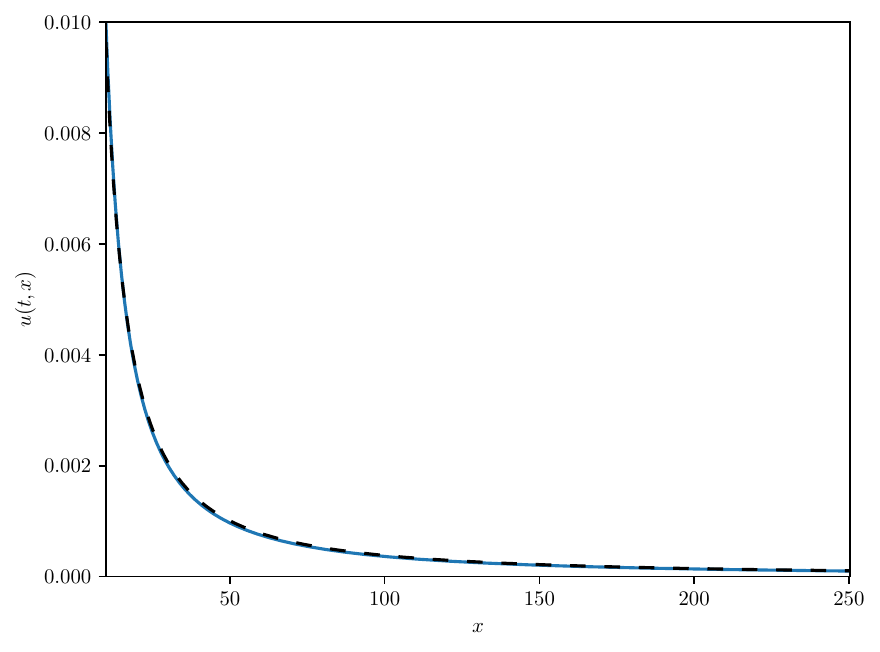}
    	\caption{$s=0.7$.}
    \end{subfigure} %

\caption{Fitting of of the function $\frac{D}{x^{2s}}$ with the part of the tail between values $10^{-2}$ and $10^{-4}$ of the numerical approximation at time $t=1$, for $\beta=1.5$ and different values of $s$. The dashed curves correspond to the fitted functions.}\label{nose-fit}
\end{figure}

Lastly, we tried to match the numerical approximation to the expected asymptotic profile of the solution at time $t=1$. To this end, we used the method of least-squares, implemented in the \texttt{scipy.optimize.curve\_fit} routine, to fit the constant in the function $\frac{D}{x^{2s}}$ to the part of the tail between values $10^{-2}$ and $10^{-4}$ of the computed numerical approximation. In Figure \ref{nose-fit}, it can be seen that, despite the imposed Dirichlet boundary conditions, the numerical approximation exhibits the correct decay for several values of the parameter $s$.

\bibliographystyle{abbrv}
\bibliography{cleaned_biblio}
\end{document}